\documentclass{ruthesis}

\usepackage{amsmath}
\usepackage{amscd}
\usepackage{latexsym}
\usepackage{amsfonts}
\usepackage{amsthm}
\usepackage{amssymb}
\usepackage{enumerate}
\usepackage{float}
\usepackage{paralist}
\usepackage[all]{xy}

\theoremstyle{definition}
\newtheorem{theorem}{Theorem}[section]
\newtheorem{lemma}[theorem]{Lemma}
\newtheorem{proposition}[theorem]{Proposition}
\newtheorem{corollary}[theorem]{Corollary}
\newtheorem*{myproof}{Proof}
\newtheorem{definition}[theorem]{Definition}
\newtheorem{example}[theorem]{Example}
\newtheorem{remark}[theorem]{Remark}

\newcommand{\Cart}{\operatorname{Cart}}
\newcommand{\Pic}{\operatorname{Pic}}
\newcommand{\Cl}{\operatorname{Cl}}
\newcommand{\Div}{\operatorname{Div}}
\renewcommand{\div}{\operatorname{div}}
\newcommand{\Hom}{\operatorname{Hom}}
\newcommand{\pHom}{\operatorname{Hom_{\mathbb{F}_1}}}
\newcommand{\Spec}{\operatorname{Spec}}
\newcommand{\MSpec}{\operatorname{MSpec}}
\newcommand{\MProj}{\operatorname{MProj}}
\newcommand{\Aut}{\operatorname{Aut}}
\newcommand{\Map}{\operatorname{Map}}
\newcommand{\Ho}{\operatorname{Ho}}
\newcommand{\Tor}{\operatorname{Tor}}
\newcommand{\Ext}{\operatorname{Ext}}

\newcommand{\Quot}{\operatorname{Quot}}

\newcommand{\coker}{\operatorname{coker}}
\renewcommand{\ker}{\operatorname{ker}}
\newcommand{\im}{\operatorname{im}}
\newcommand{\coeq}{\operatorname{coeq}}
\newcommand{\eq}{\operatorname{eq}}

\newcommand{\id}{\mathrm{id}}
\newcommand{\ann}{\mathrm{ann}}
\newcommand{\Ass}{\mathrm{Ass}}
\newcommand{\nil}{\mathrm{nil}}
\newcommand{\red}{\mathrm{red}}
\newcommand{\nor}{{\operatorname{nor}}}
\newcommand{\sn}{{\operatorname{sn}}}
\newcommand{\rank}{\mathrm{rk}}
\renewcommand{\dim}{\mathrm{dim}}

\newcommand{\ot}{\otimes}

\renewcommand{\P}{\mathbb{P}}
\newcommand{\A}{\mathbb{A}}
\newcommand{\F}{\mathbb{F}}
\newcommand{\N}{\mathbb{N}}

\newcommand{\R}{\mathbb{R}}
\newcommand{\Z}{\mathbb{Z}}
\newcommand{\bD}{\mathbf{D}}
\newcommand{\cA}{\mathcal{A}}

\newcommand{\cC}{\mathcal{C}}
\newcommand{\cD}{\mathcal{D}}

\newcommand{\cF}{\mathcal{F}}
\newcommand{\cL}{\mathcal{L}}
\newcommand{\cH}{\mathcal{H}}
\newcommand{\frakm}{\mathfrak{m}}
\newcommand{\frakp}{\mathfrak{p}}
\newcommand{\frakq}{\mathfrak{q}}

\newcommand{\D}{\Delta}
\renewcommand{\d}{\partial}
\newcommand{\Dop}{{\Delta^{op}}}
\newcommand{\snsphere}{\partial\Delta^n}

\newcommand{\nkhorn}{\Lambda^n_k}
\newcommand{\horn}[2]{\Lambda^{#1}_{#2}}

\newcommand{\set}{\text{\bf -sets}}
\renewcommand{\mod}{\text{\bf -mod}}

\newcommand{\Sets}{\mathbf{Sets}}
\newcommand{\mon}{\mathbf{Mon}_*}
\newcommand{\Ab}{\mathbf{Ab}}

\newcommand{\SSets}{\Dop\Sets}
\newcommand{\SAset}{\Dop A\set}

\newcommand{\Da}{\operatorname{Da}}
\newcommand{\Dag}{\operatorname{Da}_{\geq 0}}
\newcommand{\Dar}{\operatorname{\widetilde{D}a}}
\newcommand{\Darg}{\operatorname{\widetilde{D}a}_{\geq 0}}
\newcommand{\Ch}{\operatorname{Ch}}
\newcommand{\Chg}{\operatorname{Ch}_{\geq 0}}

\def\ra{\rightarrow}
\def\hra{\hookrightarrow}
\def\xra#1{\xrightarrow{#1}}
\newcommand{\da}[2]{\overset{#1}{\underset{#2}\rightrightarrows}}
\newcommand{\rra}{\rightrightarrows}
\def\lra{\longrightarrow}
\def\map#1{{\buildrel #1 \over \lra}}

\def\cof{\rightarrowtail}
\def\fib{\twoheadrightarrow}
\def\tcof{\overset{\sim}\rightarrowtail}
\def\tfib{\overset{\sim}\twoheadrightarrow}
\def\weq{\overset{\sim}\ra}

\begin{document} \phd
\setdefaultenum{i)}{a)}{i)}{i)}
\title{Homological Algebra for Commutative Monoids} 
\author{Jaret Flores}

\program{Mathematics} 
\director{Charles Weibel} \approvals{4}
\submissionyear{2015} 
\submissionmonth{January} 
\abstract{We first study commutative, pointed monoids providing basic definitions and results in a manner similar commutative ring theory.  Included are results on chain conditions, primary decomposition as well as  normalization for a special class of monoids which lead to a study monoid schemes, divisors, Picard groups and class groups.  It is shown that the normalization of a monoid need not be a monoid, but possibly a monoid scheme.  

After giving the definition of, and basic results for, $A$-sets, we classify projective $A$-sets and show they are completely determine by their rank.  Subsequently, for a monoid $A$, we compute $K_0$ and $K_1$ and prove the Devissage Theorem for $G_0$.  With the definition of short exact sequence for $A$-sets in hand, we describe the set $Ext(X,Y)$ of extensions for $A$-sets $X,Y$ and classify the set of square-zero extensions of a monoid $A$ by an $A$-set $X$ using the Hochschild cosimplicial set.

We also examine the projective model structure on simplicial $A$-sets showcasing the difficulties involved in computing homotopy groups as well as determining the derived category for a monoid.  The author defines the category $\Da(\cC)$ of double-arrow complexes for a class of non-abelian categories $\cC$ and, in the case of $A$-sets, shows an adjunction with the category of simplicial $A$-sets.  }
\beforepreface 
\acknowledgements{The author would like to thank Charles Weibel for years of advising including the collaboration leading to \cite{FW} which constitutes part of Chapter 2 and all of Chapter 4.  The author would also like to thank classmate Knight Fu and colleague Oliver Lorscheid for the many discussions on topics covered in this thesis.}
\dedication{This paper is dedicated to my friends, family and classmates whose support contributed to my successful navigation of graduate school and its many obstacles.  And also to my advisor whose generous patience allowed for the completion of this paper.}
\afterpreface

\chapter{Introduction}

As pointed out in \cite{deit}, Jacques Tits' question about the existence of $\F_1$, the ``field with one element,'' spawned interest in commutative, pointed monoids.  In this paper we consider $\F_1=\{0,1\}$ to be the multiplicative monoid having only an identity and zero element.  An $\F_1$ ``algebra'' $A=\F_1[S]$ consists of all monomials provided by a semigroup $S$ with ``coefficients'' in $\F_1$.  Then $A$ is a monoid with zero element and if we also impose the condition that $S$ be commutative (or abelian), then $A$ is a commutative, pointed monoid.  In this paper the term monoid will always mean commutative, pointed monoid unless otherwise stated.

The inclusion of a zero element in $\F_1$, and hence all monoids, provides a more interesting theory of ideals since it allows for the existence of zero divisors.  In any case, monoids have a rich theory of ideals including prime spectrums and their associated Zariski topology.  Perhaps the most noticeable divergence from commutative ring theory is the \emph{lack} of a Nakayama's Lemma due to the non-cancellative nature of monoids (a monoid $A$ is \emph{cancellative} when $ab=ac$ in $A$ means $b=c$).  Although this lack of cancellation does not provide significant obstacles to the commutative algebra of monoids, we will see that it leads to a more problematic homological theory.

Prime ideals lead the way to a theory of monoid schemes which is discussed briefly in the last chapter whose content is drawn from \cite{FW}.  There we provide a definition for the normalization of a special class of monoids.  It is then shown that the normalization of such monoids need not be a monoid, but rather a monoid scheme.  This is one more way in which the theory of monoids departs from the analogy with commutative rings and also shows the necessity for quickly moving into a study of the geometry of monoids.

When $k$ a field, the analog of a ``module over a $k$-algebra'' is obtained by forgetting the additive structure of those modules.  That is, a ``module'' over a monoid $A$ is a pointed set $X$ together with an action $A\times X\ra X$.  These objects, called $A$\emph{-sets}, are the primary study of this paper.  A morphism of $A$-sets is an $A$-equivariant, pointed set map and the category $A\set$ is both complete and cocomplete, though not abelian.

Perhaps the most significant feature of monoids, $A$-sets and their morphisms is that they do not admit a ``First Isomorphism Theorem.''  That is, when $f:X\ra Y$ is an $A$-set morphism, \emph{it is not true in general} that $X/\ker(f)\cong \im(f)$.  For example, if $A=\F_1$, then $\F_1\set=\Sets_*$ the category of pointed sets where it is well known that the coproduct of two pointed sets $X,Y$ is the wedge sum $X\vee Y=(X\coprod Y)/(0_X\sim 0_Y)$.  The ``fold'' map $\F_1\vee\F_1\ra \F_1$ defined by $1\mapsto 1$ in both summands is surjective with trivial kernel.  Spoken another way, \emph{we do not have} that $\ker(X\xra{f} Y)=0$ implies $f$ is injective (i.e., one-to-one).  It is for this reason that the definition of a \emph{congruence}, i.e. an equivalence relation that is also an $A$-set, plays a central role in the theory of monoids and $A$-sets.

Readers familiar with homological algebra for abelian categories are likely thinking that this failing of $A$-sets and their morphisms is bound to create many obstacles in a homological theory for $A$-sets.  Moreover, any homological theory for $A$-sets must remain valid when $A=\F_1$ and $A\set=\Sets_*$.  Hence we are drawn into a study of simplicial objects and homotopy groups which is well known to be computationally difficult.

It was the author's primary goal to investigate the possibility of a homological theory for $A\set$ that corresponds to the homotopy theory of simplicial $A$-sets in a manner analogous to the Dold-Kan Theorem for abelian categories.  When $\cA$ is an abelian category, the Dold-Kan Theorem defines an adjunction $K:\Chg(\cA)\leftrightarrows\Dop \cA:N$, the latter category being that of simplicial objects in $\cA$, that is an equivalence of categories.  Furthermore, this equivalence descends to their homotopy categories $\Ho \Chg(\cA)\cong\Ho\Dop\cA$ (a so-called Quillen equivalence).  Recall that the homotopy category $\Ho\cC$ of a model category $\cC$ is obtained from $\cC$ by localizing at the weak equivalences.  When $\cA=R\mod$, the category of modules over a commutative ring, of primary homological interest is the derived category $\mathbf{D}(R)$ of $R$ which is  the homotopy category $\Ho\Ch(R)$.  This begs the question: what should be the derived category of a monoid $A$?

Aside from being markedly simpler to work with, the category $\Ch(R)$ from which we obtain $\mathbf{D}(R)$ can be thought of as an extension of $\Dop R\mod\cong\Chg(\cA)$, and it is because of Dold-Kan that this extension is possible.
\[\SelectTips{eu}{12}\xymatrix{
\Dop R\mod \ar[d]\ar[r]^{\cong} & \Chg(R) \ar[d] \ar[r] & \Ch(R) \ar[d]\\
\Ho\Dop R\mod \ar[r]^{\ \quad\cong} & \mathbf{D}_{\geq0}(R) \ar[r] & \mathbf{D}(R)
}\]
We are unable to produce an analogous diagram for $A$-sets since there is no candidate for the category of ``complexes of $A$-sets'' bounded below by 0, say $\mathrm{Com}_{\geq0}(A)$.  Whatever $\mathrm{Com}_{\geq0}(A)$ might be, it is clear that it must have a model structure so that $\Ho\mathrm{Com}_{\geq0}(A)\cong \Ho\Dop A\set$.  The (projective) model structure on $\Dop R\mod$, hence $\Chg(R)$, and $\Dop A\set$ are both defined in terms of the model structure for the category $\Dop\Sets$ of simplicial sets (see \ref{d:ssetmodelstructure}).  In both cases, a map is a weak equivalence (resp., fibration, resp., cofibration) when it is so on the underlying simplicial set.  It is reasonable to believe that one could use this same strategy to obtain the desired model structure on $\mathrm{Com}_{\geq0}(A)$.  From this it is immediate that 
\begin{compactenum}
\item the homology of a complex of $A$-sets will have, at least, the structure of a group,
\item computing homology \emph{groups} in $\mathrm{Com}_{\geq0}(A)$ will be difficult since, contrary to $\Ch(R)$, \emph{not every} complex of $A$-sets will be fibrant.
\end{compactenum}
Looking at (ii) from a simplicial perspective, consider the following.  Since $R$-modules have an underlying abelian group structure, every simplicial $R$-module is fibrant.  Then the homotopy groups $\pi_*(M)$ of a simplicial $R$-module $M$ may be computed directly from $M$ and it is not surprising that these homotopy groups are themselves $R$-modules.  The Dold-Kan correpondence then allows us to work with chain complexes, rather than simplicial objects.

On the other hand, $A$-sets do not even possess an internal binary operation and in general, $A$ is itself less than a group.  Therefore, a general simplicial $A$-set is not fibrant and one is not be able to compute the homotopy groups $\pi_*(X)$ directly from a generic simplicial $A$-set $X$.  One must first find a fibrant replacement $\widetilde{X}$ for $X$ and find $\widetilde{X}$ is likely to be most difficult part of computing the $\pi_*(X)$.  Since $\mathrm{Com}(A)_{\geq 0}$ will have a model structure that mirrors $\Dop A\set$, one will be forced to compute fibrant replacements in $\mathrm{Com}(A)_{\geq 0}$ as well.  Of course, the necessity for fibrant objects stems from the fact that the invariants (homotopy classes of maps) provided by model categories are at least groups.  Perhaps the future will bring a theory of ``model categories'' where the invariants computed for objects in a category $\cC$ reside within $\cC$ itself.

In Chapter \ref{c:dacomplexes} we define the category $\Da(A)$ of \emph{double-arrow complexes of} $A$-sets which the author believed to fill the role of chain complexes for a special class of non-abelian categories.  Functors analogous to those used in the Dold-Kan correspondence are also defined in Chapter \ref{c:dacomplexes} (\ref{e:moore}, \ref{d:inversedk}) and are shown to give an adjunction $K:\Darg(A)\leftrightarrows\Dop A\set: N$.  However, they \emph{do not} provide an equivalence of categories.  The author has not determined if the model structure on $\SSets$ translates to double-arrow complexes, though it seems reasonable.  Therefore, it is unknown to the author if, after restricting to their respective categories of fibrant (or fibrant and cofibrant) objects, the functors $K,N$ provide a Quillen equivalence.  Whatever the case may be, the author hopes the homological definitions and strategy of this paper serve as an example of how, or how not, to determine the derived category of a monoid.

\chapter{The commutative algebra of monoids}

Here we investigate the algebraic properties of commutative, pointed monoids.
There are many similarities between these monoids and commutative rings, hence,
the outline of the theory follows suit.  We provide many basic definitions and
results which point out the difficulties which arise in objects having no
(abelian) group structure.  Included is the definition of short exact sequence for $A$-sets in Section \ref{exactsequences}.

There is much overlap between this chapter and \cite{Kobsa}.  A much more
general theory of ideals can be found in \cite{H-K}.

\section{Monoids and ideals}\label{monoidsandideals}

A \emph{monoid} is a set $A$ together with a binary operation $\cdot:A\times A\ra A$ satisfying $(a\cdot b)\cdot c=a\cdot (b\cdot c)$ (associativity) and an \emph{identity element} $1_{A}\in A$ satisfying $a\cdot 1_A=1_A\cdot a=a$ for all $a\in A$.  We will generally drop the $\cdot$ and $1_{A}$ notation and write the operation using juxtaposition $ab$ and the identity as $1$.  The monoid is \emph{commutative} if we also have $ab=ba$ for every $a,b\in A$.  A \emph{zero element}, or \emph{basepoint}, is a unique element $0\in A$ satisfying $a0=0a=0$ for every $a\in A$ and a monoid having a basepoint is \emph{pointed}.

The basepoint will always be written $0$ or $*$ except perhaps in specific examples. Many sets have multiple monoid structures and the role of ``0'' may change.  In cases when the monoid operation is not specified and the role of $0$ may be unclear, we will explicitly state identity and basepoint elements.  If necessary we denote a monoid by a ordered pair listing the set and operation.

\emph{Throughout this thesis, the term \textbf{monoid} will always mean a
commutative, pointed monoid unless otherwise stated}.

\begin{example}\label{e:pointedgroup}
If $G$ is any abelian group, we can form a monoid by adding a disjoint basepoint $G_{+}=G\coprod\{*\}$.  If $H_{n}=\{1,x,x^{2},\ldots,x^{n}=1\}$ is the cyclic group of order $n$, we let $C_{n}=H_{n}\cup\{0\}$ denote its associated pointed monoid.  Notice that $C_{1}$ is the pointed trivial group and as noted in the introduction, it is also the ``field with one element'' $\F_1=C_1$.
\end{example}

\begin{example}
After including a disjoint basepoint, the natural numbers $\N_+$ has additive monoid structure $\{-\infty,0,1,2,\ldots\}$ and is the \emph{free monoid in one variable} (see Section \ref{freeasets}).  Alternatively, we may consider the natural numbers $\{0,1,2,\ldots\}$ as a monoid using multiplication as the binary operation.  This latter multiplicative monoid structure for $\N$ is of little interest so, for convenience, the notation $\N$ will always refer to the free monoid in one variable unless otherwise stated.  The free monoid in one variable may also be written  multiplicatively as $\{0,1,x,x^2,\ldots\}$.  Since multiplicative notation is preferable, we take the \emph{monoid} $\N$ to be $\{0,1,x,x^2,\ldots\}$.
\end{example}

\begin{example}
When $X$ is any set, possibly infinite, we define the free monoid on generators the elements of $X$ to be
\[\{ 0,1,x_1^{k_1}\cdots x_n^{k_n}\ |\ x_i\in X\text{ for all } 1\leq i\leq n,\ n>0, k_i\geq 0\}.\]
The elements $x_1^{k_1}\cdots x_n^{k_n}$ are referred to as \emph{monomials} or \emph{words}.
\end{example}

Let $A$ be a monoid and $X\subseteq A$ a subset.  Then $A$ is \emph{generated} by $X$ if every non-zero $a\in A$ can be written $a=x_{1}^{m_{1}}\cdots x_{n}^{m_{n}}$ with $x_{i}\in X$ and $m_i\in\N$ for every $i$.  In this case we write $A=\langle x \ | \ x\in X\rangle$ and we do not require $X$ to contain $0, 1$ since their presence in $A$ is implied.  When $X$ can be chosen to be finite, $A$ is \emph{finitely generated}.  This definition for generators of a monoid coincides with that of groups.

\begin{example}\label{e:toricvariety}
$(\mathbb{R}^{n},+)$ is an unpointed monoid with identity element $0=(0,\ldots,0)$ containing $M=(\Z^{n},+)$ as an unpointed submonoid.  Given a finite set of vectors $v_1,\ldots,v_n$ of $M$, the $\mathbb{R}_{\geq 0}$-linear span $\sigma=\{\sum r_iv_i\ |\ r_i\in\mathbb{R}_{\geq0}\}$ is a \emph{rational polyhedral cone} when $\sigma$ contains no lines through the origin\cite{F}.  In this case, $S_{\sigma}=\sigma\cap M$ is a finitely generated, unpointed monoid, known as an \emph{affine semigroup}. Associated to $S_{\sigma}$ is the monoid ring $R=\mathbb{C}[S_{\sigma}]$ which is the coordinate ring of the affine toric variety $\Spec(R)$.
\end{example}

\begin{definition}
Let $A,B$ be monoids. Their \emph{product} is the usual cartesian product, i.e.
$A\times B=\{ (a,b) \ | \ a\in A,\ b\in B\ \}$ having identity $(1,1)$ and
basepoint $(0,0)$.  Their coproduct, also called the \emph{smash product},  is
$A\land B=(A\times B)/((A\times\{0\})\cup (\{0\}\times B))$ which is analogous
to the smash product of pointed topological spaces.  Elements of the smash
product are written $a\land b$ so that the identity element is $1\land 1$ and
the basepoint is $0=0\land 0$.
\end{definition}

\subsection{Morphisms of monoids}\label{monoidmorphisms}

Let $A,B$ be two monoids.  A morphism $f:A\ra B$ is function satisfying
\begin{center}
    \begin{inparaenum}[i)]
    \item $f(0)=0$\qquad
    \item $f(1)=1$\qquad
    \item $f(aa')=f(a)f(a')$
    \end{inparaenum}
\end{center}
Note that (i) and (ii) do not necessarily follow from (iii) as shown by the maps $A\ra A\times A$ defined by $a\mapsto (a,0)$ and $a\mapsto (a,1)$.  The image of a monoid $A$ under a monoid morphism $f$ is itself a monoid, called the \emph{image of $f$}.  Of course, the composition of two monoid morphisms is again a monoid morphism.  In general the monoid adjective will be dropped and we call a morphism between monoids simply, a \emph{morphism}.  The category of (commutative, pointed) monoids together with their morphisms will be denoted $\mon$.  The set of morphisms from $A$ to $B$ will be denoted $\Hom_{\mon}(A,B)$ or simply $\Hom(A,B)$ when there is no risk of confusion.

\begin{example}\label{e:forgetrings}
There is a functor $U:\mathbf{Rings}\ra\mon$ from commutative rings with
identity to monoids where $U(R)$ forgets the additive structure of $R$.  Then
$(U(R),\cdot)$ is a monoid with unit 1 and basepoint 0.
\end{example}

As usual, a monoid $B$ is a sub-monoid of $A$, written $B\subseteq A$, if there
exists an injection $i:B\ra A$.   In this case, $B\cong i(B)\subseteq A$.

\subsection{Ideals and quotients}\label{idealsandquotients}

When $A$ is a monoid, a subset $I\subseteq A$ is called an \emph{ideal} if
$ax\in I$ for any $a\in A$ and any $x\in I$.  The ideal is \emph{proper} if
$I\neq A$ and this occurs only when $1\not\in I$.  Let $I,J\subseteq A$ be
ideals and define the product of $I$ and $J$ to be $IJ=\{ab \ | \ a\in I,b\in J
\}$, which is again an ideal.  The intersection and union of a collection of
ideals is again an ideal; the intersection of finitely many ideals contains
their product.

An ideal $I\subseteq A$ is \emph{generated by} a subset $Y\subseteq I$ if every $x\in I$ can be written $x=ay$ for some $a\in A$ and $y\in Y$.  Here we write $I=A(y \ | \ y\in Y)$ or simply $I=(y \ | \ y\in Y)$ when there is no ambiguity. When $Y$ can be chosen finite the ideal is \emph{finitely generated}.  We will see in the next chapter that this definition is equivalent to $I$ being finitely generated as an $A$-set.  If $I$ is generated by a single element, say $x\in I$, then $I$ is a \emph{principal ideal} written using the above convention(s) or simply as $I=Ax$.

Let $f:A\ra B$ be a monoid homomorphism and $J\subseteq B$ an ideal.  Then the
inverse image $f^{-1}(J)$ is always an ideal, called the \emph{contraction} of
$J$.  On the other hand, if $I\subseteq A$ is an ideal, then the image $f(I)$
need not be.  Here we define the \emph{extension} of $I$ to be the ideal
generated by the elements of $f(I)$.

An equivalence relation, $R$, on $A$ is a subset $R\subseteq A\times A$
satisfying:
\begin{compactenum}
    \item $(a,a)\in R$ for all $a\in A$ (reflexive)
    \item $(a,b)\in R$ implies $(b,a)\in R$ (symmetric about the diagonal)
    \item $(a,b),(b,c)\in R$ implies $(a,c)\in R$ (transitive)
\end{compactenum}
When it is convenient we use the symbol $\sim$ to denote the equivalence
relation $R$ imposes on $A$ and denote an element $(a,b)\in R$ by $a\sim b$.
This notation emphasizes the role of equivalence relations play in producing
quotient monoids.

A \emph{congruence} on $A$ is an equivalence relation $R\subseteq A\times A$ such that $(x,y)\in R\times R$ implies $(ax,ay)\in R\times R$.  Given a congruence $R$ on $A$, define the equivalence class of $a\in A$ to be $[a]=\{b\in A \ | \ a\sim b\}$.  Multiplication of equivalences classes $[a][b]=[ab]$ is well defined since $a\sim c$ and $b\sim d$ means $ab\sim cb\sim cd$.  The set of equivalence classes of $A$ with respect to $R$, denoted $A/R$ or $A/\sim$, forms a monoid under the multiplication just defined.  There is a (monoid) morphism $\pi:A\ra A/\sim$ that sends every element to its equivalence class, i.e.  $\pi(a)=[a]$, called a \emph{quotient} map.  We will generally drop the $[a]$ notation and denote the equivalence classes in $A/\sim$ by $a$, or $\bar{a}$, with the understanding that the equivalence class is really what is meant.  This should not cause confusion as long as it is made explicit to which monoid the element is considered as being contained in at any given time.

Let $A$ be a monoid and $Y=\{(a_{1},b_{1}),\ldots\}\subseteq A\times A$ be \emph{any} subset.   The smallest congruence $X$ containing $Y$ is the
congruence \emph{generated} by $Y$ and the elements of $Y$ are the
\emph{generators of $X$}.   When $Y$ is finite, we say that $X$ is
\emph{finitely generated} and write the quotient $A/X$ as
$A/(a_1=b_1,\ldots,a_n=b_n)$.  We nearly always refer only to the generators of
a congruence, with the understanding that the congruence is actually what is
meant.

\begin{proposition}\label{p:monoidcong}
Every monoid morphism $f:A\ra B$ may be factored as the composition $A\xra{p}
f(A)\xra{i} B$ where $p$ is a surjection and $i$ is an inclusion.  Moreover, if
$R$ is the congruence on $A$ generated by relations $a\sim b$ when $p(a)=p(b)$,
then $A/R\cong f(A)$.  In particular, there is a one-to-one correspondence
between surjective morphisms $A\ra B$ and congruences on $A$.
\end{proposition}

When $I\subseteq A$ is an ideal, the subset $I\times\{0\}\subseteq A\times A$ generates a  congruence whose associated quotient monoid is written $A/I$ and identifies all elements of $I$ with $0$ leaving $A\backslash I$ untouched.  In this case $\pi:A\ra A/I$ has $\pi^{-1}(0)=I$ and $\pi^{-1}(a)=a$ for $0\neq a\in A/I$.  Whenever $A$ is an ideal the notation $A/I$ will always refer to the quotient of $A$ by the congruence generated by $I$.  When $I=(x_1,\ldots,x_n)$ is finitely generated, we may write the quotient monoid as $A/(x_1,\ldots,x_n)$ which is shorthand for $A/(x_1=0,\ldots,x_n=0)$.

\begin{remark}\label{r:congidealdiff}
This difference between a general congruence $R$ on a monoid $A$ and an ideal
$I\subseteq A$ is not seen in commutative ring theory.  Namely, the quotient
$A/R$ is not necessarily obtained from $A$ by identifying an ideal with 0.  The
notation for the quotient of a monoid by a congruence and ideal are written
similarly for convenience, but it is important to remember the difference.
\end{remark}

If $f:A\ra B$ is a morphism, then $f^{-1}(0)$ is an ideal of $A$ called the
kernel of $f$, written $\ker(f)$, and there is an induced morphism
$\bar{f}:A/\ker(f)\ra B$ defined by $[a]\mapsto f(a)$.  Since
$\bar{f}|_{A\backslash\ker(f)}=f$, we may write $f$ rather than $\bar{f}$ when
there will be no confusion.

We wish to stress it is not the case that a morphism $f:A\ra B$ is injective when the ideal $\ker(f)=0$.  For example, let $C_{n}$ be the pointed abelian group of order $n$ and consider $f:\N\ra C_{n}$ defined by $f(x)=x$.  Then $f(x^{k})=f(x^{k+n})$ for every $k\in\N$ but $\ker(f)=f^{-1}(0)=0$.  This example showcases the lack of a ``First Isomorphism Theorem'' for monoids, i.e. we \emph{do not have} $A/\ker(f)\cong\im(f)$.

\begin{remark}\label{r:addmultnotation}
We previously noted that $\langle x\rangle\cong\N$ and that we consider $\langle x\rangle$ a multiplicative notation for $(\N,+)$.  We use the notation $\N_{m}$ for the quotient monoid $\N/(m)\cong \langle x\rangle/(x^{m})$ having elements $\{0,1,x,\ldots,x^{m-1}\}$, not to be confused with $C_m\cong \langle x\rangle/(x^m\sim 1)$ the pointed finite \emph{group}.
\end{remark}

Given ideals $I,J\subseteq A$, one can form the \emph{ideal quotient} $(J:_{A}I)=\{a\in A \ | \ ab\in J \ \mathrm{for \ every} \ b\in I \}$ or $(J:I)$ if there is no chance for confusion.  If $I=(b)$ is principal we simply write $(J:b)$.  When $J=0$, the ideal quotient $(0:I)=\{a\in A \ | \ ax=0\ \mathrm{for \ all}\ x\in I\}$ is the \emph{annihilator} of $I$, denoted $\ann_{A}(I)$ or simply $\ann(I)$.  If $I=(b)$ and $\ann(b)\neq 0$, $b$ is a \emph{zero divisor} or \emph{torsion} element.  When $\ann(b)=0$ for every nonzero $b\in A$, $A$ is \emph{torsion free}.  A much stronger condition requires that $ab=ac$ implies $b=c$ for all $a,b,c\in A$ with $a\neq 0$.  Monoids satisfying this condition are \emph{cancellative} and consequently torsion free.  A monoid isomorphic to the quotient of a cancellative monoid by an ideal is \emph{partially cancellative} or \emph{pc}.

\subsection{Prime ideals and units}\label{primeidealsandunits}

An element $u\in A$ is a \emph{unit} if there exists $u^{-1}\in A$ with
$uu^{-1}=1$; equivalently $u$ is not contained in any proper ideal of $A$.  The
set of units, $A^{\times}$, of $A$ form an abelian group and is the largest
multiplicatively closed subset disjoint from every proper ideal.  The complement
of $A^{\times}$ is an ideal; it is the unique maximal ideal of $A$, and is written
$\frakm_A$ or simply $\frakm$.

An ideal $\frakp$ is \emph{prime} if $\frakp\neq A$ and it satisfies one of the following three equivalent conditions:
\begin{compactenum}
    \item $ab\in\frakp$ mean $a\in\frakp$ or $b\in\frakp$,
    \item $a\not\in\frakp$ and $b\not\in\frakp$ means $ab\not\in\frakp$,
    \item $A/\frakp$ is torsion free.
\end{compactenum}
The set of all prime ideals of a monoid $A$ is written $\MSpec(A)$.  The union
of prime ideals is prime, though the intersection and product are not.  For
convenience and tradition, we do not allow $A$ to be a prime ideal so that prime
ideals are always proper.  Condition (ii) implies that $A\backslash\frakp$ is a
multiplicatively closed subset of $A$ containing $A^{\times}$, namely
$A\backslash\frakp\cup\{0\}\subseteq A$ is a submonoid.  Every proper ideal is
disjoint from $A^{\times}$ so that the prime $\frakm=A\backslash A^{\times}$ is
the unique maximal ideal of $A$; then $A/\frakm\cong A^{\times}_{+}$ is a
pointed abelian group.  In commutative ring theory a ring is \emph{local} when
it has only a single maximal ideal.  Using this vocabulary we have that
\emph{every monoid is local}.  This is one way in which the structure of monoids
is simpler than that of rings.  However, with no Nakayama's Lemma (see Remark
\ref{r:noNAK}) available, monoids have their own complexities.

\begin{proposition}\label{p:finiteprimes}
If $A$ is a finitely generated monoid, then $A$ has finitely many prime ideals.
\end{proposition}

\begin{myproof}
Let $X=\{a_1,\ldots,a_n\}$ be a set of generators for $A$ and
$\frakp\subseteq A$ a nonzero prime ideal.  If $a\in\frakp$ is nonzero,
then $a=u\cdot\prod_{1\leq i\leq n}a_i^{k_i}$ with $k_i\geq 0$ and $u$ a unit.
By definition of primality, $\frakp$ must contain at least one of the $a_i$
having $k_i> 0$.  More generally, $X\cap\frakp$ generates $\frakp$ since every
element of $A$ can be written as a product of elements of $X$.
\qed
\end{myproof}

\begin{remark}
Proposition ~\ref{p:finiteprimes} provides an upper bound on the cardinality of
$\MSpec(A)$ when $A$ is a finitely generated monoid.  That is, if $A$ is
generated by $n$ elements, $\MSpec(A)$ has at most $2^n$ primes.  The upper
bound is attained when $(0)$ is prime and every generator of $A$ generates a
(principal) prime ideal.
\end{remark}

The intersection of every prime ideal of $A$ is an ideal called the
\emph{nilradical} and denoted $\nil(A)$.  Every element of $\nil(A)$ is
\emph{nilpotent}, i.e. $a^{n}=0$ for some $n\geq 1$; conversely every nilpotent
element is contained in the nilradical.  We say that $A$ is {\it reduced} if
whenever $a,b\in A$ satisfy $a^2=b^2$ and $a^3=b^3$ then $a=b$. This implies
that $A$ has no nilpotent elements, i.e., that $\nil(A)=\{ a\in A: a^n=0 \text{
for some } n\}$ vanishes.  When $A$ is a pc monoid (see after Remark
\ref{r:addmultnotation}), it is reduced if and only if $\nil(A)=0$; in this case
$A_{\red}=A/\nil(A)$ is a reduced monoid.  Note that the equivalence of these
two conditions does not hold in general, for example $A=\langle 0,1,x,y\ |\
x^2=y^2,\ x^3=y^3\rangle$ is not reduced yet $\nil(A)=0$.

A proper ideal $\frakq\subseteq A$ is {\it primary} when $xy\in\frakq$ implies
$x\in\frakq$ or $y^n\in\frakq$.  Alternatively, $\frakq$ is primary when every
zero-divisor $a\in A/\frakq$ is nilpotent.  The \emph{radical} of an ideal $I$
is $\sqrt{I}=\{a\in A\ |\ a^n\in I\}$; it is a prime ideal when $I$ is primary.
It is easy see the radical of a primary ideal $\frakq$ is prime, say $\frakp$,
and when convenient we say that $\frakq$ is $\frakp$-primary.

Here are a couple of basic results, whose proofs exactly mimic those of ring
theory, that we will require later:

\begin{lemma}\label{l:pcont.quotprim}
Let $\frakp$ be prime ideal in a monoid $A$.
\begin{compactenum}
\item[i)] If $I_1,\ldots,I_n$ are ideals such that $\cap_i I_i\subseteq\frakp$,
    then $I_i\subseteq\frakp$ for some $i$. If in addition $\frakp=\cap_i I_i$,
    then $\frakp=I_i$ for some $i$.
\item[ii)] Let $\frakq$ be a $\frakp$-primary ideal of $A$. If $a\in
    A\backslash\frakq$, then $(\frakq:a)$ is $\frakp$-primary.
\end{compactenum}
\end{lemma}

When $A$ is a monoid, $\MSpec(A)$ is a topological space using the Zariski
Topology.  This is not to be confused with the maximal ideal spectrum
$\frakm$-$\Spec(R)=\{\frakm\in\Spec(R)\ |\ \frakm\ \mathrm{maximal}\}$, where
$R$ is a commutative ring.  For any set $S\subseteq A$, the Zariski closed set
containing $S$ is $V(S)=\{\frakp\in\MSpec(A)\ | \ \frakp\supseteq S\}$.  The
proof that this forms a topology is exactly the same as for rings.

We call a strictly increasing sequence of prime ideals $\frakp_0\subset
\frakp_1\subset\cdots$ a \emph{chain}.  When the sequence is finite, say
$\frakp_n$ is the final prime, we say the chain has \emph{length} $n$.  The
\emph{(Krull) dimension} of a monoid $A$ is the supremum of the lengths of all
chains of prime ideals in $A$; this may be infinite.  Proposition
\ref{p:finiteprimes} shows that all finitely generated monoids have finite
dimension.  The \emph{height} or \emph{codimension} of a prime $\frakp$ is the supremum of the
lengths of all chains of prime ideals contained in $\frakp$; equivalently the dimension of the localization $A_{\frakp}$ of $A$ at $\frakp$ (see Section \ref{localization}).

%

\section{$A$-sets}\label{asets}

Let $A$ be a monoid and $X$ a pointed set, i.e. $X$ has distinguished basepoint
denoted $0_{X}$.  A \emph{left $A$-action} on $X$ is a binary operation
$\cdot:A\times X\ra X$ satisfying:
\begin{compactenum}
\item $1\cdot x=x$
\item $0_A\cdot x=0_{X}$ and $a\cdot 0_X=0_X$
\item $(ab)\cdot x=a\cdot(b\cdot x)$ for every $a,b\in A$, $x\in X$
    (associativity)
\end{compactenum}
A \emph{left $A$-set} is a pointed set $X$ together with a left $A$-action.  One
may define a \emph{right $A$-set} in the obvious way.  If $B$ is another monoid,
a \emph{two-sided $A,B$-set} is a pointed set $X$ that is \emph{both} a left
$A$-set \emph{and} a right $B$-set with actions satisfying $(ax)b=a(xb)$ for all
$x\in X$, $a\in A$ and $b\in B$.  When $A=B$, hence $X$ has both a left and
right $A$-action, $X$ is an (non-commutative) \emph{$A$-biset}.  The action of
an $A$-biset \emph{commutes} when $a\cdot x = x\cdot a$ for all $a\in A$ and
$x\in X$; then $A$-bisets with a commutative $A$-action are \emph{commutative}.
As usual, we will drop the $\cdot$ notation and denote the action of $A$ on $X$
by juxtaposition $ax$.  An \emph{$A$-set} is a commutative $A$-biset and these
objects are our primary concern.

The term $A$-set is fairly conventional, though the terms $A$-module, $A$-system
and $A$-polygon are also found in the literature.  We avoid the former since
most objects referred to as ``modules'' have the structure of an abelian group
and we want to avoid this confusion.  The latter two are used mainly in
semigroup theory and we find the term ``set'' more suitable (and expedient!).

\begin{example}
\begin{compactenum}
\item If $I$ is any ideal of $A$, then $I$ and $A/I$ are both $A$-sets.
\item Let $R$ be a commutative ring.  The forgetful functor
    $U:\mathbf{Rings}\ra\mon$ from Example~\ref{e:forgetrings} induces the
    forgetful functor $U:R\mod\ra U(R)\set$.  To every $R$-module $M$, the
    $U(R)$-set $U(M)$ has no addition and retains its $R$-action.
\end{compactenum}
\end{example}

\begin{definition}
Let $X,Y$ be $A$-sets.  The coproduct of $X$ and $Y$, called the \emph{wedge
product} is $X\vee Y=(X\coprod Y)/(0_X\sim 0_Y)$.  A nonzero element of $Z=X\vee
Y$ is simply an element $x\in X$ \emph{or} $y\in Y$, but not both. (There is no
nonzero ``addition'' here!)  The product of $A$-sets is the usual cartesian
product.
\end{definition}

\begin{example}\label{e:groupsets}
Let $G$ be a (abelian) group, $G_+$ the associated monoid and $X$ a $G_+$-set. The structure of $X$ is well understood.  Since $X\backslash\{0\}$ is a $G$-set, it is the disjoint union of orbits each of which is isomorphic to $G/H$ for some subgroup $H$.  Thus we have $X=\bigvee_{i=0}^n(G/H_i)_{+}$ as a $G_+$-set. 
\end{example}

Let $S\subseteq X$ a subset.  We say $X$ is \emph{generated by} $S$ if every
$x\in X$ can be written $x=az$ for some $a\in A$ and $z\in S$.  In this case we
may write $X=A(z \ | \ z\in S)$ or when $S=\{z\}$ contains only a single
element, $X=Az$.  When $S$ can be chosen finite, we say that $X$ is
\emph{finitely generated}.

\subsection{Morphisms of $A$-sets}\label{asetmorphisms}

Let $A$ be a monoid and $X,Y$ be $A$-sets.  A function $f:X\ra Y$ is an $A$-set
\emph{morphism}, or simply \emph{homomorphism}, when it satisfies:
\begin{compactenum}
\item $f(0_{X})=0_{Y}$ (based)
\item $f(ax)=af(x)$ for every $a\in A$, $x\in X$ ($A$-equivariant)
\end{compactenum}
The category of $A$-sets
together with their $A$-set morphisms will be denoted $A\set$ and the set of
morphisms from $X$ to $Y$ by $\Hom_{A}(X,Y)$ or simply $\Hom(X,Y)$ when there is
no risk of confusion.  Note that $\Hom(X,Y)$ is itself an $A$-set with basepoint
the trivial map, $x\mapsto 0$, and $A$-action $(af)(x)=f(ax)$.  This makes
$\Hom_{A}(-,-): A\set^{op}\times A\set\ra A\set$ a bi-functor; it will be
discussed further in Section \ref{s:functors}.  Two $A$-sets $X,Y$ are
\emph{isomorphic} when there is a bijective morphism $f:X\ra Y$ (equivalently,
$f$ has an inverse).

\begin{remark}\label{r:homfunctornotation}
Let $f:X\ra Y$ be an $A$-set morphism.  We use the notation
\[f_*=\Hom(Z,f):\Hom(Z,X)\ra\Hom(Z,Y)\]
so that $f_*(\alpha)=f\alpha$ when $\alpha\in\Hom(Z,X)$.  Likewise, we use the
notation
\[f^*=\Hom(f,Z):\Hom(X,Z)\ra\Hom(Y,Z)\]
so that $f^*(\beta)=\beta f$ for every $\beta\in\Hom(X,Z)$.  Even though this
notation conflicts with that of extension and contraction of scalars (see
Section \ref{tensorproduct}), this should cause no confusion as the intent
should be clear from the context.
\end{remark}

\begin{remark}\label{r:wedge.morphism.notation}
Let $f:X'\vee X''\ra Y$ be an $A$-set morphism.  It will be convenient to write
$f=f'\vee f''$ where $f'=f|_{X'}$ and $f''=f|_{X''}$ even though both morphisms
will never be ``used simultaneously.''  That is, an element of $X'\vee X''$ is
an element $x'\in X'$ \emph{or} $x''\in X''$ so that
\[f(x)=(f'\vee f'')(x) = \left\{
\begin{array}{rl}
f'(x) & \text{if } x \in X',\\
f''(x) & \text{if } x\in X''
\end{array} \right.\]
makes sense.  This notation works well with the usual meaning of $\vee$ as the
logical or operand since we will only ever use $f'$ or $f''$ (but not both, so
in our case it should be exclusive).

Conversely, given two $A$-set morphisms $f:X\ra Y$ and $f':X'\ra Y'$, we define
the map $f\vee f':X\vee X'\ra Y\vee Y'$ by $(f\vee f')|_X=f$ and $(f\vee
f')|_{X'}=f'$.
\end{remark}

\subsection{Free $A$-sets}\label{freeasets}

A pair of important adjoint functors are the \emph{forgetful}, or
\emph{underlying pointed set}, functor $U: A\set\ra\Sets_*$ and the \emph{free
functor} $F:\Sets_*\ra A\set$ defined by $F(X)=\bigvee_{0\neq x\in X}A$.  The
bijection $\Hom_A(FX,Y)\cong\Hom_{\Sets}(X,UY)$ is achieved in the usual way:
every $A$-set map $f:FX\ra Y$ is determined by its generators (the elements of
$X$).  Alternatively, a pointed simplicial set is an $\F_1$-set and the
inclusion $\iota:\F_1\ra A$ induces the adjunction
\[\F_1\set\underset{\iota^*}{\overset{\iota_*}{\rightleftarrows}} A\set\]
where $F=\iota_*$, $U=\iota^*$ are the extension and contraction of scalars
respectively (see Example \ref{e:tensorcotriple}).  It should be clear that both $\iota_*$ and $\iota^*$ are exact functors.

\begin{example}\label{e:cotriple.notation}

We will use a notation for found in group theory which emphasizes that $A$ is acting on the left.  Let $X$ be a pointed set and define $A[X]=\vee_{x\in X}A$. Then write $[x]$ for the element $1$ in the component indexed by $x$, and set $[0]=0$.  Thus every nonzero element of $A[X]$ has the form $a[x]$ for a unique $a\in A$ and $x\in X$.\vspace{2mm}\\
\textbf{(Warning)} To avoid cumbersome notation, when $X$ is an $A$-set, we define $A[X]$ to be the free $A$-set $FU(X)$, equivalently $A[U(X)]/A[\{0\}]$.  This avoids the notation $A[U(X)\backslash\{0\}]$.
\end{example}

\subsection{$A$-subsets and quotients}\label{asubsetsandquotients}

We say that $Y$ is an \emph{$A$-subset} of $X$, denoted $Y\subseteq X$, when $Y$
is a subset such that $ay\in Y$ for every $a\in A$, $y\in Y$.  In line with
monoids (see Section~\ref{idealsandquotients}), a \emph{congruence} is an
equivalence relation $R\subseteq X\times X$ that is also an $A$-subset, namely
$(x,x')\in R$ implies $(ax,ax')\in R$.  For an arbitrary congruence we use the
notation $X/R$ to denote the quotient $A$-set in which $x=x'$ when $(x,x')\in
R$.  In general we may write elements of $R$ as $x\sim x'$ to emphasize these
elements are identified in the quotient.  When $R\subseteq X\times X$ is any
subset, we call the smallest congruence $R'$ containing $R$ to be the
\emph{congruence generated by $R$} and the elements of $R$ the \emph{generators
of $R'$}.  When $R$ is finite, $R'$ is \emph{finitely generated}.

If $Y\subseteq X$ is an $A$-subset, then $Y\times\{0\}\subseteq X\times X$ generates a congruence whose quotient $A$-set, denoted $X/Y$,
identifies every element of $Y$ with $0$ and fixes $X\backslash Y$.  When
$Y\subseteq X$, the notation $X/Y$ will always refer to the quotient of $X$ by
the congruence generated by $Y$.  For general congruences we adopt the
conventions of monoids.  Namely, when the congruence $R\subseteq X\times X$ is
generated by $\{(x_1,x'_1),\ldots\}$, we may write the quotient $A$-set $X/R$ as
$X/(x_1=x'_1,\ldots)$.  When the set of generators $\{(x_i,x'_i)\}$ of $R$ can
be chosen to be finite, $R$ is \emph{finitely generated}.

\begin{proposition}
Every $A$-set morphism $f:X\ra Y$ can be factored as $X\xra{p} f(X)\xra{i} Y$
where $p$ is a surjection and $i$ is an injection. Moreover, if $R$ is the
congruence on $X$ generated by relations $x=x'$ when $p(x)=p(x')$, then
$X/R\cong f(X)$.  In particular, there is a one-to-one correspondence between
surjective morphisms $X\ra Y$ and congruences on $X$.
\end{proposition}

As with monoids, we can also define the \emph{smash product} of $A$-sets $X$ and
$Y$ as the quotient
\[X\land Y=(X\times Y)/((X\times\{0\})\cup(\{0\}\times Y))\]
which is the $A$-set consisting of $0$ and all ordered pairs $(x,y)$ with
$x,y\neq 0$.  We write the elements of $X\land Y$ as $x\land y$.  The $A$-action
of $X\land Y$ is inherited from from $X\times Y$ so $a(x\land y)=ax\land ay$.

\begin{remark}
As with Remark \ref{r:congidealdiff}, not every congruence on an $A$-set $X$ is
realized as an $A$-subset.  Every $A$-subset generates a congruence but not
every congruence corresponds to an $A$-subset of $X$.  The notation for
$A$-subsets and congruences on $X$ are similar for convenience, but it is
important to remember this difference.
\end{remark}

For a morphism $f:X\ra Y$,  the kernel and cokernel are defined by the usual
categorical notions.  The kernel, denoted $\ker(f)$, is the pullback in the
diagram on the left and the cokernel, denoted $\coker(f)$ is the pushout of the
diagram on the right.
\[\SelectTips{eu}{12}\xymatrix{
0\times_Y X \ar[d] \ar[r] & X \ar[d]^f \\
0 \ar[r] & Y }
\qquad
\SelectTips{eu}{12}\xymatrix{
X \ar[d] \ar[r]^f & Y \ar[d] \\
0 \ar[r] & 0\vee_X Y
} \]
All kernels and cokernels exist in $A\set$ but \emph{we do not} have that $f$ is
injective when $\ker(f)=0$, and \emph{we do not} have in general the isomorphism
$X/\ker(f)\cong Y$.

\begin{example}\label{nokernel}
Let $A$ be any monoid and consider the $A$-set morphism $A\vee A\ra A$ which
restricts to the identity map on each summand.  This map is obviously surjective
and the kernel is $0$, however the fiber of any element contains two points.
This simple example of the ``misbehavior'' of $A$-set morphisms showcases the
primary obstacle in constructing a homological theory.
\end{example}

\begin{remark}
In $A\set$ every monomorphism is the kernel of its cokernel, but \emph{not
every} epimorphism is the cokernel of its kernel or any other morphism.  Due to
the latter fact we can no longer conclude that a morphism is a monomorphism when
its kernel is trivial.  However, we still have the fact that the kernel of a
monomorphism is trivial.

\emph{Throughout the remainder of this thesis the terms injective, monomorphism and one-to-one will all mean for a morphism $f:X\ra Y$, that $f(x)=f(x')$ implies $x=x'$}.
\end{remark}

An interesting property of $A$-set morphisms, similar to that of continuous
functions and connected components of topological spaces, is the following.

\begin{proposition}\label{p:splitting}
Let $f:X\ra \vee_{i\in I}Y_{i}$ be an $A$-set morphism and $X_i=f^{-1}(Y_{i})$.
Then $X=\cup_{i\in I}X_i$ and $X_i\cap X_j\subseteq\ker(f)$ for all $i,j$.  In
particular, $\ker(f)=0$ means $X=\vee_{i\in I}X_i$.
\end{proposition}

\begin{myproof}
If $f(x)\in Y_{i}$ for some $i\in I$, then $f(ax)=af(x)\in Y_{i}$ since $Y_{i}$
is an $A$-set.  Hence, each $X_i\subseteq X$ is itself an $A$-set and $f(X_i\cap
X_j)\subseteq Y_{i}\cap Y_{j}=0$ when $i\neq j$.

Of course when $x\in\ker(f)$, we have $x\in f^{-1}(0)\subseteq X_i$ for all $i$
hence $x\in X_i\cap X_j$ for all $i,j$.
\qed
\end{myproof}

Let $X$ be an $A$-set and $X'\subseteq X$ any $A$-subset.  Define the
\emph{annihilator} of $X'$, denoted $\ann_{A}(X')$ or simply $\ann(X')$, to be
the ideal $(0:_{A}X')=\{a\in A \ | \ ax=0 \ \mathrm{for \ every} \ x\in X'\}$. A
general $A$-set $X$ is \emph{faithful} when $\ann(X)=0$.  When $\ann(X)\neq 0$,
then $X$ is a faithful $A/\ann(X)$-set.  For an ideal $I\subseteq A$ and an
$A$-set $X$, the set $IX=\{ax \ | \ a\in I, x\in X\}\subseteq X$ is an
$A$-subset and we can form the quotient $X/IX$.  In this case
$I\subseteq\ann(X/IX)$ so that $X/IX$ is an $A/I$-set.

Let $X$ be a noetherian $A$-set.  The set of ideals $\{\ann(x) \ | \ x\in X\}$
forms a partially ordered set under containment and the maximal elements, whose
existence is given by the ACC (see Section \ref{s:chainconditions}), are prime.  To see that maximal elements are
prime, suppose $ab\in\ann(x)$ and $b\not\in\ann(x)$.  Then $a\in\ann(bx)$ and
since $\ann(x)\subseteq\ann(bx)$, if $\ann(x)$ were maximal, we must have
$\ann(bx)=\ann(x)$ hence, $a\in\ann(x)$.  Primes occurring in this way are
called the \emph{associated primes of} $X$ and the collection of all such primes
is denoted $\Ass_{A}(X)$ or simply $\Ass(A)$ when there is no chance for
confusion.  That is, a prime ideal $\frakp\subseteq A$ is an associated prime of
$X$ when $\frakp=(0:x)$ for some $x\in X$.  Then $ax=0$ in $X$ if and only if
$a\in (0:x)\subseteq \frakp$ for some $\frakp\in\Ass(X)$ so that
$D=\cup_{\frakp\in\Ass(X)}\frakp=\{$zero divisors of $X\}$.  When $X'\subseteq
X$ is an $A$-subset,  we say the primes $\Ass(X/X')$ are \emph{associated} or
\emph{belong} to $X'$.  The usual results regarding associated primes in ring
theory carry through with monoids.  This will be covered in more detail in
Section \ref{primarydecomposition}.

\begin{remark}\label{r:noNAK}
There is no Nakayama's Lemma for general monoids.  Since all monoids are local,
such a statement would read: If $I\subseteq A$ is a proper ideal and $X$ an
$A$-set, then $IX=X$ means $X=0$.  This is far from the truth.  For example,
consider the quotient $A=\langle x\rangle/(x^n=x)$ of the free monoid in one
variable. Here $\frakm=\{0,x,\ldots,x^{n-1}\}$ and $\frakm^k\cdot\frakm=\frakm$
for every $k\geq 0$.
\end{remark}

\subsection{Admissible exact sequences}\label{exactsequences}

In previous sections we have seen a big difference between the morphisms in
$A\set$ and those of abelian categories.  We expand upon this further now.  The
power of a homological theory stems from the ability to produce properties of
morphisms from the (non-)existence of kernels.  To remedy the situation, we
simply restrict our attention to morphisms $f$ which satisfy: $\ker(f)=0$ if and
only if $f$ is one-to-one.

We will see in Chapter \ref{c:homological} that considering \emph{only} these
morphisms will not be sufficient.  In particular, using such morphisms will not
allow every $A$-set to have a projective resolution.  For now, we recall that a
morphism $f:X\ra Y$ is \emph{one-to-one, injective} or a \emph{monomorphism}
when $f(x)=f(x')$ implies $x=x'$ for all $x,x'\in X$.  Also, $\ker(f)=0$ does
\emph{not} imply $f$ is injective.

\begin{definition}
A morphism $f:X\ra Y$ is \emph{admissible} when the surjection $f:X\ra f(X)$ is
a cokernel.  In this case,  $\ker(f)=0$ implies $f$ is injective.
\end{definition}

All injections $X\hra Y$ are clearly admissible since they are the cokernel of
the zero map $0\ra X$.  For a general morphism $f:X\ra Y$, whenever
$f(x)=f(x')\neq 0$ for some $x,x'\in X$, $X/\ker(f)\not\cong f(X)$ since
$\overline{x}\neq \overline{x}'$ in $X/\ker(f)$; hence $f$ cannot be admissible.
We can then say the following:
\begin{compactenum}
    \item $f:X\ra Y$ is admissible if and only if $f|_{X\backslash\ker(f)}$ is a
        (set-theoretic) injection.
    \item An admissible morphism is an injection if and only if it has trivial
        kernel.
    \item Admissible morphisms have a ``First Isomorphism Theorem,'' namely,
        $X/\ker(f)\cong\im(f)$.
\end{compactenum}

A sequence
\[\cdots\ra X_{n+1}\xra{f_{n+1}} X_n\xra{f_n} X_{n-1}\ra\cdots\]
a morphisms is \emph{admissible} when every morphism in the sequence is
admissible.  The admissible sequence is \emph{exact} when
$\im(f_{i+1})=\ker(f_i)$ for all $i$.  An \emph{admissible short exact
sequence}, or a.s.e.s., is an admissible exact sequence of the form
\[0\ra X'\xra{g}X\xra{f}X''\ra0.\]
In general we will refer to an admissible exact sequences simply as an a.e.s.
and reserve a.s.e.s. for situations when we wish to stress the a.e.s. is short.

An a.e.s. $0\ra X'\xra{g}X\xra{f}X''\ra0$ is called an \emph{extension of $X''$ by $X'$}.  We define and study the isomorphism classes of extensions of $X''$ by $X'$ in Section \ref{extensions}, but we provide a little intuition to the structure of such sequences now.

With notation as above, admissible short exact sequences have the following  familiar properties:
\begin{compactenum}
    \item $g$ is one-to-one since $\ker(g)=\im(0\ra X')$.
    \item $\im(g)=\ker(f)$
    \item $f$ is surjective since $\im(f)=\ker(X''\ra 0)$
\end{compactenum}
We now look more closely at the structure of short exact sequences.

\begin{lemma}\label{l:splitting} (Splitting)
The a.e.s. of $A$-sets $0\ra X\xra{g}Y\xra{f}Z\ra0$ splits, i.e. $Y\cong X\vee
Z$, if and only if:
\begin{compactenum}
\item There is an $A$-set morphism $\sigma:Z\ra Y$ with $f\sigma=$id$_{Z}$ or
\item There is an \emph{admissible} $A$-set morphism $\psi:Y\ra X$ with $\psi
    g=$id$_{Y}$.
\end{compactenum}
We refer to the morphisms $\varphi$ and $\psi$ as \emph{splitting} maps.
\end{lemma}

\begin{myproof}
We first note that when the sequence splits, the existence of $\sigma$ (resp.
$\psi$) defined in (i) (resp. (ii)) is obvious.

Conversely, suppose that (i) holds.  Let $\sigma:Z\ra Y$ be an $A$-set morphism satisfying $f\sigma=\id_Z$.  Since $f$ is admissible, $y\not\in\ker(f)$ implies $y\in\im(\sigma)$.  Hence, $Y=\ker(f)\cup\im(\sigma)$.  Moreover, $\ker(f)$ and $\im(\sigma)$ are $A$-subsets of $Y$ and $\ker(f)\cap\im(\sigma)=0$ since $\sigma$ is a section of $f$.  Thus, $Y=\ker(f)\vee\im(\sigma)\cong X\vee Z$.

Now, assume (ii) and write $W$ for $\ker(\psi)$.  Then $0\ra W\xra{i} 
Y\xra{\psi}X\ra 0$ is an a.e.s. with $\psi i=\id_W$.  By (i), $Y\cong X\vee W$ 
and hence the bijection $W\ra Z$ is an isomorphism of $A$-sets.
\qed
\end{myproof}

\begin{example}\label{e:nonsplitting}
Let $A=\langle x\rangle$ be the free monoid in one variable, $x$.  Fix $n>1$ and
consider the $A$-set $Y=(A\vee A')/(x^n=x'^n)$ where the $'$ notation is used to
distinguish the summands of $Y$.   There is an a.e.s. $0\ra A\xra{i}Y\ra
A/x^nA\ra 0$, where $i(1)=1'$, and a function $\psi:Y\ra A$ defined by
$\psi(1)=\psi(1')=1$ such that $\psi\circ i = \id$.  However, $Y\not\cong A\vee
A/(x^n)$ so that $\psi$ is \emph{not} a splitting map.  This shows the
admissible condition in Lemma \ref{l:splitting}(ii) is necessary.
\end{example}

In any a.e.s. $0\ra X'\ra X\ra X''\ra 0$, there is a unique, pointed set map $\sigma:X''\ra X$ giving a decomposition $X\cong X'\vee X''$ as \emph{pointed sets}.  The $A$-set structure of the $A$-subset $X'$ is completely determined. The isomorphism $X''\cong X/X'$ determines the action of $A$ on all \emph{non-zero divisors} in $\sigma(X'')$.  When $ax''=0$ in $X''$, the isomorphism only enforces that $a\sigma(x'')$ is in the $X'$ summand (as pointed sets) of $X$.  Therefore, to a define an $A$-set structure on $X$ leaving $X'$, $X''$ fixed, the only freedom available is to assign to every pair $(a,x'')$ with $ax''=0$ in $X''$, an element $x'\in X'$ in such a way that we obtain a valid $A$-set structure on $X$.

\section{Tensor products}\label{tensorproduct}

Though $A$-sets lack an abelian group structure, the tensor product of $A$-sets
is still an important categorical construction.  Its structure is much simpler
than its counterpart in commutative ring theory.  Also, its role as a universal
object is decreased in the absence of bilinear mappings, however, by removing
the additive relations of bilinear mappings, we may still define the tensor
product in this way.

Let $X,Y,Z$ be $A$-sets and $f:X\times Y\ra Z$ a function.  We say that $f$ is a \emph{bi-equivariant} $A$-set map when $f(ax,y)=af(x,y)\quad\text{and}\quad f(x,ay)=af(x,y)$ for all $a\in A$ and $(x,y)\in X\times Y$.  That is, $f$ is $A$-equivariant in both coordinates.

\begin{definition}
The \emph{tensor product} of $X$ and $Y$ is an $A$-set $T$ satisfying the following universal property for bi-equivariant maps $f$:
\[\SelectTips{eu}{12}\xymatrix{
X\times Y \ar[dr]_f \ar[r] & T \ar @{-->} [d] ^{\exists g} \\
 & Z
} \]
As in ring theory, we may construct the tensor product in the following way.
Consider the free $A$-set (see Section \ref{freeasets}) $A[X\times Y]$ having
generators the \emph{nonzero} elements $(x,y)$ of $X\times Y$.  Let $R$ be the
congruence on $A[X\times Y]$ generated by all relations of the form
\[[ax,y]\sim a[x,y]\quad\text{and}\quad [x,ay]\sim a[x,y].\]
Then $T=A[X\times Y]/R$.  We write $X\otimes_A Y$ for the tensor product $T$, or simply $X\otimes Y$ when there is no risk of confusion, and $x\ot y$ for its elements.
\end{definition}

\begin{remark}
Elements of the form $x\otimes 0$ and $0\otimes y$ are equivalent to 0 in
$X\otimes_A Y$.  In the construction for $X\otimes_A Y$ we can replace
$A[X\times Y]$ with $A[X\land Y]$ where $X\land Y$ is the smash product (see
Section \ref{asubsetsandquotients}).  It is \emph{not true} in general
that $X\otimes_A Y$ is the $A$-set $X\land Y$ modulo relations of the form
$a(x,y)\sim(ax,y)$ and $a(x,y)\sim(x,ay)$.  The $A$-action on $X\land Y$ is
coordinate-wise, since it is a quotient of $X\times Y$, and the relations would
imply $a(x,y)=(ax,ay)=a^2(x,y)$ for every $a\in A$.  In other words, $X\land Y$
modulo the aforementioned relations is isomorphic to $X\ot_A Y$ only when every
element of $A$ is idempotent.
\end{remark}

\begin{proposition}\label{p:tensorproperties}
Let $X,Y,Z$ be $A$-sets.  Then we have the following:
\begin{compactenum}
\item $X\otimes Y\cong Y\otimes X$
\item $(X\otimes Y)\otimes Z\cong X\otimes(Y\otimes Z)$
\item $(X\vee Y)\otimes Z\cong (X\otimes Z)\vee(Y\otimes Z)$
\item $(X/Y)\ot Z=(X\ot Z)/(Y\ot Z)$
\end{compactenum}
Moreover, $-\ot_A Y$ is a functor and a left adjoint via the adjunction
$\Hom_A(X\ot_A Y,Z)\cong\Hom_A(X,\Hom_A(Y,Z))$.  Hence $-\ot_A Y$ preserves
colimits (see \cite{WH}) from which (iv) is a special case.
\end{proposition}

If $f:A\ra B$ is a monoid morphism and $X$ an $B$-set, then $X$ is naturally an
$A$-set with $A$-action $ax:=f(a)x$.  This gives a functor $f^{*}:B\set\ra
A\set$, called the \emph{restriction of scalars from $B$ to $A$}.  In
particular, $B$ is an $A$-set and when $Y$ is any other $A$-set, $B\otimes_{A}Y$
is a $B$-set with action $b(b'\otimes y)=bb'\otimes y$.  This provides a functor
$f_{*}: A\set\ra B\set$ called the \emph{extension of scalars from $A$ to $B$}
with  $f_{*}Y=B\otimes_{A}Y$.

\begin{example}\label{e:tensorcotriple}
Any pointed set is an $\F_1$-set.  If $X,Y$ are $\F_1$-sets, then
$X\otimes_{\F_1}Y\cong X\land Y$ is their smash product.  When $A$ is a monoid,
there is always a morphism $\iota:\F_1\ra A$ and when $X$ is an $A$-set,
$\iota^*X$ is its underlying pointed set (removing the $A$-action).  Note that
$\iota^{*}\iota_{*}X=A\otimes_{\F_1}X=A[X]$ is a free $A$-set which is
isomorphic to $A\land X$ as \emph{pointed sets}, but not as $A$-sets due to the
action $a(b\land x)=ab\land ax$ on $A\land X$.
\end{example}

\subsection{Tensor product of monoids}\label{tensoralgebras}

Let $C$ be a monoid.  Define a \emph{$C$-monoid} or \emph{monoid over $C$} to be a monoid morphism $f:C\ra A$.  Given another $C$-monoid, $g:C\ra B$, one can form the \emph{smash product, or tensor product, of $A$ and $B$ over $C$} as the pushout of the following diagram:
\[\SelectTips{eu}{12}\xymatrix{
C \ar[d]_g \ar[r]^f & A \ar[d]_{\lrcorner} \\
B \ar[r] & A\land_C B
} \]
namely, $A\land_{C}B$ is the coproduct (see Section \ref{monoidsandideals})
$A\land B$ modulo the relations $(af(c),b)\sim(a,g(c)b)$.  Using the usual
tensor product notation, this is written $A\otimes_C B$.  Since the smash
product of monoids $A,B$ is the tensor product as $\F_1$-monoids, it is often
denoted $A\otimes_{\F_1}B$ or simply $A\otimes B$.  In an algebraic setting the
$\otimes$ notation is standard and we follow this convention.  We will only use
the smash product notation $\land$ in situations where we wish to emphasize its
coproduct/topological nature and do not require the language of tensor products. When no clear $C$-monoid structure on two monoids $A,B$ is present, the notation $A\otimes B$ will always mean $A\otimes_{\F_1}B=A\land B$.

\begin{remark}\label{r:tensornotationwarning}
This notation for the tensor product of monoids differs from the convention for $A$-sets.  When $A,B$ are monoids and $B$ is also an $A$-set, the notation $A\otimes B$ can be ambiguous.  If the tensor product is as monoids, then $A\ot B=A\ot_{\F_1} B=A\land B$, but as $A$-sets, $A\otimes B=A\otimes_{A}B=B\neq A\land B$.  It should be clear from the context whether the tensor product is formed as ($\F_1$-)monoids or $A$-sets; we will be explicit when there is a possibility for confusion.
\end{remark}

\begin{example}\label{e:algebranotationwarning} (Notation Warning!)
Any monoid $A$ is an $\F_1$-monoid and may be written as the tensor product
$\F_1\ot A$.  For the free monoid on one variable $\F_1\ot\langle
x\rangle=\{0,1,x,x^2,\ldots\}$ we adopt the algebra notation of commutative ring theory: $\F_1[x]$.  In general, we write $\F_1[x_1,\ldots,x_n]=\F_1[x]^{\otimes n}$ for the free monoid on $n$ variables.  If $A$ is any monoid, $A[x_{1},\ldots,x_{n}]=A\otimes \F_1[x_1,\ldots,x_n]$ is the \emph{free $A$-monoid on $n$ generators} and its elements are \emph{monomials with coefficients in $A$}.

Let $m_1,\ldots,m_{2k}$, $k>0$, be monomials of $A[x_1,\ldots,x_n]$ and $R$ the
congruence generated by the relations $m_1\sim m_2,\ldots,m_{2k-1}\sim m_{2k}$.
When convenient, we may write $A[x_1,\ldots,x_n\ |\ m_1=m_2,\ldots,m_{2k-1}=
m_{2k}]$ or $A[x_1,\ldots,x_n\ |\ m_1\sim m_2,\ldots,m_{2k-1}\sim m_{2k}]$ for
the quotient $A[x_1,\ldots,x_n]/R$.

In the absence of context, the notation $A[x_1,\ldots,x_n]$ is ambiguous since it may refer to either the free $A$-monoid on $n$ variables or the free $A$-set with generators $x_1,\ldots,x_n$ (see Sections \ref{freeasets}, \ref{s:rank}). Therefore, when using this notation, the author will always be explicit about the nature of the object as an $A$-monoid or free $A$-set.
\end{example}

\begin{proposition}\label{p:tensorideals}
Let $A, B$ be monoids.  Every ideal $K\subseteq A\otimes B$ can be written as $K=\bigcup_{\lambda\in\Lambda}I_{\lambda}\ot J_{\lambda}$ where $I_{\lambda}\subseteq A$, $J_{\lambda}\subseteq B$ are (not necessarily proper) ideals and $\Lambda$ is an indexing set.  Moreover, at least one of the sets of ideals $\{I_{\lambda}\}$ and $\{J_{\lambda}\}$ can be chosen to contain distinct elements.
\end{proposition}

\begin{myproof}
First note that $K$ is the union of the ideals $Ax\ot By$, $x\ot y\in K$, since if $x\ot y\in K$, then $Ax\ot By\subseteq K$.  This proves the first assertion. 

Next, note that if $I,I'\subseteq A$ and $J\subseteq B$ are ideals, then $(I\ot J)\cup(I'\ot J)=(I\cup I')\ot J$ (likewise, if $I\subseteq A$ and $J,J'\subseteq B$ are ideals, then $(I\ot J)\cup(I\ot J')=I\ot(J\cup J')$). Now for each $y\in B$, define $I_y=\bigcup_{x\ot y\in K}Ax$ and for fixed $I_{\lambda}=I_y$, set $J_{\lambda}=\bigcup\{By\ |\ I_y=I_{\lambda}\}$.  Evidently, $K=\bigcup_{\lambda\in\Lambda}I_{\lambda}\ot J_{\lambda}$ and the set $\{I_{\lambda}\}$ of ideals contains distinct elements.  Moreover, it is clear from the construction how we could have instead made the $J_{\lambda}$ distinct, i.e. start by defining $J_x=\bigcup_{x\ot y\in K}By$.
\qed

\end{myproof}

Let $f:A\ra B$ and $g:A\ra C$ be monoid morphisms.  An \emph{$A$-monoid
morphism} $h:B\ra C$ is a monoid morphism that is also an $A$-set map, where
$B,C$ are considered $A$-sets by restriction of scalars (see Section
\ref{tensorproduct}).  We say the monoid morphism $f$ is \emph{finite} and that
$B$ is a \emph{finite $A$-monoid} when $B$ is finitely generated as an $A$-set.
The morphism is of \emph{finite type} and $B$ is a \emph{finitely generated
$A$-monoid} when there exists elements $b_1,\ldots,b_n\in B$ such that every
element of $B$ can be written as a \emph{monomial} in the $b_i$ and elements of
$f(A)$.  In this case we say that the elements of $B$ are monomials in the $b_i$
\emph{with coefficients in $f(A)$}.  Notice that a monoid is finitely generated
(see Section \ref{monoidsandideals}) when it is finitely generated as an
$\F_1$-monoid.

\begin{proposition}\label{p:monoidstructure}
Let $A$ be a monoid.
\begin{compactenum}
\item  If $B$ is a countably generated $A$-monoid, then $B\cong
    A[x_1,x_2,\ldots]/R$ for some congruence $R$ on $A[x_1,x_2,\ldots]$.
\item $A\cong G_+[x_1,\ldots]/R$ where $G_+=A/\frakm$ is the
    pointed group of units of $A$.  Moreover, when the submonoid 
    $\frakm\cup\{1\}\subseteq A$ is finitely generated, $A$ is a finitely generated $A/\frakm$-monoid.
\end{compactenum}
\end{proposition}

\begin{myproof}
i) Let $b_1,b_2,...$ be generators for $B$ not in $f(A)$.  Then every element of
$B$ can be written as a product of elements of $f(A)$ and the $b_i$. Hence,
there is a surjection $A[x_1,x_2,\ldots]\ra B$ defined by
$ax_{i_1}^{n_{i_1}}\cdots x_{i_k}^{n_{i_k}}\mapsto f(a)b_{i_1}^{n_{i_1}}\cdots
b_{i_k}^{n_{i_k}}$. The congruence $R$ occurs as the pullback of this surjection
with itself (see Proposition \ref{p:monoidcong}).\\
ii) Use (i) considering $A$ as a $A/\frakm$-monoid via the obvious inclusion.
Then $A\backslash (A/\frakm)=\frakm\backslash\{0\}$ is the nonzero elements of
maximal ideal.  When $\frakm$ is finitely generated, only finitely many $x_i$
are necessary.
\qed
\end{myproof}

\begin{remark}\label{r:createmonoid}
Proposition~\ref{p:monoidstructure} provides a simple way to create many
monoids.  Beginning with a monoid $A$, the $A/\frakm$-monoid structure is the
simplest way to describe it.  Alternatively, let $G_+$ be a pointed abelian
group (hence monoid) and let $S$ be any commutative, pointed semigroup.  Adding
a distinguished identity $S'=S\coprod \{1\}$ makes $S'$ a monoid and we may
tensor over $\F_1$: $A=G_+\otimes S'$.  Note that if $S$ is already a monoid
with identity $e$, then $e$ is an (nontrivial) idempotent element of $S'$.
Let $S'\ra\F_1$ be the map sending all non-identity elements to 0.  Tensoring 
with $G_+$ induces the map $A\ra G_+$ having kernel $S$. Since $G_+$ is a 
pointed group, $(G\times S)/(G\times 0)$ is the maximal ideal of $A$.
\end{remark}

\section{Localization}\label{localization}

Let $A$ be a monoid and $S\subseteq A$ a multiplicatively closed subset.  Define $S^{-1}A$ to be the monoid with elements $a/s$, $a\in A$ and $s\in S$, where $a/s=b/t$ if there is a $u\in S$ such that $u(at)=u(bs)$.  The multiplication in $S^{-1}A$ is induced by $A$, $(a/s)(b/t)=ab/st$.  Note that $(1/s)(s/1)=1$ so that any element of $S$ becomes a unit in $S^{-1}A$.  Clearly $S\subseteq A^{\times}$ means $S^{-1}A=A$.  The monoid $S^{-1}A$ is called the \emph{monoid of fractions of $A$ with respect to $S$} or the \emph{localization of $A$ at $S$}.

There are special cases which warrant their own notation.  When $\frakp\in\MSpec(A)$, $S=A\backslash\frakp$ is multiplicatively closed and $S^{-1}A$, denoted $A_{\frakp}$, is the \emph{localization of $A$ at $\frakp$}. When $S=\{s,s^2,\ldots\}$ is generated by a single element, we write $S^{-1}A=A_s$ or $A[\frac{1}{s}]$.  More generally, when $S$ is generated by $s_1,\ldots,s_n$ we write $S^{-1}A=A[\frac{1}{s_1},\ldots,\frac{1}{s_n}]$.  Note there is a canonical map $A\ra S^{-1}A$ defined by $a\mapsto a/1$ which is injective only when $S\cup\{0\}\subseteq A$ is a cancellative monoid.

\begin{remark}
The monoid $S^{-1}A$ satisfies the usual universal property of localization.
Namely, let $f:A\ra B$ be any morphism such that the image of every element of
$S\subseteq A$ is a unit in $B$.  Then $f$ factors through the morphism $A\ra
S^{-1}A$ as in the following diagram:
\[\SelectTips{eu}{12}\xymatrix{
A \ar[dr]_f \ar[r] & S^{-1}A \ar @{-->} [d] ^{\exists} \\
 & B
} \]
\end{remark}

\begin{remark}
A multiplicatively closed subset $S\subseteq A$ is \emph{saturated} when
$xy\in S$ implies $x,y\in S$.  If $S$ is any multiplicatively closed subset of
$A$, the \emph{saturation} $\overline{S}$ of $S$ is the intersection of all the
saturated multiplicatively closed subsets of $A$ containing $S$.  Moreover, the
complement of a saturated, multiplicatively closed subset of $A$ is a prime
ideal.
\end{remark}

\begin{proposition}\label{p:localizationprime}
Let $A$ be a monoid and $S\subseteq A\backslash\{0\}$ a multiplicatively closed
subset. Then $S^{-1}A$ exists and:
\begin{compactenum}
\item $S^{-1}A=A_{\frakp}$ for some $\frakp\in\MSpec(A)$.
\item If $A$ is finitely generated, $A_\frakp=A_s$ for some $s\in A\backslash
    \frakp$.
\item The proper ideals of $S^{-1}A$ correspond to the ideals of $A$ contained
    in $A\backslash S$.
\end{compactenum}
\end{proposition}

\begin{myproof}
i) First notice that $S^{-1}A=\overline{S}^{-1}A$, since if $ab\in S$, then $\frac{1}{a}=b\frac{1}{ab}\in S^{-1}A$.  Similarly, $\frac{1}{b}=a\frac{1}{ab}\in S^{-1}A$.  Therefore we may assume that $S$ is saturated.  Then $xy\in A\backslash S$ and $x,y\not\in A\backslash S$, that is $x,y\in S$, implies $xy\in S$.  Hence, $A\backslash S$ is a prime ideal. \\
ii) Let $a_1,\ldots,a_n$ be the generators for $A$.  Since $S=A\backslash\frakp$
is saturated, it must be generated by the $a_i\not\in\frakp$.  For
$s=\prod_{a_i\not\in\frakp}a_i$, we have $A_\frakp=A_s$.\\
iii)  If $I\subseteq A$ is an ideal and $I\cap S\neq\emptyset$, then $S^{-1}I$
contains units, hence is not proper.  Conversely, if $J\subseteq S^{-1}A$ is an
ideal, then $J=S^{-1}I$ where $I=\{a\in A\ |\ a/1\in J\}$.
\qed
\end{myproof}



There are two special monoids obtained from localization:  the \emph{total
monoid of fractions} and \emph{group completion}.  When $S$ is the set of
non-zero divisors of $A$, $G(A)=S^{-1}A$ is the \emph{total monoid of
fractions}.  When $A$ is torsion free, i.e. $0$ generates a prime ideal, the
localization $A_0=A_{(0)}$ of $A$ at $(0)$ is the \emph{group completion}.  The
canonical map $A\ra A_0$ is an inclusion precisely when $A$ is cancellative.

Let $X$ be an $A$-set and $S\subseteq A$ multiplicatively closed.  Define
$S^{-1}X$, the \emph{localization of $X$ at $S$}, to be the $S^{-1}A$-set with
elements $x/s$, $x\in X$ and $s\in S$, where $x/s=x'/t$ when $u(tx)=u(sx')$ for
some $u\in S$.  The action of $S^{-1}A$ on $S^{-1}X$ is simply
$(a/s)(x/t)=ax/st$.

The $S^{-1}A$-set $S^{-1}X$ satisfies the usual universal property.  Consider an
$S^{-1}A$-set $Y$ as an $A$-set by restriction of scalars.  Then every $A$-set
morphism $f:X\ra Y$ factors through $X\ra S^{-1}X$ as in the following diagram:
\[\SelectTips{eu}{12}\xymatrix{
X \ar[dr]_f \ar[r] & S^{-1}X \ar @{-->} [d] ^{\exists} \\
 & Y
} \]
When $S=A\backslash\frakp$, the notation $X_{\frakp}$ will denote $S^{-1}X$ and
we call $X_{\frakp}$ the \emph{localization of $X$ at $\frakp$}.  Note that the
previous notations apply when $X=I$ is an ideal of $A$, namely $S^{-1}I$ is the
image of $I$ in $S^{-1}A$ and likewise for $I_{\frakp}$.

\begin{proposition}\label{p:basicfraction}
Let $X$ be an $A$-set and $S\subseteq A$ multiplicatively closed.
\begin{compactenum}
    \item Let $0\neq x\in X$. Then $\ann(x)\cap S\neq \emptyset$ if and only if
        $x/1=0$ in $S^{-1}X$.
    \item $S^{-1}X\cong X\otimes_A S^{-1}A$.
\end{compactenum}
\end{proposition}

\begin{myproof}
i) Suppose $a\in\ann(x)\cap S$.  Then $\frac{1}{a}\in S^{-1}A$ and
$x=(\frac{1}{a}a)x=\frac{1}{a}0=0$. Conversely, suppose $0\neq x\in X$ but
$x/1=0$ in $S^{-1}X$. By definition, there exists $u\in S$ with $ux=u0=0$, hence
$u\in\ann(x)\cap S$.\\
ii)  The map $X\times S^{-1}A\ra S^{-1}X$ defined by $(x,a/s)\mapsto ax/s$ is
bi-equivariant and, by the universal property of tensor products, induces
$f:X\ot_A S^{-1}A\ra S^{-1}X$.  The induced map is surjective as $x\ot
1/s\mapsto x/s$ for any $x\in X$ and $s\in S$.  Now if $f(x\ot
a/s)=ax/s=bx'/t=f(x'\ot b/t)$, there exists $u\in S$ such that $u(tax)=u(sbx')$
in $X$.  Hence,
\[x\ot \frac{a}{s}=x\ot\frac{aut}{sut}=u(tax)\ot\frac{1}{sut}=u(sbx')\ot
\frac{1}{sut}=x'\ot\frac{sub}{sut}=x'\ot\frac{b}{t}.\]
\end{myproof}



%

\section{Chain conditions}\label{s:chainconditions}

We say an $A$-set $X$ is \emph{noetherian} when it satisfies the ascending chain condition (ACC) on $A$-subsets.  That is, every increasing chain of $A$-subsets $X_1\subseteq X_2\subseteq\cdots$ stabilizes so that $X_n=X_{n+1}=\cdots$ for some $n>0$.  Similarly, $X$ is \emph{artinian} when it satisfies the descending chain condition (DCC) on $A$-subsets, namely every descending chain of $A$-subsets $X_1\supseteq X_2\supseteq\cdots$ stabilizes.  A monoid $A$ is \emph{noetherian} (respectively \emph{artinian}) when it is noetherian (respectively artinian) when considered as an $A$-set.  Since $A$-subsets of $A$ are simply ideals, $A$ is noetherian (respectively artinian) when it satisfies the ACC (respectively DCC) on ideals.  The proof of the following result is an exact replica of the analogous result in ring theory.

\begin{proposition}
An $A$-set $X$ is noetherian if and only if every $A$-subset is finitely
generated.  Hence, $A$ is noetherian if and only if every ideal is finitely
generated.
\end{proposition}

%
%


We now verify some standard results for noetherian $A$-sets and monoids.

\begin{proposition}\label{p:f1xnoeth}
The free monoid on one variable $\F_1[x]$ is noetherian.
\end{proposition}

\begin{myproof}
Every ideal is generated by its element of lowest degree.
\qed
\end{myproof}

\begin{proposition}\label{p:aesnoetherian}
Let $0\ra X'\xra{f}X\xra{g}X''\ra0$ be an a.e.s. of $A$-sets. Then
\begin{compactenum}
    \item $X$ is noetherian if and only if $X'$ and $X''$ are noetherian.
    \item $X$ is artinian if and only if $X'$ and $X''$ are artinian.
    \item If $X_i$ is a noetherian (resp. artinian) $A$-set for $1\leq i\leq
        n$, then $\vee_{i=1}^n X_i$ is noetherian (resp. artinian).
\end{compactenum}
\end{proposition}

\begin{myproof}
i) Suppose $X$ is noetherian.  Any ascending chain of $A$-subsets of $X'$ is
also an ascending chain in $X$ since $X'\subseteq X$ and thus, stabilizes.
Likewise, an ascending chain in $X''$ corresponds to an ascending chain in $X$
under $g^{-1}$.

Conversely, let $\{X_i\}_{i\geq 0}$ be an ascending chain in $X$.  Then the
$\{f^{-1}(X_i)\}$ and $\{g(X_i)\}$ form ascending chains in $X'$ and $X''$ which
stabilize.  Since each $X_i$ is the union of the image of $f^{-1}(X_i)$ and the
inverse image of $g(X_i)\backslash\{0\}$, the $\{X_i\}$ also stabilize.\\
ii) Similar to (i).\\
iii) Use induction with (i),(ii) on the a.e.s. $0\ra X_n\ra \bigvee_{i=1}^n
X_i\ra \bigvee_{i=1}^{n-1}X_i\ra 0$.
\qed
\end{myproof}

\begin{proposition}\label{p:asetimagenoetherian}
Let $f:X\ra Y$ be an $A$-set morphism.  If $X$ is noetherian, so is $f(X)$.
\end{proposition}

\begin{myproof}
Let $X'\subseteq f(X)$ be an $A$-subset.  Let $f^{-1}(X')\subseteq X$ denote the
set theoretic inverse of elements in $X'$.  Then  $f^{-1}(X')$ is finitely
generated by, say, $x_1,\ldots,x_n$.  Hence, the $f(x_i)$ form a finite set of
generators for $X'$.
\qed
\end{myproof}


\begin{proposition}\label{p:tensornoetherian}
Let $A,B$ be noetherian monoids and $A[x]$ the free $A$-monoid in one variable.
Then
\begin{compactenum}
    \item $A\otimes B$ is noetherian.
    \item $A[x_1,\ldots,x_n]$ is noetherian.
\end{compactenum}
\end{proposition}

\begin{myproof}
i)  By Proposition \ref{p:tensorideals}, every ideal $K\subseteq A\ot B$ is of
the form $K=\bigcup_{\lambda\in\Lambda}I_{\lambda}\ot J_{\lambda}$ for some
indexing set $\Lambda$.  Since $A,B$ are noetherian, for fixed $\lambda$, both
$I_{\lambda}=(a_1,\ldots,a_n)$ and $J_{\lambda}=(b_1,\ldots,b_m)$ are finitely
generated so that $I_{\lambda}\ot J_{\lambda}$ is finitely generated by
$\{a_i\otimes b_j\}$.  Thus, we need only show that the indexing set $\Lambda$
can be made finite.

Proposition \ref{p:tensorideals} also shows that at least one of the sets
$\{I_{\lambda}\}$ and $\{J_{\lambda}\}$ have distinct ideals, say
$I_{\lambda}\neq I_{\lambda'}$ for $\lambda\neq\lambda'$.  If $\Lambda$
were infinite, we can construct a strictly increasing chain of ideals
$I_{\lambda_1}\subseteq I_{\lambda_1}\cup I_{\lambda_2}\subseteq\cdots$, where
$\lambda_{n+1}\not\in\{\lambda_1,\ldots,\lambda_n\}$, in $A$ which contradicts
the noetherian assumption.\\
ii) Use Proposition \ref{p:f1xnoeth} and (i) inductively.
\qed
\end{myproof}

\begin{proposition}
Let $A$ be a noetherian monoid.
\begin{compactenum}
\item Any homomorphic image of $A$ is noetherian.
\item Finitely generated $A$-monoids are noetherian.
\item If $S\subseteq A$ is multiplicatively closed, then $S^{-1}A$ is
    noetherian.
\end{compactenum}
\end{proposition}

\begin{myproof}
i) Use Proposition \ref{p:asetimagenoetherian} on the ideals of $A$. \\
ii) A finitely generated $A$-monoid is the homomorphic image of
$A[x_1,\ldots,x_n]$ for some $n>0$.  The latter monoid is noetherian by
Proposition \ref{p:tensornoetherian}(ii).  Now use (i).\\
iii)  By Proposition \ref{p:localizationprime}(iii), the ideals of $A$ contained
in $S^{-1}A$ correspond to ideals contained in $A\backslash S$. If
$I=(x_1,\ldots,x_n)$ is finitely generated, $S^{-1}I$ is generated by $x_i/1$.
\qed
\end{myproof}

\begin{proposition}\label{p:noetherianaset}
Let $A$ be a noetherian monoid and $X$ an $A$-set.  If $X$ is finitely
generated, then it is noetherian.
\end{proposition}

\begin{myproof}
Let $x_1,\ldots,x_n$ be a set of generators for $X$.  Then there is a surjection
$\bigvee_{i=1}^nA_i\ra X$ defined by $1_i\mapsto x_i$.  Now use Proposition
\ref{p:asetimagenoetherian}.
\qed
\end{myproof}

Before proceeding we recall a definition.  Consider the following diagram of
$A$-sets:
\[X\da{f}{g}Y\xra{h}Z.\]
We say that $h$ \emph{coequalizes} $f$ and $g$ when $hf=hg$ and that $h$, or
$Z$, is the \emph{coequalizer} of $f,g$ when it satisfies the following
universal property:
\[\SelectTips{eu}{12}\xymatrix{
   &   & Z'\\
    X \ar@<2pt>[r]^f \ar@<-2pt>[r]_g & Y \ar [ur]^{h'} \ar [r]^h & Z
    \ar@{-->}[u]_{\exists\, p}
}\]
namely, if $h':Y\ra Z'$ is any other morphism that coequalizes $f$ and $g$,
there exists a morphism $p:Z\ra Z'$ such that $h'=ph$.  Of course, this is just
the universal product for colimits applied to the diagram
$\cdot\rightrightarrows \cdot$.  We sometimes represent the codomain of the
coequalizer map by $Z=\coeq(f,g)$.  Coequalizers play a larger role in the
theory of monoids (resp.  $A$-sets) since every monoid (resp. $A$-set) occurs as
a coequalizer, but not as a quotient by an $A$-subset.  See Chapter
\ref{c:homological} for more details.


In general there are many more congruences on a monoid than there are ideals;
obviously every ideal defines a congruence.  Then why not
define a monoid $A$ to be noetherian when it satisfies the DCC on quotient
monoids?  Equivalently, when $A$ has the ACC on congruences.  We will see in this section that monoids that have the ACC on congruences coincide with finitely generated monoids.

The following result is clear, as $A$-subsets determine congruences:

\begin{lemma}\label{p:conimpliessubset}
Let $X$ be an $A$-set.  If $X$ has the ACC on congruences, then $X$ is
noetherian.
\end{lemma}

A simple example of a noetherian monoid which does not have the ACC on
congruences is the following: let $G$ be the free abelian group on countably
many generators $x_i$.  Then $G_+$ is a monoid and $(x_1\sim 0)\subseteq
(x_1\sim 0, x_2\sim 0)\subseteq(x_1\sim 0,x_2\sim 0,x_3\sim 0)\subseteq\cdots$
is an ascending chain of congruences (whose generators are shown) on $G_+$ which
clearly does not stabilize.  However, $G_+$ is noetherian since the only proper ideal is $(0)$.


We now briefly introduce a functor we require for subsequent results.  For more
information see Section \ref{s:functors}.  Given a commutative ring $k$ and
monoid $A$, we may form the $k$-algebra $k[A]$ which is the quotient of the free
$k$-module with generators the elements of $A$ by the submodule generated by
$0_A$.  A generic element of $k[A]$ is a finite $k$-linear sum of nonzero elements of $A$ with multiplication provided by both $A$ and $k$.  This construction induces a functor $k[-]:\mon\ra k{\text{\bf -algebras}}$, called the \emph{$k$-realization} of $A$,
defined by $A\mapsto k[A]$. In particular $k[\F_1]=k$.  We can extend this
definition to obtain $k[-]: A\set\ra k[A]\mod$ defined by $X\mapsto k[X]$, where $k[X]$ is the
quotient of the free $k[A]$-module with generators the elements of $X$ by the
submodule generated by $0_X$.  Of course, a generic element of $k[X]$ is a
finite $k$-linear sum of nonzero elements of $X$.

\begin{proposition}\label{p:fgkrealizationnoeth}
If $A$ is a finitely generated monoid and $k$ is a noetherian ring, then $k[A]$
is a noetherian ring.
\end{proposition}

\begin{myproof}
Since $A$ is finitely generated by, say $a_1,\ldots,a_n$, we may write it as the quotient of a free monoid in $n$ variables $\langle 0,1,x_1,\ldots,x_n\rangle$. Then $k[A]$ is the homomorphic image of the polynomial ring $k[x_1,\ldots,x_n]$, which is noetherian by the Hilbert Basis Theorem in commutative ring theory.
\qed
\end{myproof}

\begin{lemma}\label{l:equalcongruence}
Let $X$ be an $A$-set and $R_1\subseteq R_2$ be congruences on $X$.  Let
$p_i:X\ra X/R_i$ and $k[p_i]:k[X]\ra k[X/R_i]$ be the usual projection maps and
their $k$-realizations respectively.  Then $R_1= R_2$ if and only if
$\ker(k[p_1])= \ker(k[p_2])$.
\end{lemma}

\begin{myproof}
Consider the following diagram:
\[\SelectTips{eu}{12}\xymatrix{
    0 \ar[r]  & \ker(k[p_1]) \ar@{-->}[d]_f \ar[r] & k[X] \ar[r]^{k[p_1]\quad}
    \ar@{=}[d] & k[X/R_1] \ar[r] \ar@{-->}[d]^{k[q]} & 0  \\
    0 \ar[r]  & \ker(k[p_2]) \ar[r] & k[X] \ar[r]_{k[p_2]\quad} & k[X/R_2]
    \ar[r] & 0  \\
} \]
where $q:X/R_1\ra X/R_2$ is the surjection induced by the containment
$R_1\subseteq R_2$.  If $R_1=R_2$, then $q$, and hence $k[q]$, is the identity
map so that $f$ is the identity map.  If $\ker(k[p_1])=\ker(k[p_2])$, then $f$
is the identity map and so is $k[q]$.  Certainly $q:X/R_1\ra X/R_2$ is the
identity map and hence, $R_1=R_2$.
\qed


\end{myproof}

\begin{proposition}\label{p:noethcon}
An $A$-set $X$ has the ACC on congruences when $k[X]$ is a noetherian
$k[A]$-module.
\end{proposition}

\begin{myproof}
Let $R_1\subseteq R_2\subseteq\cdots$ be an ascending chain of congruences on
$X$ and $X/R_1\ra X/R_2\ra\cdots$ the associated descending chain of quotient
$A$-sets.  Let $p_i:X\ra X/R_i$ denote the projection map and consider the
$k$-realizations $k[p_i]:k[X]\ra k[X/R_i]$.

As in Lemma~\ref{l:equalcongruence}, the $k[p_i]$ give rise to an increasing
chain of $k[A]$-submodules $\ker(k[p_1])\subseteq \ker(k[p_2])\subseteq \cdots$
that stabilizes since, by assumption, $k[X]$ is noetherian.  Hence, there is an
$N\geq 0$ such that
\[\ker(k[p_i])= \ker(k[p_{i+1}])=\cdots\]
for all $i\geq N$.  By Lemma \ref{l:equalcongruence}, we have
$R_i=R_{i+1}=\cdots$ for all $i\geq N$.
\qed
\end{myproof}

\begin{corollary}\label{c:accnoethequiv}
Let $A$ be a finitely generated monoid and $X$ an $A$-set.  Then $X$ has the ACC
on congruences if and only if $X$ is noetherian.
\end{corollary}

\begin{myproof}
When $X$ has the ACC on congruences, it also has the ACC on $A$-subsets by Lemma \ref{p:conimpliessubset}.  Conversely, suppose $X$ is noetherian and let $k$ be any commutative, noetherian ring.  By Proposition \ref{p:fgkrealizationnoeth}, $k[A]$ is a noetherian ring and since $X$ is noetherian, it is finitely generated.  Then $k[X]$ is finitely generated as a $k[A]$-module, hence noetherian.  Now use Proposition \ref{p:noethcon}.
\qed
\end{myproof}

We now summarize the results:

\begin{theorem}\label{t:noetheriansummary}
Let $A$ be a monoid.  Then $A$ is finitely generated if and only if $A$ has the ACC on congruences.
\end{theorem}

\begin{myproof}
Assuming $A$ is finitely generated, the result follows from Proposition \ref{p:tensornoetherian} and Corollary \ref{p:noethcon} with $X=A$.  Now, assume $A$ has the ACC on congruences and suppose it is not finitely generated, say $A=\F_1[a_1,a_2,\ldots]$.  Then we can construct an non-stabilizing chain of congruences $R_1\subseteq R_2\subseteq\cdots$ where $R_i$ is generated by relations $\{a_j\sim\epsilon(a_j)\}_{j\leq n}$ and $\epsilon(a_j)=1$ when $a_j\in A^{\times}$ and 0 otherwise.
\qed
\end{myproof}

It is not clear whether the importance of monoids having the ACC on congruences
will outweigh noetherian monoids.  Thus we leave the definition of noetherian
fixed and consider the situation nothing more than an inconvenience.


\begin{lemma}\label{rad.power}
Let $A$ be a noetherian monoid.  For any ideal $I$ of $A$,
$(\sqrt{I})^n\subseteq I$ for some $n$.
\end{lemma}

\begin{myproof}
Let $\sqrt{I}=(x_1,\cdots,x_k)$ with $x_i^{n_i}\in I$ for $n_i>0$, $1\leq i\leq
k$.  Let $m=\sum_i(n_i-1)+1$.  Then $(\sqrt{I})^m$ is generated by the products
$x_1^{r_1}\cdots x_k^{r_k}$ with $\sum r_i=m$.  From the definition of $m$ we
must have $r_i\geq n_i$ for some $1\leq i\leq k$, hence every monomial
generating $\sqrt{I}$ is contained in $I$.
\end{myproof}

\subsection{Artinian monoids}\label{artinianmonoids}

The theory of artinian monoids is not as simple as that of artinian rings due to the lack of cancellation.  It is a well known result that a ring is artinian only when it is noetherian and has dimension 0.  This is not the case for artinian monoids, in fact there are artinian monoids of any dimension.

\begin{example}\label{e:anydimartinian}
Let $A=\F_1[x_1,\ldots,x_n]/(x_1^2=x_1,\ldots,x_n^2=x_n)$.  Then $A$ is artinian since $A=\{0,1,x_1^{\epsilon_1}\cdots x_n^{\epsilon_n}\ | \ \epsilon_i = 0,1\}$ is a finite set.  Here $A$ has $2^n$ primes and Krull dimension $n$ (see Section
\ref{primeidealsandunits}).

It is also not the case that artinian monoids need are noetherian.  For example,
the maximal ideal of the monoid $A=\F_1[x_1,x_2,\ldots\ |\ x_i x_j=x_i\text{
when } i\leq j]$ is not finitely generated, but the ideals of $A$ satisfy
the DCC.  The ideals of $A$ are $(0), Ax_i$, $i\geq 1$, and the maximal ideal $\frakm=(x_1,x_2,\ldots)$; also $Ax_i\subseteq Ax_j$ only when $i\leq j$.  Hence, any fixed ideal $Ax_n$ contains only finitely many ideals. 
\end{example}





\begin{lemma}
Let $A$ be a monoid, $\frakm\subseteq A$ its maximal ideal and $X$ an $A$-set. If $X$ is finitely generated as an $A/\frakm$-set, then $X$ is artinian and noetherian.
    \label{l:artin.finite}
\end{lemma}

\begin{myproof}
Since $G_+=A/\frakm$ is a pointed abelian group, $X$ is a finite wedge sum of 
orbits of $G$ considered as a $G_+$-set; namely $X=\bigvee_{i=1}^n (G/H_i)_+$ 
where every $H_i$ is a subgroup of $G$.  Write $x_1,\ldots,x_n$ for the 
generators.  Then every $A$-subset $X'\subseteq X$ is generated by a subset of 
the $x_i$ so that every descending chain $X_1\supseteq X_2\supseteq\cdots$ must 
stabilize.
\qed
\end{myproof}

\begin{theorem}\label{t:artinian}
Let $A$ be a monoid and $\frakm\subseteq A$ its maximal ideal.  If 
\begin{compactenum}
\item $\frakm$ is finitely generated as an ideal,
\item $\frakm^n=\frakm^{n+1}$ for some $n>0$, and
\item $\frakm^n$ is finitely generated as an $A/\frakm$-set
\end{compactenum}
then $A$ is artinian.  Of course (i) may be removed if $A$ is noetherian.
\end{theorem}

\begin{myproof}
If $x_1,\ldots,x_n$ are a set of generators for $\frakm$, then $\frakm/\frakm^n$ 
is generated as an $A/\frakm$-set by all monomials in the $x_i$ of degree less 
than $n$.  Since, as $A/\frakm$-sets, $\frakm\cong (\frakm/\frakm^n)\vee 
\frakm^n$, $\frakm$ is an artinian $A$-set by Lemma \ref{l:artin.finite}.  
Hence, $A$ is an artinian monoid.
\qed
\end{myproof}

Although Example \ref{e:anydimartinian} shows not every artinian monoid is 
0-dimensional, the converse holds.

\begin{corollary}
A 0-dimensional noetherian monoid $A$ is artinian.
\end{corollary}

\begin{myproof}
As every monoid is local, $A$ has at least one prime ideal and this must be the
maximal ideal $\frakm$ of $A$.  Moreover, $\frakm$ is (trivially) the
intersection of all prime ideals in $A$ so by \ref{t:p.decomp} we have $\frakm=\nil(A)$.  By Lemma \ref{rad.power}, $\frakm^n=0$ for some $n>0$  and since $A$ is noetherian, $\frakm$ is finitely generated.  Hence $A$ is artinian 
by Theorem \ref{t:artinian}.
\qed
\end{myproof}

\section{Primary decomposition}\label{primarydecomposition}

Much of the theory of primary decomposition in commutative ring theory carries
through to monoids since it does not rely upon the underlying abelian group structure of the ring. Below
we further develop the theory of noetherian monoids since these are of primary
interest.

\begin{lemma}\label{l:irred.prim}
Every irreducible ideal of a noetherian monoid is primary.
\end{lemma}

\begin{myproof}
An ideal $I\subseteq A$ is primary if and only if the zero ideal of $A/I$ is
primary.  Therefore we need only show when $(0)$ is irreducible, it is primary.
Let $xy=0$ with $y\neq 0$; we will show $x^n=0$.  Consider the ascending chain
$\ann(x)\subseteq \ann(x^2)\subseteq\cdots$, where $\ann(x)=\{a\in A\ | \
ax=0\}$.  By the noetherian property, this chain must stabilize, say
$\ann(x^n)=\ann(x^{n+1})=\cdots$ for some $n>0$.

We claim that $0=(x^n)\cap (y)$.  Let $a\in(x^n)\cap(y)$, say $a=bx^n=cy$. Then
$0=c(xy)=ax=(bx^n)x=bx^{n+1}$.  Hence $b\in\ann(x^{n+1})=\ann(x^n)$ giving
$a=bx^n=0$.  Since $(0)$ is irreducible and $y\neq 0$, we must have $x^n=0$
proving $(0)$ is primary.
\end{myproof}

%

\begin{theorem}[Primary decomposition]\label{t:p.decomp}
In a noetherian monoid every ideal $I$ can be written as the finite intersection
of irreducible primary ideals $I=\cap_i\frakq_i$.
\end{theorem}

\begin{myproof}
Suppose the result is false.  Since $A$ is noetherian, the set of ideals which
cannot be written as a finite intersection of irreducible ideals has a maximal
element, say $I$.  Since $I$ is not irreducible, it can be written $I=J\cap K$
where $J,K$ are ideals of $A$ containing $I$.  By maximality, both $J$ and $K$
(and hence $I$) can be written as a finite intersection of irreducible ideals.
This is a contradiction, and the theorem follows via Lemma \ref{l:irred.prim}.
\end{myproof}

\begin{remark}\label{r:nil(A)}
Recall that $A$ is {\it reduced} if whenever $a,b\in A$ satisfy $a^2=b^2$ and
$a^3=b^3$ then $a=b$. This implies that $A$ has no nilpotent elements, i.e.,
that $\nil(A)=\{ a\in A: a^n=0 \text{ for some } n\}$ vanishes.  By Theorem
\ref{t:p.decomp}, $\nil(A)$ is the intersection of all the prime ideals in $A$;
cf.\ \cite[1.1]{chww-p}.
When $A$ is a pc monoid, it is reduced if and only if $\nil(A)=0$; in this case
$A_\red=A/\nil(A)$ is a reduced monoid (see \cite[1.6]{chww-p}).  Note that the
equivalence of these two conditions does not hold in general.
\end{remark}

Let $A$ be a noetherian monoid and $I\subseteq A$ an ideal.  Given a minimal
primary decomposition of $I$, $I=\cap_i\frakq_i$, $\Ass(I)$ denotes the set of
prime ideals occurring as the radicals $\frakp_i=\sqrt{\frakq_i}$; the
$\frakp_i$ are called the associated primes of $I$. Although the primary
decomposition need not be unique, the set $\Ass(I)$ of associated primes of $I$
is independent of the minimal primary decomposition, by \ref{p:ap.radquot}
below.

\begin{proposition}\label{p:ap.radquot}
Let $I\subseteq A$ be an ideal with minimal primary decomposition
$I=\cap_{i=1}^n\frakq_i$ where $\frakq_i$ is $\frakp_i$-primary.  Then $\Ass(I)$
is exactly the set of prime ideals which occur as $\sqrt{(I:a)}$, where $a\in
A$.  Hence $\Ass(I)$ is independent of the choice of primary decomposition.

In addition, the minimal elements in $\Ass(I)$ are exactly the set of prime
ideals minimal over $I$.
\end{proposition}

\begin{myproof} (Compare \cite[4.5, 4.6]{AM}.)
First note that
\[\sqrt{(I:a)} = \sqrt{(\cap_i\frakq_i:a)} = \sqrt{\cap_i(\frakq_i:a)} = \cap_i
\sqrt{(\frakq_i:a)}. \]
By Lemma~\ref{l:pcont.quotprim}(ii), this equals
$\cap_{a\not\in\frakq_i}\frakp_i$.  If $\sqrt{(I:a)}$ is prime, then it is
$\frakp_i$ for some $i$, by Lemma~\ref{l:pcont.quotprim}(i).  Conversely, by
minimality of the primary decomposition, for each $i$ there exists an
$a_i\not\in\frakq_i$ but $a_i\in\cap_{j\neq i}\frakq_j$.  Using
Lemma~\ref{l:pcont.quotprim}(i) once more, we see $\sqrt{(I:a_i)}=\frakp_i$.

Finally, if $I\subseteq\frakp$ then $\cap\frakp_i\subseteq\frakp$, so $\frakp$
contains some $\frakp_i$ by Lemma \ref{l:pcont.quotprim}(i).  If $\frakp$ is
minimal over $I$ then necessarily $\frakp=\frakp_i$.
\end{myproof}

\begin{proposition}\label{p:ap.ann}
Let $A$ be a noetherian monoid and $I\subseteq A$ an ideal.  Then the associated
prime ideals of $I$ are exactly the prime ideals occurring in the set of ideals
$(I:a)$ where $a\in A$.
\end{proposition}

\begin{myproof}
The ideals $I_i=\cap_{j\neq i}\frakq_j$ strictly contain $I$ by minimality of
the decomposition.  Since $\frakq_i\cap I_i=I$, any $a\in I_i\backslash I$ is
not contained in $\frakq_i$, hence $(I:a)$ is $\frakp_i$-primary by Lemma
\ref{l:pcont.quotprim}(ii).  Now, by Lemma~\ref{rad.power} we have
$\frakp_i^n\subseteq\frakq_i$ for some $n>0$, hence
\[\frakp_i^nI_i\subseteq \frakq_iI_i\subseteq \frakq_i\cap I_i=I. \]
Choose $n$ minimal so that $\frakp_i^nI_i\subseteq I$ (hence in $I_i$) and pick
$a\in\frakp_i^{n-1}I_i$ with $a\not\in I$. Since $\frakp_ia\subseteq I$ we have
$\frakp_i\subseteq (I:a)$; as $(I:a)$ is $\frakp_i$-primary, we have
$\frakp_i=(I:a)$.  Conversely, if $(I:a)$ is prime then it is an associated
prime by Proposition~\ref{p:ap.radquot}.
\end{myproof}

\section{Normal and factorial monoids}

In this section, we establish the facts about normal monoids needed for the
theory of divisors.

The vocabulary for integral extensions of monoids mimics that for commutative
rings. If $A$ is a submonoid of $B$, we say that an element $b\in B$ is
\emph{integral} over $A$ when $b^n\in A$ for some $n>0$, and the {\it integral
closure} of $A$ in $B$ is the submonoid of elements integral over $A$. If $A$ is
a cancellative monoid, we say that it is {\it normal} (or \emph{integrally
closed}) if it equals its integral closure in its group completion.  (See
\cite[1.6]{chww}.)

\begin{example}\label{e:normal}
It is elementary that all factorial monoids are normal.  Affine toric monoids
are normal by \cite[4.1]{chww}, and so are arbitrary submonoids of a free
abelian group closed under divisibility.  By \cite[4.5]{chww}, every finitely
generated normal monoid is $A\wedge U_*$ for an affine toric monoid $A$ and a
finite abelian group $U$.
\end{example}

One difference between integrality in commutative ring theory and monoids again
results from lack of cancellation.  Let $A, B$ be monoids and $x\in B$.
Certainly, if $x$ is integral over $A$, then $x^n\in A$ for some $n>0$ so that
$A[x]$ is finitely generated as an $A$-set.  However, the converse \emph{does
not hold}.  Consider the monoid $A[x]/(x^n=x)$ which is finitely generated as an
$A$-set by $1,x,x^2,\ldots,x^{n-1}$ even though $x$ is \emph{not} integral over
$A$.

\begin{lemma}\label{l:local.normal}
Let $A\subseteq B$ be monoids and $S\subseteq A$ be multiplicatively closed.  We
have the following:
\begin{compactenum}
\item[i)] If $B$ is integral over $A$, then $S^{-1}B$ is integral over
    $S^{-1}A$.
\item[ii)]  If $B$ is the integral closure of $A$ in a monoid $C$, then
    $S^{-1}B$ is the integral closure of $S^{-1}A$ in $S^{-1}C$.
\item[iii)]  If $A$ is normal, then $S^{-1}A$ is normal.  More generally, if $B$
    is the normalization of $A$, then $S^{-1}B$ is the normalization of
    $S^{-1}A$.
\end{compactenum}
\end{lemma}

\begin{myproof}
Suppose that $b$ is integral over $A$, i.e., $b^n\in A$ for some $n>0$.  Then
$b/s\in S^{-1}B$ is integral over $S^{-1}A$ because $(b/s)^n\in S^{-1}A$. This
proves (i). For (ii), it suffices by (i) to suppose that $c/1\in S^{-1}C$ is
integral over $S^{-1}A$ and show that $c/1$ is in $S^{-1}B$. If $(c/1)^n=a/s$ in
$S^{-1}A$ then $c^nst=at$ in $A$ for some $t\in S$. Thus $cst$ is in $B$, and
$c/1=(cst)/st$ is in $S^{-1}B$.  It is immediate that (ii) implies (iii) (see
also \cite[1.6]{chww}).
\end{myproof}

Recall that when $A$ is a torsion free monoid, the group completion $A_0$ is
obtained from $A$ by localizing at the trivial prime ideal $(0)$.  When $A$ is
cancellative, $A_0$ contains $A$ as a submonoid.

\begin{lemma}\label{l:cayley.int}
Let $A\subseteq B$ be monoids with $A$ cancellative.  If there is a finitely generated ideal $I\subseteq A$ and a $b\in B$ such that $b I\subseteq I$, then $b$ is integral over $A$.
\end{lemma}

\begin{myproof}
Let $X=\{x_1,\ldots,x_r\}$ be the set of generators of $I$.  Since $bI\subseteq I$, there is a function $\phi:X\to X$ such that $bx\in A\phi(x)$ for each $x\in X$.  Since $X$ is finite, there is an $x\in X$ and an $n$ so that $\phi^n(x)=x$. For this $x$ and $n$ there is an $a\in A$ such that $b^nx=ax$ and hence, $b^n=a$ by cancellation.
\end{myproof}

\begin{lemma}\label{l:wedge.normal}
If $A$ and $B$ are normal monoids, so is $A\wedge B$
\end{lemma}

\begin{myproof}
The group completion of $A\wedge B$ is $(A\wedge B)_0=A_0\wedge
B_0$. If $a\wedge b\in A_0\wedge B_0$ is integral over $A\wedge
B$, then $(a\wedge b)^n=a^n\wedge b^n$ in $A\wedge B$, hence
$a^n\in A$, $b^n\in B$.  By normality, $a\in A$ and $b\in B$.
\end{myproof}

Recall from \cite[8.1]{chww} that a \emph{valuation monoid} is a cancellative
monoid $A$ such that for every non-zero $\alpha$ in the group completion $A_0$,
either $\alpha\in A$ or $\alpha^{-1}\in A$. Passing to units, we see that
$A^{\times}$ is a subgroup of the abelian group $A_0^\times$, and the {\it value
group} is the quotient $A_0^\times/A^\times$.  The value group is a totally
ordered abelian group ($x\geq y$ if and only if $x/y\in A$).  Following
\cite[8.3]{chww}, we call $A$ a {\em discrete valuation monoid}, or DV monoid
for short, if the value group is infinite cyclic.  In this case, a lifting
$\pi\in A$ of the positive generator of the value group generates the maximal
ideal $\frakm$ of $A$ and every $a\in A$ can be written $a=u\pi^n$ for some
$u\in A^{\times}$ and $n\ge0$.  Here $\pi$ is called a \emph{uniformizing
parameter} for $A$.

It is easy to see that valuation monoids are normal, and that noetherian
valuation monoids are discrete \cite[8.3.1]{chww}. We now show that
one-dimensional, noetherian normal monoids are DV monoids.

\begin{proposition}\label{p:dvm.normal}
Every noetherian, one-dimensional, normal monoid is a discrete valuation monoid
(and conversely).
\end{proposition}

\begin{myproof}
%
Suppose $A$ is a one-dimensional noetherian normal monoid, and choose a nonzero
$x$ in the maximal ideal $\frakm$.  By primary decomposition \ref{t:p.decomp},
$\sqrt{xA}$ must be $\frakm$ and (by Lemma \ref{rad.power}) there is an $n>0$
with $\frakm^n\subseteq xA$, $\frakm^{n-1}\not\subseteq xA$.  Choose $y\in
\frakm^{n-1}$ with $y\not\in xA$ and set $\pi=x/y\in A_0$.  Since
$\pi^{-1}\not\in A$ and $A$ is normal, $\pi^{-1}$ is not integral over $A$.  By
Lemma~\ref{l:cayley.int}, $\pi^{-1}\frakm\not\subseteq \frakm$; since
$\pi^{-1}\frakm\subseteq A$ by construction, we have $\pi^{-1}\frakm=A$, or
$\frakm=\pi A$.

Lemma~\ref{l:cayley.int} also implies that $\pi^{-1}I\not\subseteq I$ for every
ideal $I$. If $I\ne A$ then $I\subseteq\pi A$ so $\pi^{-1}I\subseteq A$. Since
$I=\pi^{-1}(\pi I)\subset\pi^{-1}I$, we have an ascending chain of ideals which
must terminate at $\pi^{-n}I=A$ for some $n$. Taking $I=aA$, this shows that
every element $a\in A$ can be written $u\pi^n$ for a unique $n\ge0$ and $u\in
A^\times$. Hence every element of $A_0$ can be written $u\pi^n$ for a unique
$n\in\Z$ and $u\in A^\times$, and the valuation map $\mathrm{ord}:A_0\ra
\Z\cup\{\infty\}$ defined by $\mathrm{ord}(u\pi^n)=n$ makes $A$ a discrete
valuation monoid.
\end{myproof}

\begin{corollary}\label{c:DVmonoid}
If $\frakp$ is a height one prime ideal of a noetherian normal monoid $A$ then
$A_{\frakp}$ is a discrete valuation monoid (DV monoid).
\end{corollary}

\begin{myproof}
The monoid $A_{\frakp}$ is one dimensional and normal by
Lemma~\ref{l:local.normal}.  Now use Proposition~\ref{p:dvm.normal}.
\end{myproof}


\begin{lemma}\label{l:prime.principal}
If $A$ is a noetherian normal monoid, and $\frakp$ is a prime ideal associated
to a principal ideal, then $\frakp$ has height one and $\frakp_{\frakp}$ is a
principal ideal of $A_{\frakp}$.
\end{lemma}

\begin{myproof}
Let $a\in A$ and $\frakp$ a prime ideal associated to $aA$ so that by
Proposition \ref{p:ap.ann}, $\frakp=(a:b)$ for some $b\in A\backslash aA$.  To
show $\frakp_{\frakp}\subseteq A_{\frakp}$ is principal, we may first localize
and assume that $A$ has maximal ideal $\frakp$.  Let $\frakp^{-1}=\{u\in A_0 \ |
\ u\frakp\subseteq A\}$.  Since $A\subseteq \frakp^{-1}$, we have
$\frakp\subseteq \frakp^{-1}\frakp\subseteq A$, and since $\frakp$ is maximal,
we must have $\frakp^{-1}\frakp=\frakp$ or $\frakp^{-1}\frakp=A$.

If $\frakp^{-1}\frakp=\frakp$, every element of $\frakp^{-1}$ must be integral
over $A$ by Lemma~\ref{l:cayley.int}.  Since $A$ is normal,
$\frakp^{-1}\subseteq A$, hence $\frakp^{-1}=A$ and $\frakp b\subseteq aA$
implies $b/a\in\frakp^{-1}=A$.  This is only the case if $b\in aA$, since $a$ is
not a unit, contradicting the assumption.  Therefore $\frakp^{-1}\frakp=A$ and
there exists $u\in\frakp^{-1}$ with $u\frakp=A$, namely $\frakp=u^{-1}A$.
\end{myproof}

To finish the section we show that any noetherian, normal monoid is the
intersection of its localizations at height one primes.  As with Corollary
\ref{c:DVmonoid}, this result parallels the situation in commutative rings.

\begin{lemma}\label{l:local.zero}
Let $A$ be a noetherian monoid and $I\subseteq A$ an ideal.
\begin{compactenum}
\item[i)] If $I$ is maximal among ideals of the form $(0:a)$, $a\in A$, then $I$
    is an associated prime of $0$.
\item[ii)] If $a\in A$, then $a=0$ in $A$ if and only if $a=0$ in $A_{\frakp}$
    for every prime $\frakp$ associated to $0$.
\end{compactenum}
\end{lemma}

\begin{myproof}
Suppose that $I=(0:a)$ is maximal, as in (i).  If $xy\in I$ but $y\not\in I$,
then $axy=0$ and $ay\ne0$.  Hence $I\subseteq I\cup Ax\subseteq (0:ay)$; by
maximality, $I=I\cup Ax$ and $x\in I$. Thus $I$ is prime; by
Proposition~\ref{p:ap.radquot}, $I$ is associated to $(0)$.

Suppose that $0\neq a\in A$, and set $I=(0:a)$.  By (i), $I\subseteq\frakp$ for
some associated prime $\frakp$.  But then $a\ne0$ in $A_{\frakp}$.
\end{myproof}

\begin{theorem}\label{t:intersect.Ap}
A noetherian normal monoid $A$ is the intersection of the $A_{\frakp}$ as
$\frakp$ runs over all height one primes of $A$.
\end{theorem}

\begin{myproof}
That $A$ is contained in the intersection is clear.  Now, suppose $a/b\in
A_0\backslash A$ so that $a\not\in bA$.  Any $\frakp\in\Ass(b)$ has
$\frakp_{\frakp}$ principal by Lemma~\ref{l:prime.principal}, hence height one,
and $a/b\not\in A_{\frakp}$ when $a\not\in bA_{\frakp}$.  Therefore to find an
associated prime $\frakp$ of $bA$ with $a\not\in bA_{\frakp}$ will complete the
proof.  But this is easy since $a\in bA_{\frakp}$ for every $\frakp\in\Ass(b)$
if and only if $a=0$ in $A_{\frakp}/bA_{\frakp}=(A/bA)_{\frakp}$ for every
$\frakp\in\Ass(b)$, which happens if and only if $a=0$ in $A/bA$ by Lemma
\ref{l:local.zero}, which happens if and only if $a\in bA$.  Since $a\neq 0$,
such
a prime must exist.
\end{myproof}

\chapter{Homological algebra}\label{c:homological}

In this chapter we examine the possibilities for a computable homological theory
on $A\set$.  To this end we first examine what it would mean to extend the
homological theory of abelian categories to $A\set$ using simplicial objects and
Quillen's homotopical algebra \cite{Quil}.  As per the Introduction, computing
homotopically is extremely difficult in categories whose objects do not have an
underlying (abelian) group structure since general simplicial objects will not
be fibrant.  To remedy the situation we extend the usual definition of
\emph{chain complex} to \emph{double-arrow complex} to explore as a replacement for
simplicial objects.

All results in this chapter will be proved for $A\set$ but it is clear that they
can be extended to more general categories.  \emph{All results of this chapter
apply to categories $\cC$ which}:
\begin{compactenum}
\item are concrete, namely every object of $\cC$ has an underlying set.  More
    precisely, there is a faithful functor $u:\cC\ra\Sets$.
\item have all equalizers and coequalizers.
\item have all pullbacks and pushouts.
\end{compactenum}
It is possible that (i) may be too strong of a condition and that all the
results of this section can be reworked without using the underlying set of the
objects in $\cC$; instead relying upon the existence of objects like the image of
a morphism.  We do not pursue this type of generalization any further.

\section{Preliminaries}\label{preliminaries}

We briefly review some definitions given earlier in the thesis.  Let $f:X\ra Y$
be a morphism of $A$-sets.  Then $f$ is \emph{admissible} when it is the
cokernel of a morphism, equivalently $X/\ker(f)\cong \im(f)$.  A sequence of
morphisms $X\xra{f}Y\xra{g}Z$ is said to be \emph{exact at $Y$} when
$\ker(g)=\im(f)$.  An exact sequence of the form $0\ra X'\ra X\ra X''\ra0$ is
called a \emph{short exact sequence} or \emph{s.e.s.} for brevity.  When every
morphism of a s.e.s.  is admissible, we say call it an \emph{admissible short
exact sequence} (\emph{a.s.e.s.}) of a.e.s. for brevity.  In this case, using
the notation as above, we have
\begin{compactenum}
\item $X'\ra X$ is an injection.
\item $X\ra X''$ is a surjection and $X/X'\cong X''$.
\end{compactenum}
The adjectives ``exact'' and ``admissible'' apply equally well to more general
sequences:
\[\cdots\ra X_{n+1}\ra X_n\ra X_{n-1}\ra\cdots\]
When it will cause no confusion we forgo the usage of ``short'' in s.e.s.  and
a.s.e.s., simply referring sequences being exact and/or admissible.  In this
case we will abbreviate admissible exact sequence by a.e.s..

There are many $A$-sets which cannot be written as the quotient of a free
$A$-set by an $A$-subset, e.g. see \ref{e:nonsplitting} .  This motivates the
following definitions.  The equalizer (resp.  coequalizer) of two morphisms is
the limit (resp.  colimit) of the following diagram:
\[\cdot\rra\cdot\]
The \emph{equalizer} of $f,g:X\ra Y$ is the morphism $i:Z\ra X$ where
\[Z=\{ x\in X \ | \ f(x)=g(x)\}\subseteq X\]
and the \emph{coequalizer} is
\[\coeq(f,g)=Y/(f(x)\sim g(x) \ | \ x\in X).\]
To avoid confusion, we may write $\eq(f,g)$ (respectively $\coeq(f,g)$) for the
equalizer (respectively coequalizer) of $f$ and $g$.  In general, we say that a
diagram
\[Z\xra{i}X\da{f}{g}Y\]
\emph{commutes} when $fi=gi$.  In this case $i$ \emph{equalizes} $f$
and $g$ (but $i$ \emph{need not be the equalizer}!).  Similarly, a diagram
\[X\da{f}{g}Y\xra{h}Z\]
\emph{commutes} when $hf=hg$ and in this case, $h$ \emph{coequalizes} $f$ and
$g$ (but $h$ \emph{need not be the coequalizer}!).  For more on limits and
colimits, see \cite{WH}.  It is shown in \cite{cls} that $A\set$ contains all
small limits and colimits so that we need not worry about the existence of such
objects in this paper.

\begin{proposition}\label{p:regularmorphism} \begin{compactenum} \item Every
surjection of $A$-sets is a coequalizer. \item Every morphism of $A$-sets
$g:X\ra Y$ can be factored $g=i\circ f$ where     $f$ is a surjection and $i$
is injective. \end{compactenum} \end{proposition}

\begin{myproof}
i) Let $f:X\ra Y$ be a surjection and consider the pullback:
\[\SelectTips{eu}{12}\xymatrix{
    Z \ar[d]_{p_2} \ar[r]^{p_1}_{\ulcorner} & X \ar[d] \\
X \ar[r]  & Y
} \]
Then $Z$ is an $A$-set and $f$ is the coequalizer
$Z\overset{p_1}{\underset{p_2}\rightrightarrows} X\xra{f} Y$.\\
ii)  The map $f:X\ra g(X)$ defined by $f(x)=g(x)$ is a surjection.  If 
$i:g(X)\hookrightarrow Y$ is the
inclusion, clearly
$g=i\circ f$.
\qed
\end{myproof}

By Proposition \ref{p:regularmorphism} every $A$-set $X$ can be realized as the
quotient of a free $A$-set $F$ as in $R\rra F\ra X$. Continuing this diagram to
the left provides the definition for ``double-arrow complex'' which extends the
usual definition of ``complex'' in abelian categories (see Section
\ref{dacomplexes}).  Then we are able to construct projective resolutions for
general $A$-sets allowing a more computable theory of derived functors.

Readers familiar with the language of regular and/or (Barr) exact categories
will be very familiar with the previous concepts. Though the concepts defining
these categories are closely related to the properties we require, it is not
clear that regular categories provide the right context.

\section{Functors and exactness}\label{s:functors}

With short exact sequences in hand we may now turn to functors.  For a functor $F:\cC\ra\cD$ and an object $X\in \cC$, we sometimes write $FX$ instead of $F(X)$ when the role of $F$ as a functor is clear from the context.  When $\cC$, $\cD$ are pointed categories, $F$ is \emph{based} when $F(0_{\cC})=0_{\cD}$. The definition of admissible morphism in $A\set$ may be defined in any pointed category with images and cokernels: a morphism $f:X\ra X'$ is \emph{admissible} when $X\ra\im(f)$ is a cokernel.  We say that $F$ \emph{preserves admissibles} if $f\in\Hom_{\cC}(X,X')$ is admissible means $F(f)$ is also.

To clarify why we are only concerned with admissible exact sequences, let $A=\F_1[x]$ and consider the following exact, but not admissible, sequence of $A$-sets:
\[0\ra \bigvee_{i=0}^{\infty} A_i\ra A\ra 0,\]
where each generator $1_i$ of $A_i$ is mapped to $x_i$.  This sequence is exact even though the cardinalities of their generator(s) differ.  Exactness alone does not say very much about objects in the sequence.

Let $\cC$ be a pointed, concrete category so that images and cokernels exist.
Then we may use the adjectives admissible, exact and short in reference to a
sequence $\cdots\ra X_{n+1}\ra X_n\ra X_{n-1}\ra\cdots$ in $\cC$.

\begin{definition}\label{d:lrexact}
Let $\cC$ be a concrete, pointed category.
A based functor $F:\cC\ra A\set$ is said to be \emph{left} (resp.$\,$\emph{right}) \emph{exact} when $0\ra X'\ra X\ra X''\ra 0$ is an \emph{admissible} exact sequence in $\cC$ implies
\[0\ra F(X')\ra F(X)\ra F(X'') \ \ \mathrm{(resp.}\ F(X')\ra F(X)\ra F(X'')\ra 0\mathrm{)}\]
is an \emph{admissible} exact sequence.  
In $A\set$ we are only ever concerned with admissible short exact sequences, hence when $\cC= A\set$ we \emph{add the additional requirement} that $F$ preserve admissibles.
A functor is \emph{exact} when it is both left and right exact.
\end{definition}

\begin{remark}
For an $A$-set $Y$, we note that $\Hom(Y,-):A\set\ra\Sets_*$, is ``left exact'' in the following sense.  If $0\ra X'\xra{f} X\xra{g} X''\ra 0$ is an a.e.s. of $A$-sets, then $0\xra{}\Hom(Y,X')\xra{f_*}\Hom(Y,X)\xra{g_*}\Hom(Y,X'')$ is exact and $f_*$ is injective.  However, $g_*$ need not be admissible.  Note also that $\Hom(-,Y):A\set\ra\Sets_*$ is contravariant left exact in the same sense.
\end{remark}

This definition for left and right exactness is strictly for functors from
$A\set$ to itself and we will see below that it does not conflict with the
definition for functors taking values in an abelian category.  Thus, although it
is tempting to leave out the ``preserves admissibles'' property and define left
(and right) exact functors strictly in terms of the exactness property of
morphisms, this will only cause more repetition in the future.  Within $A\set$
we will only ever be concerned with admissible exact sequences and we define
exact functors accordingly.

\begin{remark}
A morphism $f:Y\ra Z$ is admissible when $Y\ra\im(f)\ra Z$ is an admissible
sequence, i.e. $Y\ra\im(f)$ is a cokernel, hence colimit.  Therefore,
any concrete functor preserving colimits, or even cokernels, and injectives will 
preserve the admissible property.  In particular, any functor which is a left 
adjoint which preserves injectives preserves admissibles.
\end{remark}

Since $-\ot X$, $X\ot -$ do not preserve injective maps and $\Hom_A(X,-)$, 
$\Hom_A(-,X)$ do not preserve colimits, we should not expect any of these 
functor to preserve admissible morphisms.

\begin{example}
To see $-\ot X$ does not preserve admissible morphisms, let $A=\F_1[t]$, 
$X=A\vee A'/(t\sim t')$ and consider the admissible inclusion $A\xra{t} A$.  
Applying $-\ot X$, the morphism $X\xra{t} X$ is no longer admissible since both 
generators $1, 1'\mapsto t(=t')$.
\end{example}

\begin{example}\label{e:homnotadm}
Let $A=\F_1[t]$ be the free monoid in one variable, $X=A\vee A$ and $p:X\ra A$ the admissible morphism defined by $p=1\vee 0$ (see Remark \ref{r:wedge.morphism.notation} about notation).  
Define $\alpha,\beta:X\ra X$ by $\alpha=1\vee t$ and $\beta=1\vee t^2$.  
Then $\alpha\neq \beta\in\Hom(X,X)$, but $p_*(\alpha)=p_*(\beta)=1\vee 0$ in $\Hom(X,A)$ so that $p_*:\Hom(X,X)\ra\Hom(X,A)$ is not admissible.
That is, the fiber over the \emph{nonzero} morphism $1\vee 0$ in $\Hom(X,A)$ contains multiple elements (at least $\alpha$ and $\beta$).
\end{example}

A similar example can be constructed to show $\Hom(-,Z)$ also does not preserve admissibles.  
Furthermore, $\Hom(-,Z)$ is a contravariant functor and we do not discuss contravariant theories in this thesis.  
See Section \ref{extensions} and Chapter \ref{geometry.monoids} for more information on cohomology.



\subsection{$k$-realization}\label{krealization}

Let $k$ be a commutative ring with identity and $A$ a monoid.  The \emph{$k$-realization} functor $k:\mon\ra k\text{\bf{-algebras}}$ assigns to every monoid $A$ the free $k$-module $k[A]$ whose basis is the nonzero elements of $A$ with multiplication induced by the multiplication of $A$.  This is easily extended to the $k$-realization of $A$-sets, $k: A\set\ra k[A]\mod$, for if $X$ is an $A$-set, let $k[X]$ be the free $k[A]$-module whose basis is the set of nonzero elements of $X$ and with $k[A]$-action given by the $A$-action on $X$ (together with the action of $k$).  The following proposition is clear.

\begin{proposition}
The functor $k: A\set\ra k[A]\mod$ is exact in the sense that it carries admissible exact sequences to exact sequences (in the usual sense of abelian categories).
\end{proposition}

In the other direction is the \emph{forgetful}, or \emph{underlying monoid}
(resp.  \emph{$A$-set}), functor $U:\mathbf{Rings}\ra\mon$ (resp.
$U:k[A]\mod\ra A\set$), that simply forgets the addition of $k$.  For example,
using $k=\Z$ and $A=\F_1$, we have $U(k[\F_1])=U(\Z)$ is the free monoid on
countably many generators (the primes) with coefficients in $\{0,\pm 1\}$.

The functors $k$ and $U$ (for both $\mon$ and $ A\set$) form an adjunction with
$k$ left adjoint to $U$.  That is, the map
$\Hom_{\mon}(A,U(R))\ra\Hom_{k\text{\bf{-algebras}}}(k[A], R)$ extends a monoid
morphism $f$ to addition, $f\left(\sum_i r_ia_i\right)=\sum_ir_if(a_i)$.  In the
opposite direction, a ring map $k[A]\ra R$ simply forgets it is a homomorphism
with respect to addition.


\begin{example}
Let $C$ be a monoid and let $A,B$ be $C$-monoids.  Then $k[A\otimes_C B]\cong
k[A]\otimes_{k[C]} k[B]$.  When $C=\F_1$ and $B=\F_1[x_1,\ldots,x_n]$ is the
free monoid in $n$ variables, $A\otimes B=A[x_1,\ldots,x_n]$ is the free
$A$-monoid on $n$ variables and $k[A[x_1,\ldots,x_n]]=k[A][x_1,\ldots,x_n]$ is
the free polynomial algebra on $n$ variables with coefficients in $k[A]$.
\end{example}


\section{Projective $A$-sets}\label{s:projectives}

Free objects are special in that every morphism whose domain is a free object is
uniquely determined by where it maps the generators.  Thus $\Hom_A(F,X)$ is
plentiful when $F, X$ are $A$-sets and $F$ is free.  Projectives are a larger
class of objects that which generalize this idea.  In any category an object $P$
is \emph{projective} when it satisfies the following universal lifting property.
\[\SelectTips{eu}{12}\xymatrix{
    & P \ar[d]^g \ar@{-->}[dl]_{\exists\,\varphi} \\
 X \ar[r]^f  & Y
} \]
Given any epimorphism $f:X\ra Y$ and any
morphism $g:P\ra Y$, there exists $\varphi:P\ra X$ such that $g=f\varphi$ (here
$\varphi$ is called a \emph{lifting} map).  Equivalently, if $X\ra Y$ is a
epimorphism, then so is $\Hom(P,X)\ra\Hom(P,Y)$.

A projective $A$-set may also be realized as the retract of a free object.
Given a morphism $f:X\ra Y$, we say that $Y$ is a \emph{retract} of $X$ when
there exists a morphism $\sigma:Y\ra X$ such that $\id_Y=f\sigma$.  In this
case, we say that $\sigma$ is a \emph{section} of $f$.  Any retract of a 
projective $A$-set is projective.

\begin{lemma}
Let $A$ be a monoid, $e\in A$ idempotent and $K$ a set.  Then $eA$, $A[K]$ and
$eA[K]$ are projective $A$-sets.
\end{lemma}

\begin{myproof}
Any morphism $f:A[K]\ra Y$ is determined by the (set-theoretic) function $f|_K:K\ra Y$.  To lift any surjective $A$-set map $X\ra Y$, simply lift $f|_K$. In the latter case, note that $eA$ is a retract of $A$.  Namely, the surjection $e:A\ra eA$ defined by $1\mapsto e$ has a section $\sigma:eA\ra A$  which is the inclusion $\sigma(e)=e$.  Since the retract of a projective $A$-set is projective, the result follows.
\qed
\end{myproof}

\begin{remark}
(Warning) If $e$ is an idempotent element of $A$, the surjection $A\ra eA$ may
not be admissible so we do not necessarily have $A=eA\vee Q$ for some
$A$-set $Q$.
\end{remark}

If $P=\vee_i P_i$, then $P$ is projective if and only if every $P_i$ is projective.  This is easily deduced from the fact that $\Hom(P,X)=\prod_i\Hom(P_i,X)$ for any $A$-set $X$.  We can now classify projective $A$-sets.

\begin{theorem}\label{t:projectives}
Let $A$ be a monoid and $P$ a projective $A$-set.  Then $P=\vee_{i\in I}Ae_{i}$
where $e_{i}\in A$ is idempotent for each $i$.
\end{theorem}

\begin{myproof}
Consider the diagram:
\[\SelectTips{eu}{12}\xymatrix{
    & P \ar @{-->} [dl]_{\exists\,\varphi} \ar @{=} [d] & \\
    A[P] \ar[r]^{\ \ \ f}  & P \ar[r] & 0
} \]
where $f[x]=x$ (see Example \ref{e:cotriple.notation} about the notation
$A[P]$).  Then $f$ is a surjection and the section $\varphi$  monic since $f\circ\varphi=\mathrm{id}_{P}$.  Lemma~\ref{p:splitting} shows
$P=\vee_{x\in P}P_{x}$ where $P_x=\varphi^{-1}(\varphi(P)\cap A_{x})$ and the containment $P_x=f(\varphi(P)\cap A_x)\subseteq f(A_x)$ implies $P_x$ is generated by $f(1)$, $1\in A_x$.  The restrictions $\varphi:P_x\ra A_x$ are themselves sections of $f$ hence $P_x\cong A(\varphi\circ f)(1)=A\varphi(x)$.  Finally, the equation $\varphi(x)=(\varphi\circ f)(\varphi(x))=\varphi(\varphi(x)x)=\varphi(x)^2$ shows $\varphi(x)$ is idempotent.
\qed
\end{myproof}

\begin{corollary}
Let $k$ be a commutative ring and $P$ a projective $A$-set.  Then $k[P]$ is a
projective $k[A]$-module.
\label{c:k.proj}
\end{corollary}

\begin{myproof}
Given a diagram
\[ \SelectTips{eu}{12}\xymatrix{
    & k[P] \ar[d] & \\
 N \ar[r]^f  & M \ar[r] & 0
} \]
apply the forgetful functor $U:k[A]\mod\ra A\set$ and use the lifting property
of $P$ to obtain the following diagram:
\[ \SelectTips{eu}{12}\xymatrix{
    P\ar[r]\ar@{-->}[d]_{\exists\,\varphi}  & Uk[P] \ar[d] & \\
 UN \ar[r]^{Uf}  & UM \ar[r] & 0
} \]
Then $k[\varphi]:k[P]\ra N$ provides a lifting for the original diagram.
\qed
\end{myproof}

Of course, the usual splitting result holds.

\begin{proposition}\label{p:projectivessplit}
When $P$ is a projective $A$-set, any a.e.s. $0\ra X'\ra X\ra P\ra 0$ splits.
\end{proposition}

\begin{myproof}
Lemma \ref{l:splitting} shows that the sequence splits when there is a section
$\sigma:P\ra X$.  Since $P$ is projective and $X\ra P$ is surjective, we are
done.
\qed
\end{myproof}

\begin{proposition}\label{p:pushforwardprojective}
Let $f:A\ra B$ be a morphism of monoids.  If $P$ is a projective $A$-set, then
$f_*P=B\otimes_A P$ is a projective $B$-set.
\end{proposition}

\begin{myproof}
If $e\in A$ is idempotent, then $f(e)$ is idempotent in $B$.  Let $P=\vee_i
Ae_i$ be a projective $A$-set with each $e_i$ idempotent. Then
\[B\otimes_A P=B\otimes_A (\vee_i Ae_i)\cong\vee_i (B\otimes_A Ae_i)\cong \vee_i
Bf(e_i)\]
which is also projective.
\qed
\end{myproof}

\subsection{Admissibly projective $A$-sets}\label{ss:relativelyprojective}

Proposition \ref{p:projectivessplit} shows the standard result that
projective $A$-sets split admissible exact sequences.  They are not alone.

\begin{definition}
An $A$-set $Q$ is \emph{admissibly projective} if it satisfies the
lifting property
\[\SelectTips{eu}{12}\xymatrix{
    & Q \ar[d]^g \ar@{-->}[dl]_{\exists\,\varphi} & \\
 X \ar[r]_f  & Y \ar[r] & 0
} \]
whenever $f$ is an \emph{admissible} surjection.  From Lemma \ref{l:splitting}
it is immediate that every a.e.s. $0\ra X\ra Y\ra Q\ra 0$ with $Q$ relatively
projective, splits.  The converse is shown in Proposition \ref{p:split.relproj}.
\end{definition}

\begin{proposition}\label{p:split.relproj}
Let $Q$ be an $A$-set such that every a.e.s. $0\ra X'\ra X''\ra Q\ra 0$ splits.  
Then $Q$ is admissibly projective.
\end{proposition}

\begin{myproof}
Consider the diagram
\[\SelectTips{eu}{12}\xymatrix{
    & & P \ar@{-->}[d]_{p_1}\ar@{-->}[r]^{p_2} & Q \ar[d]^g & \\
    0\ar[r] & K \ar[r] & X \ar[r]_f  & Y \ar[r] & 0
} \]
where $f$ is an admissible surjection, $K=\ker(f)$ and $P=X\times_Y Q$ is the 
pullback of $f$ and $g$.  The projection map $p_2:P\ra Q$ is simply 
$(x,q)\mapsto q$ and therefore, if $(x,q)\neq(x',q')$ and $p_2(x,q)=p_2(x',q')$, then $q=q'$ and by commutativity, $f(x)=f(x')$.  Since $f$ is admissible, we must have $x,x'\in K$, so
$q,q'\in\ker(p_2)$, showing that $p_2$ is also admissible.  By assumption, the 
a.e.s.  $0\ra \ker(p_2)\ra P\xra{p_2} Q\ra 0$ splits, say with section  
$\varphi:Q\ra P$, and the composition $p_1\varphi:Q\ra X$ provides the lifting 
map since $fp_1\varphi=gp_2\varphi=g$.
\qed
\end{myproof}

\begin{example}\label{e:relativelyprojective} 
Let $A=\F_1[x]$ be the free monoid in one variable and set $Q=(A\vee A)/(x^n\vee 0\sim 0\vee x^n)$ for some $n>1$. For any a.e.s.  
\[0\ra Z\ra Y\xra{f} Q\ra0\]  
consider the pointed set theoretic section $\sigma:Q\ra Y$ sending every nonzero element of $Q$ to the unique element in its fiber.  Since $Y=Z\vee \sigma(Q)$ \emph{as pointed sets} and the image of $Z$ is an $A$-subset of $Y$, if $\sigma(Q)$ is also an $A$-subset of $Y$ we will have $Y\cong Z\vee\sigma(Q)$ as $A$-sets so that the sequence splits.  But $a\sigma(x)\in Z$ if and only if $ax=0$ in $Q\cong Y/Z$.  Since $Q$ is torsion free, we have $Y\cong Z\vee Q$ as an $A$-set.

Now, $Q$ is \emph{not} projective since there is no lifting of the identity map
$Q\ra Q$ over the surjection $A\vee A\ra Q$ which maps the two generators of
$A\vee A$ to the two generators of $Q$.
\end{example}

Recall that an element $x$ of an $A$-set $X$ is a torsion element when $ax=0$
for some nonzero $a\in A$.  Note that this notion of torsion differs from that
in group theory.  For example, the pointed, cyclic group
$G_+=\{0,1,x,\ldots,x^n\}$ is a torsion free $\F_1[x]$-set.

\begin{proposition}\label{p:relprojtorsionfree}
Any torsion free $A$-set is admissibly projective.
\end{proposition}

\begin{myproof}
Let $Q$ be a torsion free $A$-set.  By Proposition \ref{p:split.relproj}, it suffices to show that any a.e.s. $0\ra Z\ra Y\xra{f} Q\ra 0$ splits. 
Let $\sigma:Q\ra Y$ be the pointed set theoretic section of $f$.  
When $Q$ is torsion free, $0\neq x\in Q$ and $0\neq a\in A$,  we have $f(a\sigma(x))=af(\sigma(x))=ax\neq 0$ so that $a\sigma(x)$ is the unique element in the fiber of $ax$, i.e. $\sigma(ax)$. 
Hence $\sigma$ is an $A$-set morphism and $Q$ is an $A$-subset of $Y$.
\qed
\end{myproof}

Note that Proposition \ref{p:relprojtorsionfree} is also a consequence of Proposition \ref{p:extensions}.  Also, the converse to Proposition \ref{p:relprojtorsionfree} is certainly not true. Indeed, every projective $A$-set is admissibly projective and when $A$ has zero divisors, any free $A$-set provides a counterexample.

\subsection{Rank}\label{s:rank}

Let $G$ be an abelian group and $G_+$ the associated monoid (see Example 
\ref{e:pointedgroup}).  As we noted in Example \ref{e:groupsets}, every 
$G_+$-set can be written as $X=\vee_{i\in I}X_i$, $I$ an indexing set, where 
each $X_i=(G/H_i)_+$ is a (pointed) quotient of $G$ by a subgroup $H_i$.  In 
this situation the cardinality of $I$ is well defined and we may define the 
\emph{rank} of a $G_+$-set as the cardinality of $I$.  Of course in classical 
group theory, the $X_i$ are called the \emph{orbits} of $X$ and $|I|=|X/G|$ 
where $X/G$ is the quotient $G$-set obtained by considering ``$X$ modulo the 
action of $G$.''  In the case of monoids, the basepoint will always be its own 
orbit.  If $G$ is finite, Burnside's Lemma provides the formula
\[|X/G|=\frac{1}{|G|}\sum_{g\in G}|X^g|\]
where $X^g=\{ x\in X \ | \ gx=x\}$.  This formula can be useful when $G$ is a
\emph{finite} abelian group and may not be as helpful in the study of finitely
generated abelian groups.

The rank of an $A$-set may now be defined analogously to that of modules over 
$R$, a commutative ring.  Recall that the \emph{group completion} of a 
cancellative monoid $A$ is the monoid $A_0$ obtained by inverting all non-zero 
elements of $A$.  When $A$ is non-cancellative, then $ab=ac$ for some $a,b,c\in 
A$ and $b/1=c/1\in A[\frac{1}{a}]$.  Borrowing notation from commutative algebra 
again, we may use the notation $\Quot(A)$ for $A_0$.  If $\frakp\in\MSpec(A)$, 
then $A/\frakp$ has no zero-divisors but $A/\frakp\ra\Quot(A/\frakp)$ need not 
be injective.  However, the monoid $G(\frakp)=\Quot(A/\frakp)$ is nonetheless a 
pointed (abelian) group called the \emph{residue group at} $\frakp$.  

With these definitions in hand, we are able to define the rank of a general
$A$-set.  Let $A$ be a monoid and $X$ a finitely generated $A$-set.  The
\emph{rank} of $X$ over $A$ is the function $\rank:\MSpec(A)\ra\N\cup\{\infty\}$
defined by
\[\rank_\frakp X=\rank_{G(\frakp)} (X\otimes_A G(\frakp))\]
where $G(\frakp)=\Quot(A/\frakp)$ is the residue group at the prime $\frakp$.

\begin{lemma}
If $e$ is an idempotent element of $A$, then
\begin{equation*}
    \rank_{\frakp}(eA) = \left\{
\begin{array}{rl}
0 & e\in\frakp \\
1 & e\not\in\frakp
\end{array} \right.
\end{equation*}
\label{l:rank.idem}
\end{lemma}

\begin{myproof}
Since $eA\ot G(\frakp)=eA\ot (A/\frakp)_{\frakp}\cong (eA/e\frakp)_{\frakp}$, 
the prime $\frakp$ contains $e$ if and only if $eA\ot G(\frakp)=0$.  Moreover, 
if $\frakp$ does not contain $e$, then $e=1$ in $G(\frakp)$ so that $eA\ot 
G(\frakp)\cong G(\frakp)$.
\qed
\end{myproof}

\begin{remark}
Unlike commutative rings, non-trivial idempotent elements do not disconnect $\MSpec(A)$.  This should be clear since all monoids are local and a disconnected spectrum requires at least two maximal ideals. It is a standard result in commutative ring theory that for a commutative ring $R$, $\Spec(R)$ is connected if and only if $R$ contains only the trivial idempotents.  When $R$ contains non-trivial idempotents and $\Spec(R)$ is disconnected, the underlying monoid has $\MSpec(U(R))$ connected.  It is also known that when $P$ is a projective $R$-module, the rank function for $P$ is locally constant.  However, since the connected components of $\Spec(R)$ come together to form the connected topological space $\MSpec(U(R))$, we should not expect the rank function for the $U(R)$-set $U(P)$ to be locally (in fact, globally) constant.

A simple example: $A=\F_1[x,y]/(x^2=x,y^2=y)$ consists of only idempotents and 
has $\MSpec(A)=\{(0), (x), (y), (x,y)\}$.  Then the projective $A$-set 
$X=Ax=\{0,x,xy\}$ has $\rank_{(0)}X=\rank_{(y)}X=1$ and 
$\rank_{(x)}X=\rank_{(x,y)}X=0$.

\end{remark}

Let $A$ be a monoid and $E\subseteq A$ the submonoid of idempotent elements of
$A$.  Define a relation $\leq$ on $E$ by $e \leq f$ whenever $ef=e$.  Then
$\leq$ is a partial order, and $E$ is a poset, since:
\begin{compactenum}
\item $e\leq e$ since $e^2=e$ (reflexive)
\item $e\leq f$ and $f\leq e$ implies $e=f$ since $e=ef=fe=f$ (antisymmetric)
\item $e\leq f$ and $f\leq h$ implies $e\leq h$ since $eh=(ef)h=e(fh)=ef=e$
    (transitive)
\end{compactenum}
In fact, $E$ is a semi-lattice where the meet, or greatest lower bound, of $e$
and $f$ is $e\land f=ef$.  Also, $E$ has 1 as its greatest element and 0 as its
least element.

\begin{lemma}\label{l:invertidempotents}
Let $e, f\neq 1$ be idempotent elements of $A$ and $\frakp$ the unique 
prime ideal maximal with respect to the condition $f\not\in\frakp$.  Then 
$e\not\in\frakp$ if and only if $e\geq f$.
\end{lemma}

\begin{myproof}
First, assume $e\geq f$ and suppose that $e\in\frakp$.  It immediately follows that $ef=f\in\frakp$, a contradiction.  

For the converse, we prove the contrapositive, namely if $e\not\geq f$ then $e\in\frakp$.  Notice that $\frakp$ my be realized as the contraction of the maximal ideal in the map $A\ra A_f$.  Also note that $e\in\frakp$ when $e$ is not a unit of $A_f$, equivalently $f\neq ae$ for any $a\in A$.  To make matters worse, suppose (in order to obtain a contradiction) that $f=ae$ for some $a\in A$.  Then $ae$ is an idempotent and $ef=e(ae)=ae=f$ implies $f\leq e$.  Hence, we must have $f\neq ae$ so that $e\in\frakp$.
\qed
\end{myproof}

\begin{theorem}\label{t:rankprojiso}
Let $P,Q$ be finitely generated projective $A$-sets.  Then $P\cong Q$  if and
only if $\rank_{\frakp}(P)=\rank_{\frakp}(Q)$ for every $\frakp\in\MSpec(A)$.
\end{theorem}

\begin{myproof}
If $P\cong Q$, it is clear that $\rank_{\frakp}(P)=\rank_{\frakp}(Q)$ for every
$\frakp$.

Now assume $\rank_{\frakp}P=\rank_{\frakp}Q$ for every prime $\frakp$.  By 
Theorem \ref{t:projectives}, $P=\bigvee_{i=1}^n Ae_i$ and we will show that for each distinct idempotent generator $e$, there is a prime $\frakp\subseteq A$ such that $\rank_{\frakp}(P)$ is the number of summands of $P$ generated by all $e_i\geq e$.  We may assume each $e_i$ is distinct and ordered so that $j>i$ means $e_i\not\leq e_j$.  For each $1\leq i\leq n$, let $\frakp_i$ be the prime ideal of $A$ maximal with respect to the condition $e_i\not\in\frakp_i$ so that by Lemma \ref{l:invertidempotents}, $j>i$ implies $e_j\in\frakp_i$ and hence, $\rank_{\frakp_i}(Ae_j)=0$.

Using Lemma \ref{l:rank.idem}, $k_1=\rank_{\frakp_1}(P)=\rank_{\frakp_1}(Q)$ is the number of $Ae_1$ summands in $P$ (and hence $Q$).  Let $P_1=P$ and for $i>1$, inductively define $P_i=P_{i-1}/(\bigvee_{1}^{k_i}Ae_i)$.  Making similar definitions for $Q_i$, we see that $k_i=\rank_{\frakp_i}(P_i)=\rank_{\frakp_i}(Q_i)$ is the number of $Ae_i$ summands in both $P$ and $Q$.  Therefore, $P\cong Q$.
\qed
\end{myproof}

\subsection{Observations on $K_0, G_0$ and $K_1$}\label{lowerkgroups}

With a good understanding of projective sets in hand, we show a few basic
$K$-theory results if only to become more familiar with $A\set$.  Given
\emph{any} (not necessarily commutative or pointed) monoid $M$, we construct the
\emph{group completion} of $M$ as follows.  First consider the cartesian
product $M\times M$ together with coordinate-wise addition:
\[(m_1,m_2)+(m'_1,m'_2)=(m_1+m'_1,m_2+m'_2).\]
Next define an equivalence relation $\sim$ on $M\times M$ by declaring
$(m_1,m_2)\sim(m'_1,m'_2)$ whenever there exists $n\in M$ such that
$m_1+m'_2+n=m_2+m'_1+n$.  This equivalence relation is compatible with the
additive structure.  It is easy to see that elements of the form $(m,m)$ are
identity elements and $(m,m')$ is the inverse of $(m',m)$.

The functor $K_0: \mon\ra\Ab$ assigns to each monoid $A$ the abelian group $K_0(A)$, called the \emph{Grothendieck group}, defined to be the group completion of the monoid whose elements are the isomorphism classes of finitely generated, projective $A$-sets, with binary operation given by the wedge sum $\vee$.

Theorem \ref{t:rankprojiso} shows that isomorphism classes of projective $A$-sets are completely determined by the number of summands generated by each idempotent of $A$.  The following proposition shows that $K_0(A)$ also has a multiplicative structure given by $\otimes$ so that it is actually a commutative ring with identity.

\begin{proposition}\label{p:tensorprojectives}
If $P$ and $Q$ are projective $A$-sets, so is $P\otimes_A Q$.
\end{proposition}

\begin{myproof}
Note that when $e,f\in A$ are idempotent their product $ef$ is also and $Ae\ot Af\cong Aef$.   Using Theorem \ref{t:projectives} and Proposition \ref{p:tensorproperties} we have
\[P\otimes Q=\bigvee_{i\in I}Ae_i\otimes\bigvee_{j\in J}Af_j\cong
\bigvee_{i,j\in I\times J}Ae_i\otimes Af_j\cong \bigvee_{i,j\in I\times
J}Ae_if_j\]
which proves the claim.
\qed
\end{myproof}

The next theorem summarizes the previous remarks.

\begin{theorem}
If $A$ is a monoid, then $K_0(A)$ is a commutative ring with addition given by
$\vee$ and multiplication given by $\otimes$.  The additive identity is the
trivial projective set $0$ and the multiplicative identity is the free rank one
set.  Moreover, there is a ring isomorphism $K_0(A)\cong\Z[E]$ where
$E=E_A$ is the submonoid of $A$ generated by all the idempotent elements.
\end{theorem}

\begin{myproof}
Recall that $\Z[E]$ is the monoid ring of $E$.   By Theorems \ref{t:projectives} and \ref{t:rankprojiso}, every projective $A$-set $P=\vee_{i\in I}Ae_i$ is determined, up to isomorphism, by the $e_i$.  Therefore, the $A$-sets $Ae$, where $e\in A$ is idempotent, form a generating set for $K_0(A)$ when considered as an abelian group. The function $f:K_0(A)\ra\Z[E]$ defined on the generators by $Ae\mapsto e$ extends linearly to a group homomorphism $\vee_{i\in I}Ae_i\mapsto\sum_{i\in I}e_i$.  It is clearly multiplicative, hence a ring homomorphism, and surjective.  For injectivity, after combining coefficients in 
\[f(P)=f(\vee_{i\in I}Ae_i)=\sum_{i\in I}e_i,\] 
we have $f(P)=0$ if and only if $e_i=0$ for all $i\in I$.
\qed
\end{myproof}

\begin{example}
The previous theorem implies that every monoid $A$ with only the trivial
idempotents (0 and 1) has $K_0(A)=\Z$.  In fact, in this case, every projective
$A$-set is free.  For instance, every $\F_1$-set $X$ is free and is completely
determined by its cardinality $|X|=\rank_{(0)}(X)+1$.
\end{example}

Let $f:A\ra B$ be a monoid morphism.  Recall that Proposition
\ref{p:pushforwardprojective} shows that the push-forward $f_*P$ is a projective
$B$-set when $P$ is a projective $A$-set.  It is then clear that $f$ induces a
ring homomorphism $f_*:K_0(A)\ra K_0(B)$.  Since we have the isomorphism
$K_0(A)\cong\Z[E_A]$ for every monoid, the following proposition is clear.

\begin{proposition}
Let $f:A\ra B$ be a monoid morphism.   If $f|_{E_A}:E_A\ra E_B$ is an
isomorphism, then so is $f_*:K_{0}(A)\ra K_{0}(B)$.
\end{proposition}

Perhaps more interesting than $K_0$ is $G_0$.  The definition of $K_0$ for 
finitely generated $A$-sets is extended from that of abelian categories (see 
\cite[II.6.2]{WK}).  We define the \emph{Grothendieck group} 
$G_0(A)$ of $A\set$ to be the abelian group having one generator $[X]$ for each 
finitely generated $A$-set $X$, modulo one relation $[X]=[X']+[X'']$ for every 
\emph{admissible} exact sequence $0\ra X'\ra X\ra X''\ra 0$ (see \cite[Section 5]{cls}).  
As with abelian categories, we immediately having the following identities:
\begin{compactenum}
    \item $0\ra X\ra X\ra 0$ implies $[0]=0$, that is, the generator given by
        the trivial $A$-set is the (additive) identity of $G_0(A)$.
    \item If $X\cong X'$, then $[X]=[X']$ from the a.e.s. $0\ra X\ra X'\ra 0$.
    \item $[X\vee X']=[X]+[X']$ from the a.e.s. $0\ra X\ra X\vee X'\ra X'\ra 0$.
\end{compactenum}

\begin{example}\label{e:k0.f1}
When $A=\F_1$, any finitely generated $\F_1$-set is simply a finite (pointed) set and two finite sets $X,X'$ are isomorphic if and only if $|X|=|X'|$.  The a.e.s. $0\ra \F_1\ra \vee_{1}^{n}\F_1\ra\vee_{1}^{n-1}\F_1\ra 0$ shows that $[\F_1]$ generates $G_{0}(\F_1)$.  Consequently, $G_{0}(\F_1)\cong \Z$.
\end{example}

\begin{example}(Burnside Ring)
Let $A=G_+$ be the pointed cyclic group of order $n$.  Every $G$-set is 
isomorphic to a wedge sum of (pointed) of cosets $(G/H)_+$ and two cosets $G/H$, 
$G/H'$ are isomorphic if and only if $H=H'$.  Thus, the relations shown above 
encompass all the relations of $G_0(A)$.  If $H_1,\ldots,H_m$ lists all possible
subgroups of $G$, including $G$ itself, we have 
$G_0(A)=\oplus_{i=1}^{m}\Z[G/H_i]$.  In particular, when $n=p$ is prime, 
$G_0(A)=\Z$.
\end{example}

\begin{example}
Let $A=\F_1[t]$ be the free monoid in one variable and $X$ a finitely generated
$A$-set with generators $x_1,\ldots,x_n$.  For every $x_i$, we have one of the
following possibilities for $Ax_i$:
\begin{compactenum}
\item it is free,
\item it is isomorphic to $A/At^k$ for some $n>0$,
\item there is a relation $ax_i=a'x_j$ for some $a,a'\in A$ and generator $x_j$
    of $X$.
\end{compactenum}
For (ii), the a.e.s. $0\ra A\xra{t^k}A\ra A/At^k\ra 0$ shows that $[A/At^k]=0$.
Now suppose we are in case (iii) and $i\neq j$.  Let $X_i\subseteq X$ be the
$A$-subset generated by $x_1,\ldots,x_{i-1},x_{i+1},\ldots,x_n$.  We have an
a.e.s. $0\ra X_i\hookrightarrow X\ra A/At^k\ra 0$ for some $k>1$ depending on
$a$; thus $[X]=[X_i]$.  Repeating this process if necessary, we can remove all
generators sharing a relation with second, distinct generator.

In case (iii) and when $i=j$, each generator $x_i$ satisfying a relation of the form $ax_i=a'x_i$ generates its own summand, so we may restrict our attention to $A$-sets $X$ generated by a single element of this form.  Here $X\cong A/(t^n=t^m)$ for some $n\geq m$.  When $m>0$, the a.e.s. 
\[0\ra A/(t^{n-m}=1)\xra{t^m}X\ra A/At^m\ra 0\] 
shows $[X]=[A/(t^{n-m}=1)]$. Then $[X]= [C_{n-m}]$ is equivalent to the pointed cyclic group of order $n-m$.

Thus, the only non-trivial generators of $G_0(A)$ are $[C_n]$, $n\geq 1$, and
$[A]$.  We now show there are no relations shared between these generators.

Recall that an $A$-set is torsion free when it has no zero-divisors.  Each $C_n$ is torsion free, hence admissibly projective by Proposition \ref{p:relprojtorsionfree}, and $A$ is projective.  Therefore, every a.e.s. $0\ra Y'\ra Y\ra Y''\ra 0$ with $Y''=C_n$ or $A$, splits and provides no nontrivial relations.  Moreover, there are \emph{no} morphisms $C_n\ra A$ and no
injective $A$-set morphisms of the form $A\ra C_n$ or $C_n\ra C_m$, for any $n,m>0$ and 
$n\neq m$, so that $A$ and the $C_n$ do not admit any admissible sequences in 
the $Y', Y$ positions. (To see there are no injective maps $f:C_n\ra C_m$ when 
$n\neq m$, the relation $t^j=f(1)=f(t^n)=t^nf(1)=t^{n+j}$ in $C_m$ implies 
$n=km$ for some $k\geq 1$.  Then $f$ can be injective only when $k=1$.)

Therefore the $[C_n]$, $n>0$, are non-trivial generators and
\[G_0(A)\cong \Z[A]\bigoplus\left( \bigoplus_{i=1}^{\infty}\Z[C_n]\right)\]
is an infinitely generated, free abelian group!  This varies greatly from the 
analogous result in commutative ring theory result that $G_0(k[t])\cong \Z$ 
where $k[t]$ is the polynomial ring in one variable with coefficients in a field 
$k$.  We can attribute this oddity of monoids once again to the lack of 
cancellation.
\end{example}

Let $A,B$ be noetherian monoids.  Recall that for a functor $F: A\set\ra B\set$
to be exact, it must preserve admissible morphisms.  Every exact functor $F$
induces a group homomorphism $G_0(A)\ra G_0(B)$.  If $A$ is a $B$-monoid via the monoid map $f:B\ra A$, it is immediate that $f^*$ is always exact.  

Before proceeding we recall the following definitions.  A \emph{filtration} for an $A$-set $X$ is a descending sequence of $A$-sets $0=X_n\subseteq\cdots \subseteq X_1\subseteq X_0=X$.  For convenience, we may write the filtration as $\{X_i\}_{i=1}^n$.  A \emph{refinement} of $\{X_i\}$ is a filtration $\{X'_j\}$ such that every $X_i$ occurs as one of the $X'_j$.  Note that Lemma \ref{l:devissage.refinement} (i) and (ii) are modifications of Zassenhaus' Lemma and the Schneider Refinement Theorem respectively, while (iii) and Theorem \ref{t:devissage} modify \cite[II.6.3]{WK}.  

\begin{lemma}\label{l:devissage.refinement}
\begin{compactenum}
\item If $X_2\subseteq X_1, X'_2\subseteq X'_1$ are $A$-subsets of an $A$-set $X$, then
\[\Big((X_1\cap X'_1)\cup X_2\Big)/\Big((X_1\cap X'_2)\cup X_2\Big)\cong \Big((X_1\cap X'_1)\cup X'_2\Big)/\Big((X_2\cap X'_1)\cup X'_2\Big).\]
\item Let $\{X_{i}\}_{i=1}^n$ and $\{X'_{j}\}_{j=1}^m$ be filtrations of an $A$-set $X$.  Then there are filtrations $\{Y_{i,j}\}$ and $\{Y'_{j,i}\}$ of $X$ such that the collections $\{Y_{i,j}/Y_{i,j+1}\}$ and $\{Y'_{j,i}/Y'_{j,i+1}\}$ are equivalent.
\item With notation as in (ii), we have $[X]=\sum_{i=1}^{n-1}[X_i/X_{i+1}]=\sum_{j=1}^{m-1}[X'_j/X'_{j+1}]$ in $G_0(A)$.
\end{compactenum}
\end{lemma}

\begin{myproof}
i)  Upon inspection, we see that both sides are equivalent to
\[(X_1\cap X'_1)/\Big((X_1\cap X'_2)\cup(X'_1\cap X_2)\Big).\]
ii)  Set $Y_{i,j}=(X_i\cap X'_j)\cup X_{i+1}$ and $Y'_{j,i}=(X_i\cap X'_j)\cup X'_{j+1}$ and notice that $\{Y_{i,j}\}$ and $\{Y'_{j,i}\}$ form filtrations when ordered lexicographically.  Moreover, $\{Y_{i,j}\}$ is a refinement of $\{X_{i}\}$, and $\{Y'_{j,i}\}$ is a refinement of $\{X'_{j}\}$, since $X_i=Y_{i,1}$ and $X'_j=Y'_{j,1}$.  Finally, (i) shows that $Y_{i,j}/Y_{i,j+1}\cong Y'_{j,i}/Y'_{j,i+1}$ so that the collections $\{Y_{i,j}/Y_{i,j+1}\}$ and $\{Y'_{j,i}/Y'_{j,i+1}\}$ are equivalent.\\
iii) We show $\sum_{i=1}^{n-1}[X_i/X_{i+1}]=\sum_{k=1}^{t-1}[Z_i/Z_{i+1}]$ where $\{Z_k\}_{k=1}^t$ is \emph{any} refinement of $\{X_i\}$. Then, applying (ii) gives
\[\sum [X_i/X_{i+1}]=\sum [Y_{i,j}/Y_{i,j+1}]=\sum [Y'_{j,i}/Y'_{j,i+1}]=\sum[X'_j/X'_{j+1}].\]
To show the invariance of the initial sum under refinement, we need only consider a refinement which adds a single term, namely $X_{i+1}\subseteq Z\subseteq X_i$ for some $1\leq i\leq n-1$.  This invariance follows immediately from the a.e.s. $0\ra Z/X_{i+1}\ra X_i/X_{i+1}\ra X_i/Z\ra 0$ which shows $[X_i/X_{i+1}]=[X_i/Z]+[Z/X_{i+1}]$.
\qed
\end{myproof}

\begin{theorem}\label{t:devissage}(Devissage)
Let $A$ be a monoid and $I\subseteq A$ a nilpotent ideal, i.e. $I^{n}=0$ for
some $n\geq1$.  Then the projection $\pi:A\ra A/I$ induces an isomorphism
$\pi^{*}:G_{0}(A/I)\ra G_{0}(A)$.
\end{theorem}

\begin{myproof}
Since $I$ is nilpotent, $\{I^iX\}_{i=0}^n$ is a filtration of $X$ which refines the trivial filtration $\{0,X\}$.  Then applying Lemma \ref{l:devissage.refinement}(iii) we obtain $[X]=\sum_{i=0}^{n-1}[I^iX/I^{i+1}X]$ which shows $\pi^*$ is surjective.  We now show that the map $\varphi:G_0(A)\ra G_0(A/I)$ given by $[X]\mapsto\sum_{i=0}^{n-1}[I^iX/I^{i+1}X]$ is a well defined inverse for $\pi^*$.

The fact that the compositions $\pi^*\varphi$ and $\varphi\pi^*$ are identity morphisms is clear from Lemma \ref{l:devissage.refinement}(iii).  We have that $\varphi$ defines a group homomorphism on the free abelian group generated by the isomorphism classes of finitely generated $A$-sets.
To show that $\varphi$ descends to a homomorphism on $G_0(A)$, it suffices to show that for every a.e.s. $0\ra X'\xra{f} X\xra{g} X''\ra 0$ of $A$-sets providing a relation $[X]=[X']+[X'']$ on $G_0(A)$, we must have $\varphi[X]=\varphi[X']+\varphi[X'']$ in $G_0(A/I)$.

Write $\{X_j\}_{j=0}^{2n+2}$ for the filtration
\[\{0=f(I^nX'),\ldots,f(IX'),f(X')=g^{-1}(I^nX''),\ldots,g^{-1}(IX''),g^{-1}(X'')=X\}\]
of $X$.  First notice that for $n+1\leq j\leq 2n+1$, $X_j/X_{j+1}\cong I^jX'/I^{j+1}X'$ is an $A/I$-set.  For the initial half of the filtration, it is clear that $X'\cup I^iX\subseteq g^{-1}(X'')$.  To see the opposite containment, let $ix''\neq 0$ be in $I^iX''$, $i\in I^i$, and $x=g^{-1}(ix'')$ be the unique element in the fiber over $ix''$ ($g$ is admissible).  Then $x''\neq 0$ and $g(ig^{-1}(x''))=ix''$, so it must be that $x=ig^{-1}(x'')\in I^iX$.  Hence, for $0\leq j\leq n-1$, each $g^{-1}(I^jX'')/g^{-1}(I^{j+1}X'')\cong I^jX''/I^{j+1}X''$ is an $A/I$-set so $\sum_{j=0}^{2n}[X_j/X_{j+1}]$ is an element in $G_0(A/I)$.  Of course, $X_n/X_{n+1}=0$ is an $A/I$-set.  By appealing to Lemma \ref{l:devissage.refinement}(iii) once more we obtain the first equality in our desired result.
\[\varphi[X]=\sum_{j=0}^{2n}[X_j/X_{j+1}]=\sum_{j=n+1}^{2n+1}[X_j/X_{j+1}]+0+\sum_{j=0}^{n-1}[X_j/X_{j+1}]=\varphi[X']+\varphi[X'']\]

\qed
\end{myproof}

\begin{example}
When $A=\F_1[x]/(x^n)$, the ideal generated by $x$ is nilpotent so the map
$\pi:A\ra A/(x)\cong \F_1$ induces an isomorphism $G_{0}(A)\cong G_0(\F_1)$ by Theorem \ref{t:devissage}.  By Example \ref{e:k0.f1}, we have $G_0(A)\cong \Z$.
\end{example}

For a ring $R$, the group $K_1(R)$ may be defined to be the abelianization of the \emph{infinite linear group} $GL(R)$.  The group $GL_n(R)$ consists of all linear automorphisms of $\oplus_{1}^{n}R$ and there is an inclusion $GL_n(R)\hra GL_{n+1}(R)$ defined by 
\[g\mapsto \left(\begin{array}{cc} g & 0\\ 0 & 1 \end{array}\right).\]
Then $GL(R)$ is union of the $GL_n(R)$. Let $A$ be a monoid and $X$ the free $A$-set on $n$ generators; then $GL_n(A)=\Aut(X)$ is easily computable.  An $A$-set morphism $A\ra A$ is completely determined by the image of $1$ and is invertible if and only if $1\mapsto u$ where $u\in A^{\times}$ is a unit.  From this it easily follows that $GL_n(A)\cong (\prod_1^n A^{\times})\times\Sigma_n$ where $\Sigma_n$ is the symmetric group of $n$ elements.  We then realize $GL_n(A)$ as the set of all $n\times n$ matrices having only one non-zero entry in each row and column with coefficients in $A^{\times}$.  The next result then follows from \cite[IV.1.27,4.10.1]{WK}.

\begin{proposition}
For any monoid $A$, we have $K_{1}(A)\cong A^{\times}\times \Z_{2}$.
\end{proposition}




\section{Model categories}\label{homotopytheory}

When $A=\F_1$, the category $A\set$ is equivalent to $\Sets_*$, the category of pointed sets.  Thus, when looking for a homological theory for $A\set$, investigating the well known homotopy theory associated to $\Sets_*$ is a good place to start.  In \cite{Quil} Quillen provided the theory of homotopical algebra which axiomatized the homotopy theory of topological spaces. Once a category meets the stated requirements, a homotopy theory is immediate, providing a theory of long exact sequences and derived functors.  It is possible for a category to have multiple homotopy theories associated to it each with its own corresponding \emph{model structure} (defined below).

It is not always clear what model structure will provide the ``best'' homotopy
theory, so we investigate the most natural model structure and the invariants,
i.e.  homotopy groups and derived functors, it provides.  Luckily, the
definition of the standard (Kan) model structure on $\Dop\Sets_*$ and the
projective model structure on $\Chg(R)$ can be made equally well for $A\set$.
This does not mean, a priori, the (projective) model structure will be the most
useful, but we will see that it provides invariants similar to what we expect.

This section provides the most basic definitions we require as a reference.  In
the next section we review the category of simplicial sets and the standard
model structure defined there.  Throughout, one may also use the categories
$\Ch(R)$ and $\mathbf{Top}$ to guide their intuition.  Interested readers are
directed to \cite{Quil} and \cite{Hov} for a more complete introduction to the
theory of model categories.

For a category $\cC$, let $\Map\cC$ denote the collection of all morphisms in
$\cC$ which we may consider a category by defining morphisms to be commutative
squares.  A \emph{functorial factorization} is an ordered pair of functors
$(\alpha,\beta): \Map\cC\ra\Map\cC$ such that $f=\beta(f)\circ\alpha(f)$ for all
$f$ in $\Map\cC$.

\begin{definition}
A \emph{model structure} on a category $\cC$ consists of three subcategories of
$\Map\cC$ called weak equivalences ($\weq$), cofibrations ($\cof$) and
fibrations ($\fib$), and two functorial factorizations $(\alpha,\beta)$ and
$(\gamma,\delta)$ satisfying the following properties:
\begin{compactenum}
\item If $gf$ is a composition of morphisms in $\cC$ and any two of $f$, $g$ or
    $gf$ are weak equivalences, then so is the third.
\item If $f$ and $g$ are morphisms such that $f$ is a retract of $g$ and $g$ is
    a weak equivalence (resp. cofibration, resp. fibration), then so is $f$.
\item Call a morphism that is both a cofibration (resp. fibration) and weak
    equivalence a \emph{trivial cofibration} (resp. \emph{trivial fibration}).
    In the diagram
\[ \SelectTips{eu}{12}\xymatrix{
    A \ar [r] \ar@{>->} [d]_f & B \ar@{->>} [d]^g\\
    C \ar [r] \ar@{-->} [ur]^h & D
} \]
the morphism $h$ exists when $f$ is a trivial cofibration and $g$ is a fibration or, $g$ is a trivial cofibration and $f$ is a fibration.  In the former case we say that trivial cofibrations have the \emph{left lifting property} (LLP) with respect to fibrations and, in the latter case, trivial fibrations have the \emph{right lifting property} (RLP) with respect to cofibrations.  The morphism $h$ is called a \emph{lift}.
\item For any morphism $f$, $\alpha(f)$ is a trivial cofibration, $\beta(f)$ is
    a fibration, $\gamma(f)$ is a cofibration and $\delta(f)$ is a trivial
    cofibration.
\end{compactenum}
A \emph{model category} is a category $\cC$ that is both complete and
cocomplete, together with a model structure.
\end{definition}

This is the definition for model structure (and category) given by Hovey.  As noted in \cite{Hov}, this definition differs from that of \cite{Quil} since Quillen merely required the existence of factorizations, not that they be functorial.  Also, Quillen does not require the category to be (co)complete, but that it need only contain all \emph{finite} limits and colimits.  It is shown in \cite{cls} that $A\set$ is both complete and cocomplete so this distinction is inconsequential.  The functorial factorizations are generally difficult to write down explicitly and to use for computations.  Since our primary goal is to find a \emph{computable} homological theory, this distinction is also unimportant to us.

Let $\cC$ be a model category with zero object $0$.  We say that an object $C$
of $\cC$ is \emph{fibrant} when $C\fib 0$ is a fibration and $C$ is
\emph{cofibrant} when $0\cof C$ is a cofibration.  If $C$ is \emph{any} object,
we can factor the (not necessarily cofibrant) morphism $0\ra C$ into $0\cof
X\tfib C$ where the object $X$ is cofibrant and weak equivalent to $C$.
Similarly, we may factor the (not necessarily fibrant) morphism $C\ra 0$ into
$C\tcof Y\fib 0$ so that $Y$ is fibrant and weak equivalent to $C$.  The object
$X$ (resp. $Y$) is called a \emph{cofibrant replacement} (resp.  \emph{fibrant
replacement}) for $C$.  Let $\cC_c$ and $\cC_f$ denote the full subcategories of $\cC$ whose objects are the cofibrant and fibrant objects in $\cC$ respectively. Given the functorial nature of the factorizations, there are functors $Q:\cC\ra\cC_c$ and $R:\cC\ra \cC_f$ such that $QC$, respectively $RC$, is cofibrant, respectively fibrant, replacement for $C$.

Let $\cC$ be any category and $S$ a collection of morphisms in $\cC$.  The
\emph{localization of} $\cC$ \emph{with respect to } $S$ is a category
$S^{-1}\cC$ together with a functor $q:C\ra S^{-1}\cC$ satisfying:
\begin{compactenum}
\item $q(s)$ is an isomorphism for every $s\in S$.
\item Any functor $F:\cC\ra\cD$ such that $F(s)$ is an isomorphism for every
    $s\in S$ factors through $q$ uniquely.
\end{compactenum}
This definition of localization is exactly analogous to that of 
rings and monoids (Section \ref{localization}) together with the universal
property.  Now let $\cC$ be a model category.  The localization of $\cC$ at the
collection of weak equivalences is the \emph{homotopy category} $\Ho\cC$ of
$\cC$.  Thus, $\Ho\cC$ is a category where every weak equivalence becomes an
isomorphism.  In particular every object $C$ is isomorphic to its functorial
cofibrant and fibrant replacements.  It is shown in \cite[1.2.10]{Hov}, and a
consequence of \cite[10.3.7]{WH}, that whenever $\cC$ is a model category,
$\Ho\cC$ exists and is itself a category.  The notation
$[A,B]=\Hom_{\Ho\cC}(A,B)$ is standard.

The following definitions are provided by \cite[I.4]{Quil}.  Let $\cC$ be a 
model category and $\gamma:\cC\ra\cC'$ and $F:\cC\ra \cD$ be functors.  Should 
it exist, the \emph{left derived functor} $L_{\gamma}F$ of $F$ with respect to 
$\gamma$ is the functor forming the diagram:
\[ \SelectTips{eu}{12}\xymatrix{
    \cC \ar [r]^{\gamma} \ar [dr]_{F} & \cC' \ar@{-->} [d]^{L_{\gamma}F}\\
     & \cD
} \]
together with a natural transformation $\varepsilon:L_{\gamma}F\circ\gamma\ra F$ which satisfies the following universal property.  Given any functor $G:\cC'\ra \cD$ and natural transformation $\eta:G\circ\gamma\ra F$, there is a unique natural transformation $\Theta:G\ra L_{\gamma}F$ making the following diagram commute.
\[ \SelectTips{eu}{12}\xymatrix{
    G\circ\gamma \ar [r]^{\eta} \ar@{-->} [d]_{\Theta} & F \\
    L_{\gamma} F\circ\gamma  \ar [ur]_{\varepsilon}
} \]
For the special case when $\gamma$ is the localization $\cC\ra\Ho\cC$, we write $LF$ for the left derived functor.  When $LF$ exists, $F$ maps weak equivalent objects of $\cC$ to isomorphic objects in $\cD$.  If $\cD$ is a model category, the \emph{total left derived functor}, also written $LF$, is the functor making the following diagram commute.
\[ \SelectTips{eu}{12}\xymatrix{
    \cC \ar [r]^{F} \ar [d]_{} & \cD \ar[d]^{} \\
    \Ho\cC \ar@{-->} [r]^{LF} & \Ho\cD
} \]
For the total left derived functor to exist, $F$ must map weak equivalences of $\cC$ to weak equivalences of $\cD$.  Hence, if $C$ is any object of $\cC$, we are free to compute $LF(C)$ using any object weak equivalent to $C$, e.g. a cofibrant replacement $QC$.  Notice that the use of the notation $LF$ should be clear from the context as the total left derived functor takes values in the homotopy category associated to a model category, rather than the (model) category itself.  It is shown in \cite[I.4.2,I.4.4]{Quil} that when $F$ carries weak equivalences in $\cC_c$ to:
\begin{compactenum}
    \item isomorphisms in $\cD$, the left derived functor $L_{\gamma}F$ exists.
    \item weak equivalences in $\cD$, the total left derived functor
        $LF:\Ho\cC\ra\Ho\cD$ exists.
\end{compactenum}
Similar definitions can be made for \emph{right}/\emph{total right} derived functors; see \cite[I.4]{Quil}.

Given two model categories $\cC$ and $\cD$, it is natural to ask when their homotopy categories are equivalent.  This motivates the following definitions. A functor $F:\cC\ra \cD$ is a \emph{left Quillen functor} when it is a left adjoint and preserves cofibrations and trivial cofibrations.  A functor $G:\cD\ra \cC$ is a \emph{right Quillen functor} when it is a right adjoint and preserves fibrations and trivial fibrations.  Now suppose $F,G$ form an adjunction with $F$ left adjoint to $G$.  It is a \emph{Quillen adjunction} when $F$ is a left Quillen functor.  It is shown in \cite[1.3.4]{Hov} that the adjunction $F,G$ is a Quillen adjunction if and only if $G$ is a right Quillen functor.

\subsection{Simplicial sets}\label{simplicialsets}

It would be helpful to find a homological theory for $A\set$ that does not require the use of simplicial $A$-sets.  Simplicial objects are generally cumbersome to work with and, in categories whose objects have less structure than groups, difficult to compute with.  The model structure on simplicial $A\set$ is a direct generalization of that on $\Dop\Sets$.  Before we introduce a new notion of ``chain complex'' we will present a homotopy theory for simplicial $A\set$ provided by Quillen's homotopical algebra (see \cite{Quil}). Here we recall the basic definitions and notation conventions before defining the standard model structure on $\Dop\Sets$.

Let $\Delta$ denote the category whose objects are the ordered sets
\[[n]=\{0<1<\cdots<n-1<n\}\]
and whose morphisms are order preserving set maps.  Of particular importance are
the face maps $\varepsilon_i:[n-1]\xra{}[n]$ and degeneracy maps
$\eta_i:[n+1]\xra{}[n]$, $0\leq i\leq n$, defined by
\begin{equation*}
\varepsilon_i(j)= \left\{
    \begin{array}{ll}
    j & \text{if} \ j<i,\\
    j+1 & \text{if} \ j\geq i
    \end{array} \right.
\end{equation*}

\begin{equation*}
\eta_i(j)= \left\{
    \begin{array}{ll}
    j & \text{if} \ j\leq i,\\
    j-1 & \text{if} \ j>i.
    \end{array} \right.
\end{equation*}
In short, the image of $\varepsilon_i$ misses $i$ and the image of $\eta_i$ doubles $i$.

The following is a standard result which can be found in nearly any text
containing simplicial theory (e.g. \cite{WH}).

\begin{proposition}\label{p:epimonic}
Every morphism $\alpha:[m]\ra[n]$ in $\Delta$ has a unique epi-monic
factorization $\alpha=\varepsilon\eta$ such that
\begin{equation*}
\begin{array}{ll}
\varepsilon=\varepsilon_{i_1}\cdots\varepsilon_{i_s} & \text{with} \ \, 0\leq
i_s\leq\cdots \leq i_1\leq m \vspace{2mm}\\
\eta=\eta_{j_{1}}\cdots\eta_{j_t} & \text{with a} \ \, 0\leq j_1<\cdots<j_t<n
\end{array}
\end{equation*}
\end{proposition}

Let $\cC$ be any category.  A simplicial object of $\cC$ is a contravariant
functor $X:\Dop\ra\cC$.  We use the following notation conventions:
\[X_n = X([n]), \qquad \partial_i = X(\varepsilon_i), \qquad \sigma_i =
X(\eta_i).\]
The simplicial face and degeneracy maps $\partial_i, \sigma_i$ satisfy the
following simplicial identities:
\[\partial_i\partial_j=\partial_{j-1}\partial_i \ \text{if} \ i<j\]
\[\sigma_i\sigma_j=\sigma_{j+1}\sigma_i \ \text{if} \ i\leq j\]
\begin{equation*}
\partial_i\sigma_j= \left\{
\begin{array}{ll}
\sigma_{j-1}\partial_i & \text{if} \ i<j\\
\text{id} & \text{if} \ i=j \ \text{or} \ i=j+1 \\
\sigma_j\partial_{i-1} & \text{if} \ i>j+1
\end{array} \right.
\end{equation*}

\begin{remark}
Suppose $\cC$ is a concrete category so that there is a faithful functor
$U:\cC\ra\Sets$.  In this setting we are able to talk about the elements an
object of $\cC$ contains.  An element $x\in X_n$ is called an \emph{$n$-cell},
\emph{$n$-simplex}, or simply a \emph{cell/simplex}.  Any cell $x$ in the image
of some $\sigma_i$ is \emph{degenerate} (as they do not contribute to homotopy
groups).  If $X_n$ contains a non-degenerate cell and for $i>n$ all cells in
$X_i$ are degenerate, we say $n$ is the \emph{dimension} of $X$ and call the
non-degenerate elements of $X_n$ \emph{top dimension} (or \emph{level})
\emph{cells} or simply \emph{top cells}.  If $x$ is a top cell in $X$, the
element $\partial_ix$ is its $i^{th}$ \emph{face}.  Also, the elements of $X_0$
are called \emph{vertices}.
\end{remark}

A morphism of simplicial objects $X, X'$ is a natural transformation
$f:X\Rightarrow X'$.  The category of simplicial objects of $\cC$ together with
these morphisms will be denoted $\Dop\cC$.  When $\cC$ has the necessary
products (coproducts), we may define the product (coproduct) of two simplicial
objects by
\[(A\times B)_n=A_n\times B_n \quad\text{and}\quad (A\amalg B)_n=A_n\amalg B_n\]
with the obvious induced face and degeneracy maps.  Note that to give a map
$f:\D^n\ra Y$ of simplicial sets, it is enough to define $f$ on the unique
non-degenerate top cell of $\D^n$ since $f$ commutes with both the face and
degeneracy maps.  For this reason, we may use the notation $\D^n\xra{x}X$ to
mean the map sending the unique non-degenerate $n$-cell in $\D^n_n$ to $x\in
X_n$.

Perhaps the most ubiquitous simplicial category is $\Dop\Sets$, whose objects
can be used to build simplicial objects in any category.  If $X$ is a simplicial
set and $C$ is a simplicial object of a category $\cC$ with all small
coproducts, construct the simplicial object $C\boxtimes X$ in $\Dop\cC$ by
\[(C\boxtimes X)_n=\coprod_{x\in X_n} C_n[x], \ \ \ C_n[x]:=C_n,\]
with face maps $\partial_i$ defined by $c[x]\mapsto \partial_ic[\partial_ix]$
and degeneracy maps defined analogously.  The notation $C_n[x]$ merely provides
the index $x$ in the notation for convenience.  That is, the face map sends the
element $c$ in the summand indexed by $x$ to the element $\partial_ic$ in the
summand indexed by $\partial_ix$.

There is an inclusion $\cC\hra\Dop\cC$ sending every object $C$ to the
simplicial object, also denoted $C$, with $C_n=C$ and all face/degeneracy maps
defined to be the identity map. Simplicial objects of this kind are called
\emph{constant}.  Thus, for any object $C$ and simplicial set $X$, we have the
simplicial object $C\boxtimes X$ of $\cC$.  Constant simplicial objects are a
special case of a broader class of nice simplicial objects.

\begin{definition}
A simplicial object $C$ is \emph{split} if there exist subobjects $N(C_m)$ of
$C_m$ such that the map
\[\coprod_{\eta:[n]\twoheadrightarrow [m]}N(C_m)\ra C_n,\]
where each morphism $\eta$ in $\Delta$ is surjective, is an isomorphism for all
$n\geq 0$.
\end{definition}

Split objects tend to separate the degenerate and non-degenerate cells as is the
case in $\Dop\Sets$ and $\Dop\cA$ where $\cA$ is an abelian category.  In fact,
the well known Dold-Kan theorem constructs a functor $K$ which assigns to every
chain complex $C$ a split simplicial object $KC$.  In this paper we attempt to
extend the Dold-Kan theorem to a nice class of non-abelian categories, like
$A\set$, taking advantage of $K$ along the way.

One of the most fundamental simplicial sets is the \emph{simplicial $n$-simplex}
$\D^n:\Dop\ra\Sets$ defined by $[k]\mapsto\Hom_{\D}([k],[n])$.
The \emph{boundary} $\snsphere \subseteq\D^n$ is the simplicial set consisting
of all non-identity injections $[k]\ra[n]$.  This effectively removes the
single, non-degenerate $n$-cell of $\D^n$ and may be thought of as a
\emph{simplicial $n-1$ sphere}. Finally, the \emph{$k$-horn}
$\nkhorn\subseteq\D^n$ is the simplicial set consisting of all non-identity,
injective morphisms $[k]\ra[n]$ \emph{except} the inclusion $[n-1]\ra[n]$ whose
image misses $k$.  This effectively removes the non-degenerate, $n$-cell and its
$k^{th}$ face.

Let $|\D^n|\subseteq \R^n$ denote the convex hull of the points $e_0,\ldots,e_n$
where $e_0=(0,\ldots,0)$ and $e_i$ has $i^{th}$ coordinate 1 and all others 0.
Then $|\D^n|$ consists of all points $(t_1,\ldots,t_n)\in\R^n$ such that
$t_i\geq 0$ for all $i$ and $\sum_i t_i\leq 1$. Thus $|\D^n|$ is a topological
space (with the subspace topology inherited from $\R^n$) called the
\emph{geometric realization} of $\D^n$.  A map $\alpha:[m]\ra[n]$ induces a map
$|\alpha|:|\D^m|\ra|\D^n|$ defined by $\alpha(e_i)=e_{\alpha(i)}$ and extending
linearly.

We may extend definition of geometric realization to any simplicial set $X$ as
follows.  For $n\geq 0$ define a topology on the product $X_n\times\D^n$ by
considering it as the disjoint sum $\amalg_{x\in X_n}\D^n$ and using the
topology of $|\D^n|$.  On the disjoint union
\[\coprod_{n\geq 0}X_n\times \D^n\]
define an equivalence relation $\sim$ by declaring $(x,t)\in X_m\times\D^m$ and
$(x',t')\in X_n\times\D^n$ equivalent when there is a morphism
$\alpha:[m]\ra[n]$ in $\D$ such that $X(\alpha)(x')=x$ and
$|\alpha(t')|=|\alpha(t)|$.  Note that $X(\alpha):X_n\ra X_m$ is the set map
obtained by applying $X:\Dop\ra\Sets$ itself to $\alpha$.  The quotient
topological space $\amalg(X_n\times\D^n)/\sim$ is the \emph{geometric
realization} of $X$.

A morphism $f:X\ra Y$ of simplicial sets induces a morphism $|f|:|X|\ra|Y|$
defined on cells by extending each map $x\mapsto f(x)$ linearly to
$|x|\ra|f(x)|$.  In this way the geometric realization defines a functor
$\Dop\Sets\ra\mathbf{Top}$ where $\mathbf{Top}$ is the category whose objects
are topological spaces and whose morphisms are continuous functions.

\begin{definition} \label{d:ssetmodelstructure}
Define a morphism $f:X\ra Y$ in $\Sets$ to be a weak equivalence when the
corresponding continuous map $|f|:|X|\ra|Y|$ of topological spaces induces
isomorphisms $\pi_n(|X|)\ra\pi_n(|Y|)$ for all $n\geq 0$.  A morphism $p$ is a
fibration when it has the RLP for all diagrams of the form
\[ \SelectTips{eu}{12}\xymatrix{
    \nkhorn \ar [r] \ar [d] & X \ar [d]^p\\
\D^n \ar [r] \ar@{-->} [ur]^{\exists} & Y
} \]
where the left vertical map is the inclusion.  A morphism is a cofibration when it has the LLP with respect to all trivial fibrations.  It is shown in \cite[3.6.5]{Hov} that these collections of morphisms define a model structure on $\SSets$.  By \cite[3.2.1]{Hov}, a map $p$ is a trivial fibration when it has the RLP for all diagrams of the form
\[ \SelectTips{eu}{12}\xymatrix{
\snsphere \ar [r] \ar [d] & X \ar [d]^p\\
\D^n \ar [r] \ar@{-->} [ur]^{\exists} & Y
} \]
where the left vertical map is the inclusion, and consequently, a morphism is a
cofibration if and only if it is an injection.  Trivial cofibrations are also
called \emph{anodyne extensions}.
\end{definition}

For an exhaustive study of $\SSets$ and this model structure, see chapter 3 of \cite{Hov}.  As pointed out in \cite[3.6.6]{Hov}, the definitions of this model structure are made equally well for pointed simplicial sets and do form a model structure on $\SSets_*$.  When discussing pointed simplicial sets, we use $*$ to denote the basepoint.

\begin{remark}\label{concretems}
This model structure on $\Dop\Sets$ can be used to define a model structure on the category of simplicial objects in many concrete categories.  Let $\cC$ be a concrete category with faithful (forgetful) functor $U:\cC\ra\Sets$.  Define a model structure on $\Dop\cC$ by declaring a morphism $f$ to be a weak equivalence (resp., cofibration, resp., fibration) whenever $U(f)$ is a weak equivalence (resp., cofibration, resp., fibration) in $\Dop\Sets$.  When $\cC=R\mod$, we obtain a model structure on $\Dop R\mod$ and the Dold-Kan correspondence (see \cite[8.4]{WH}) provides an equivalence of categories $\Chg(R)\cong\Dop R\mod$ which descends to $\bD_{\geq0}(R)\cong\Ho \Dop R\mod$. Since computing homology in $\Ch(R)$ is equivalent computing homotopy groups in $\Dop R\mod$, there is no need to work with simplicial $R$-modules.
\end{remark}


A fibration in $\SSets$ is traditionally called a \emph{Kan fibration} and may be defined element-wise in the following manner.  A morphism $p:X\ra Y$ is a Kan fibration when: for every $n\geq 0, y\in Y_{n+1}$ and $k\leq n+1$, if $x_0,\ldots,x_{k-1},x_{k+1},\ldots,x_n+1\in X_n$ are such that $\d_iy=p(x_i)$ and $\d_ix_j=\d_{j-1}x_i$ for all $i<j$ and $i,j\neq k$, then there exists $x\in X_{n+1}$ such that $p(x)=y$ and $\d_ix=x_i$ for all $i\neq k$.  Applying this definition to the basepoint map $X\ra *$ in $\SSets_*$ gives the lifting properties $X$ must satisfy to be fibrant.


Recall that for two simplicial sets $X$ and $Y$, the \emph{mapping} or \emph{function complex} $\Map(X,Y)$ is the simplicial sets with $\Map(X,Y)_n=\Hom_{\SSets}(X\times\D^n,Y)$ and boundary map $\d_kf$ the composition $X\times \D^{n-1}\xra{1\times i_k} X\times\D^n\xra{f} Y$ where $i_k$ is the inclusion of $\D^{n-1}$ into $\D^n$ as the $k^{th}$ face.

Let $X$ be a simplicial set and $x,y\in X_0$ be vertices.  Let $\pi_0(X)$ denote $X/\sim$ where $\sim$ is the equivalence relation on $X_0$ generated by the relations $x\sim y$ whenever there exists $z\in X_1$ with $\d_1z=x$ and $\d_0z=y$.  Choosing a distinguished vertex $*\in X_0$ makes $\pi_0(X)$ a pointed set with basepoint the homotopy (equivalence) class of $*$.  We use the notation $\pi_0(X,*)$ when we wish to make the basepoint explicit.  Note that the set of equivalence class of $\pi_0(X)$ are in one-to-one correspondence with the connected components of $X$.  In particular, $X$ is connected when $\pi_0(X)$ has a single element.

In general, when $X$ is fibrant, the fiber $F$ over $*$ in $\Map(\D^n,X)\ra\Map(\snsphere,X)$ is fibrant (see \cite[3.3.1]{Hov}).  In this case we define (see \cite[3.4.4]{Hov}) the $n^{th}$ \emph{homotopy group of} $X$ to be $\pi_n(X,*)=\pi_0(F,*)$ .  That is, an element of $\pi_n$ is a homotopy class of a map $\D^n\ra X$ such that the composition $\snsphere\ra\D^n\ra X$ is the constant map $*$.  As we typically identify the sphere $S^n$ with the cofiber of $\snsphere\hra\D^n$, this definition agrees with our intuition.

When $X$ is \emph{any} simplicial $A$-set, we define the homotopy groups
$\pi_n(X)$ to be $\pi_n(RX)$ where $RX$ is a fibrant replacement for $X$.

\begin{remark}\label{r:homotopycomp}
Let $X$ be a fibrant, pointed simplicial set.  To give a very convenient
element-wise description of the homotopy groups, let
\[C_n(X,*)=\bigcap_{i=0}^n\ker(\d_i:X_n\ra X_{n-1})=\{x\in X_n \ |\ \d_ix=*,
0\leq i\leq n\}\]
and
\[B_n(X,*)=\bigcap_{i=0}^{n-1}\ker(\d_i:X_{n+1}\ra X_n)=\{x'\in X_{n+1}\ |\
\d_ix=*, 0\leq i\leq n-1\}.\]
Define a relation on $C_n(X,*)$ by declaring $x\sim x'$ whenever there is a
$y\in B_n(X,*)$ such that $\d_{n+1}y=x$ and $\d_ny=x'$.  When $X$ is fibrant,
$\sim$ is an equivalence relation and $\pi_n(X,*)$ is the quotient
$C_n(X,*)/\sim$.  Equivalently, we may realize $\pi_n$ in the following way.  In
the diagram
\[B_{n+1}(X,*)\da{\d_{n+1}}{\d_n}B_n(X,*)\da{\d_n}{d_{n-1}} B_{n-1}(X,*)\]
both $\d_n,\d_{n-1}:B_n(X,*)\ra B_{n-1}(X,*)$ coequalize
$\d_{n+1},\d_n:B_{n+1}(X,*)\ra B_n(X,*)$, hence induce
$\bar{\d}_n,\bar{d}_{n-1}:\coeq(\d_{n+1},\d_n)\ra B_{n-1}(X,*)$ by the universal
property for coequalizers.  Then $\pi_n(X,*)=\ker(\bar{\d}_n)\cap\ker
(\bar{\d}_{n-1})$ since $\sim$ is an equivalence relation.  That is, given
$x_1,x_2$ in $C_n(X,*)$, we have $x_1\sim x_2$ if and only if $x_1$ and $x_2$
are identified in $\coeq(\d_{n+1},\d_n)$.
\end{remark}

We now work in $\SSets_*$ and use the notation $\pi_n(X)$ for $\pi_n(X,*)$.  To 
obtain long exact sequences we must consider fiber sequences $F\xra{i} E 
\overset{p}\fib B$ where $F$ is the fiber over $*$.  Note that the previous 
sequence need only be isomorphic in $\Ho\SSets$ to a sequence where each term is 
fibrant.  Thus we may replace the given sequence with the new sequence where 
each term is now fibrant (though this may change the homotopy groups of the 
fiber).  It is a standard result that when $B$ is fibrant and $p$ is a 
fibration, $E$ and $F$ are also fibrant.  In any case we obtain (see 
\cite[3.4.9]{Hov}) a long exact sequence
\[\cdots\ra\pi_{n+1}(B)\xra{\d}\pi_n(F)\xra{i}\pi_n(E)\xra{p}\pi_n(B)\xra{\d}\cdots\]
where the \emph{boundary map} $\d$ is defined as follows (see 
\cite[3.4.8]{Hov}).  Given $b\in C_n(B,*)$, the fibration condition provides 
$e\in E_n$ with $p(e)=b$ and $\d_ie=*$ for all $i<n$.  The equivalence class of 
$\d_ne$ in $\pi_{n-1}(F)$ is independent of the choice of $e$ and thus induces a 
map $\d:\pi_{n}(B)\ra\pi_{n-1}(F)$ for all $n\geq 1$.

\begin{remark}
Given a functor $F:\SSets_*\ra\SSets_*$, we are interested in computing, should
it exist, the total left derived functor
\[\Ho\SSets_*\xra{LF}\Ho\SSets_*\]
By analogy with $R\mod$, if we wish to compute the $n^{th}$ left derived functor
$L_nF$ of $F$, we may include $\Sets_*\hra\SSets_*$ by considering any pointed
set $X$ as a constant simplicial set.  Since all simplicial sets are cofibrant,
finding a cofibrant replacement $QX$ for $X$ is unnecessary so
$L_nF(X)=\pi_n(F(X))$. Note that finding a fibrant replacement of $F(X)$ is
required for computing $\pi_n$.

Notice how the situation here is opposite from that in $\Ch(R)$ with the
projective model structure.  Every chain complex is fibrant so that the
fibrant replacement functor $R$ is superfluous.  However, not every chain
complex is cofibrant so that we must compute a cofibrant replacement (e.g.
projective resolution) to compute the $n^{th}$ left derived functor of a right
Quillen functor.  Computing the projective resolution of a chain complex
concentrated in degree 0 is not terribly difficult which allows a computable
theory of derived functors for $R\mod$.  We will see in the next section that
$\SAset$ inherits the worst of both situations. \end{remark}

\section{Simplicial $A$-sets}\label{simplicialasets}

A monoid $A$ may be realized as a category $\cA$ having a single object $*$ with 
morphisms $\operatorname{End}(*)=A$.  Then a \emph{simplicial} $A$\emph{-set}  
$X$ is a functor $X:\cA\ra\SSets_*$ with $X(*)=X$ and $A$-action realized by the 
simplicial set map $X(a):X\ra X$, $a\in \operatorname{End}(*)$;  more 
explicitly, for $x\in X_n$, $ax=X(a)(x)$.  In this way $\SAset$ is identified 
with the functor category $\SSets^{\cA}$ and \cite[11.7.3]{Hir} provides a model 
structure (Definition \ref{d:sasetmodelstructure} below) on $\SAset$ defined 
just as the model structure on $\SSets$ provided in \ref{d:ssetmodelstructure}.  
We now describe the morphisms which make up this model structure.

\begin{definition}\label{d:sasetmodelstructure}
The \emph{projective model structure} on $\SAset$ is defined as follows.  A map 
$p:X\ra Y$ is a \emph{fibration} when the underlying set map $U(p)$ is a 
fibration in $\SSets_*$.  Equivalently, $p$ is a fibration when it has the RLP 
for all diagrams of the form
\[ \SelectTips{eu}{12}\xymatrix{
A\boxtimes \nkhorn \ar [r] \ar [d] & X \ar [d]^p\\
A\boxtimes \D^n \ar [r] \ar@{-->} [ur]^{\exists} & Y
} \]
where $A$ is the free, rank 1, constant simplicial $A$-set and $\boxtimes$ is 
defined as in Section \ref{simplicialsets}.  A map is a \emph{weak equivalence} 
when it induces isomorphisms an all homotopy groups and a \emph{trivial fibration} when it is a fibration \emph{and} weak equivalence.  Furthermore, a \emph{cofibration} has the LLP with respect to all trivial fibrations and a \emph{trivial cofibration} is a map which is both a cofibration and weak equivalence.  In particular, 
$A\boxtimes\d\D^n\hra A\boxtimes\D^n$ is a cofibration and 
$A\boxtimes\nkhorn\hra A\boxtimes\D^n$ is a trivial cofibration.  Note that 
$\pi_n(X)=\pi_n(UX)$ where $UX$ is the underlying simplicial set of $X$.  Hence, the homotopy groups of $X$ may be computed by analyzing geometric realization $|UX|$ (see Section \ref{simplicialsets}) or by finding a fibrant replacement $RX$ for $X$ and applying Remark \ref{r:homotopycomp}.  It is evident from Remark \ref{r:homotopycomp} that the homotopy groups of a simplicial $A$-set are themselves $A$-sets.
\end{definition}

\begin{remark}
Call a morphism $X\ra Y$ of simplicial $A$-sets \emph{free} when for each $n\geq
0$, the map $X_n\ra Y_n$ is an inclusion of the form $X_n\hra X_n\vee F_n$ where
each $F_n$ is a free $A$-set.  It is shown in \cite[11.5.36]{Hir} that free
morphisms, as well as their retracts, are cofibrations.  In particular, a
simplicial $A$-set $F$ with $F_n$ free for all $n$ is cofibrant.  

Recall that
$f:X\ra Y$ is a \emph{retract} of $g$ when there is a commutative diagram
\[ \SelectTips{eu}{12}\xymatrix{
    X \ar [r] \ar [d]_{f} & W \ar [d]^g \ar [r] & X \ar [d]^{f}\\
Y \ar [r] & Z \ar [r] & Y
} \]
where the horizontal compositions are the identity map.  When $F$ is a level-wise free simplicial $A$-set, $0\ra F$ is a cofibration and a retract diagram
\[ \SelectTips{eu}{12}\xymatrix{
    0 \ar [r] \ar [d] & 0 \ar [d] \ar [r] & 0 \ar [d]\\
    X \ar [r]^{} & F \ar [r]^{} & X
} \]
implies $X_n$ is a retract of $F_n$, hence projective $A$-set, for all $n$.  Level-wise projective, simplicial $A$-sets occuring in this way are cofibrant.  However, in contrast to $R\mod$, a general level-wise projective, simplicial $A$-set need not be cofibrant.
\end{remark}


Hence, when $X$ is an $A$-set considered as a constant simplicial $A$-set, to
find a cofibrant replacement $QX$ for $X$ it is enough for $QX$ to be level-
wise  free with $\pi_0(QX)\cong X$ and $\pi_i(X)=0$ for $i>0$.  This situation
is  similar to that in $\Ch(R)$ as we are finding a ``projective resolution''
for  $X$.  When $A=\F_1$,  this agrees with the standard model structure on
$\SSets_*$.  In this setting, every inclusion is a cofibration so every object
of $\Dop\F_1\set$ is cofibrant.  Of course, every $\F_1$-set is free so that
\emph{every} inclusion $X\hra Y$ is a free morphism.  We now say what little
we can about fibrations.


\begin{proposition}\label{p:fibration}
Let $X\xra{p} Y$ be a fibration in $\Dop A\set$.  If $p_*:\pi_0(X)\ra\pi_0(Y)$ is onto, then $p$ is onto.  We say that fibrations are \emph{surjective on connected components}.
\end{proposition}

\begin{myproof}
We may assume that $Y$ consists of a single connected component and show that
$p$ is surjective.  It is enough to show that $p_0:X_0\ra Y_0$ is surjective since we may then use the LLP of the trivial cofibration $A\boxtimes\D^0\xra{i} A\boxtimes \D^n$ in the diagram
\[ \SelectTips{eu}{12}\xymatrix{
    A\boxtimes \D^0 \ar [r]^{} \ar [d]_i & X \ar [d]^p\\
    A\boxtimes \D^n \ar [r]_{\ \ \ y} \ar@{-->} [ur]^{\exists x} & Y
} \]
to lift all $n$-cells with $n>0$.  That is, if any single vertex of an $n$-cell $y$ can be lifted, then $y$ itself can be lifted.

Let $y_0,y_1\in Y_0$ be vertices connected by $y\in Y_1$ with $\d_0y=y_0$ and $\d_1y=y_1$.  If $p(x_0)=y_0$ for $x_0\in X_0$, we may use the lifting
\[ \SelectTips{eu}{12}\xymatrix{
    A\boxtimes \horn{1}{1} \ar [r]^{x_0} \ar [d] & X \ar [d]^p\\
    A\boxtimes \D^1 \ar [r]_{\ \ \ y} \ar@{-->} [ur]^{\exists x} & Y
} \]
to obtain $x\in X_1$ with $p(x)=y$, hence $y_1=\d_1y=\d_1p(x)=p(\d_1x)$.  
Similarly, if $p(x_1)=y_1$, replacing the map $x_0$ with 
$x_1:A\boxtimes\horn{1}{0}\ra X$, we obtain $x'\in X_1$ with $p(x')=y$ and 
$y=p(\d_1x')$.  Thus $y_0$ is in $p(X_0)$ if and only if $y_1$ is in $p(X_0)$.  
Since $Y$ is connected, we may do this for all vertices, so that $p_0:X_0\ra 
Y_0$ is surjective.
A very similar argument shows $p_1:X_1\ra Y_1$ is 
surjective.  Namely, given $y\in Y_1$ and a lift $x_0$ of $\d_0y$, the diagram 
above provides $x\in X_1$ mapping onto $y$.
\qed
\end{myproof}

\begin{remark}\label{r:fibrationadjunction}
The free $A$-set functor $F$ is left adjoint to the forgetful functor $U$ and
provides a Quillen adjunction
\[\SSets\underset{U}{\overset{F}{\rightleftarrows}}\SAset.\]
This shows that a map $X\ra Y$ in $\SAset$ is a fibration if and only if the 
underlying simplicial set map is a fibration in $\SSets$.  Also, for every 
cofibration (i.e. injection) $K\ra L$ of simplicial sets, $A\boxtimes K\ra 
A\boxtimes L$ is a cofibration in $\SAset$.  This provides an alternate proof of 
Propositon \ref{p:fibration}.
\end{remark}

It immediately follows that trivial fibrations $p:X\tfib Y$ are surjections 
since $\pi_0(X)\ra\pi_0(Y)$ is surjective.  Unfortunately, not all surjections 
are fibrations.  Fibrations and fibrant objects do not have a nice, explicit 
description as is the case in $\Ch(R)$.

Outside of surjectivity and Remark \ref{r:fibrationadjunction} there is little
we can say explicitly about fibrations $X\fib Y$.  The best description seems to
be given by the lifting property itself.  Generally, one does not compute
explicitly with fibrations and fibrant objects, rather uses their existence and properties to prove more abstract theorems.

Note that the homotopy groups 
$\pi_n(X)=\pi_n(UX,0)$ of a simplicial $A$-set $X$ are computed using the 
underlying simplicial set $UX$.  However, every face and degeneracy morphism of 
$X$ is an $A$-set morphism so that the formulas in Remark \ref{r:homotopycomp} 
hold equally well for simplicial $A$-sets.  That is, the kernel of each morphism 
is an $A$-set and the equivalence relation on $C_n(X,0)$ induced by the 
coequalizer is, in fact, a congruence.  This implies that the homotopy groups are themselves $A$-sets.

\chapter{Complexes and computation}\label{c:dacomplexes}

Using the homological theory of $R$-modules as a guide defining the derived category of a monoid $A$, our first task is to define what is a ``chain complex'' of $A$-sets.  In this chapter we present the definition of the category of double-arrow complexes $\Da(A)$.  It is shown that all $A$-sets have a projective resolution which is a double-arrow complex and there is an adjunction between the category of reduced, double-arrow complexes and $\SAset$.  We then compute $\Tor_0^A$ and $\Tor_1^A$ in a special case.

To finish the chapter we investigate the theory of $A$-set extensions.  Though the author has not investigated a cohomology or cohomotopy theory for $A$-sets, we are able to describe extensions and prove a correspondence between extensions and an invariant obtained produced by a simplicial set traditionally used to compute Hochschild homology.

\section{Double-arrow complexes}\label{dacomplexes}

Working with simplicial objects and homotopy groups explicitly is very tedious 
and motivates alternative methods for computing invariants.  Abelian categories 
$\cA$ avoid this mess through the Dold-Kan theorem which states that $\Dop 
\cA\cong \Chg(\cA)$ and that this equivalence descends to their homotopy 
categories.  In light of the more general homotopical algebra of Quillen, we can 
think of the traditional homological algebra of abelian categories, i.e.,  the 
study of chain complexes, as a \emph{preferable}  alternative to a theory using 
the language of simplicial objects.  That is, all the ``homotopical 
information'' of a simplicial object in $\cA$ is carried equally well by a chain 
complex and its homology.

In the proof of the Dold-Kan theorem, the functor which provides a (Moore) chain 
complex for every simplicial object can be generalized slightly.  Combining this 
with the observations made in Remark \ref{r:homotopycomp}, it seems that the key 
to coming up with a complex suitable for carrying the homotopical information of 
a simplicial $A$-set lies in extending a coequalizer sequence $\cdots X\rra Y\ra 
Z$ to the left.

We say that a diagram
\[X\da{r}{s}Y\da{u}{v} Z\]
\emph{commutes} when both
\[X\da{r}{s}Y\xra{u}Z\text{ \ \ and \ \ }X\da{r}{s}Y\xra{v}Z\]
commute (see Section \ref{preliminaries}).  That is, $ur=us$ and $vr=vs$.

\begin{definition}
A \emph{double-arrow complex} $(X_\cdot,r_\cdot,s_\cdot)$, or simply $X$, with
objects in a category $\cC$ is a commutative diagram
\[\cdots\da{r_{n+2}}{s_{n+2}}
X_{n+1}\da{r_{n+1}}{s_{n+1}}X_n\da{r_n}{s_n}X_{n-1}\da{r_{n-1}}{s_{n-1}}\cdots\]
The subscripts on the morphisms will generally be dropped, as they should be
clear from the context, so the commutativity conditions may be written
\[rr=rs \text{\ \ and\ \ } sr=ss.\]
We may refer to $r$ and $s$ as the \emph{boundary maps} of $X$.  A double-arrow
complex is \emph{reduced} when $sr=ss=0$.
\label{d:dac}
\end{definition}

Note that the commutativity conditions imply that each $r_n, s_n$ induces morphisms
$\bar{r}_n,\bar{s}_n:\coeq(r_{n+1},s_{n+1})\ra C_{n-1}$ in any category $\cC$ where the coequalizers exist.  We say that a double
arrow complex $X$ is \emph{bounded below} (resp.  \emph{above}) when there is an
$N\in\Z$ with $X_n=0$ for all $n<N$ (resp.  $n>N$).  A complex is bounded when
it is both bounded above and below.

\begin{definition}\label{d:dac.morphism}
Let $X$ and $Y$ be double-arrow complexes.  A \emph{morphism of double-arrow
complexes}, or \emph{complex map}, $f:X\ra Y$ is a collection of morphisms
$f_n:X_n\ra Y_n$, $n\in\Z$ such that the following diagrams commute:
\begin{equation}\label{eq:dacmorph}
\SelectTips{eu}{12}\xymatrix{
X_n \ar [r]^{f_n} \ar [d]_{r_n} & Y_n \ar [d]^{r_n}\\
X_{n-1} \ar [r]_{f_{n-1}}  & Y_{n-1} }
\qquad
\SelectTips{eu}{12}\xymatrix{
X_n \ar [r]^{f_n} \ar [d]_{s_n} & Y_n \ar [d]^{s_n}\\
X_{n-1} \ar [r]_{f_{n-1}}  & Y_{n-1} }
\end{equation}
namely, $rf=fr$ and $sf=fs$.  In general, we condense the previous commutative
squares into a single diagram
\[ \SelectTips{eu}{12}\xymatrix{
    X_{n}\ar[r]^{f_{n+1}}\ar@<.5ex>[d]^{r} \ar@<-.5ex>[d]_{s} &
    Y_{n}\ar@<.5ex>[d]^{r} \ar@<-.5ex>[d]_{s} \\
    X_{n-1}\ar[r]^{f_{n}}  & Y_{n-1}
}\]
which we say \emph{commutes} when both diagrams in \ref{eq:dacmorph} commute.
\end{definition}

Again, the subscripts will generally be dropped when there is no chance of 
confusion.  Let $\Da(\cC)$ (resp., $\Dar(\cC)$, resp., $\Dag(\cC)$) denote the 
category of double-arrow complexes (resp., reduced complexes, resp.,  complexes 
bounded below by 0) with objects in a category $\cC$ together with morphisms 
defined in \ref{d:dac.morphism}.  For simplicity we write $\Da(A)$ for the 
category of double-arrow complexes of $A$-sets.

As with chain complexes in abelian categories, we may define the 
\emph{translation} $X[p]$ of a double-arrow complex $X$ to be the complex with 
$X[p]_n=X_{n+p}$ and the obvious boundary maps.  In particular, $X[1]_0=X_1$ so 
that $X[1]$ is obtained by shifting the terms of $X$ to the right.

A map $f:X\ra Y$ in $\Da(A)$ is \emph{admissible} when each $f_n:X_n\ra Y_n$ is
admissible.  A sequence  $X\ra Y\ra Z$ of double-arrow complexes of
$A$-sets is \emph{exact at} $Y$ when each sequence of $A$-sets $X_n\ra Y_n\ra Z_n$ is exact for every $n$.  Therefore, a sequence $0\ra X\ra Y\ra Z\ra 0$ in $\Da(A)$ is an \emph{admissible exact sequence} (a.e.s.) when each sequence $0\ra X_n\ra Y_n\ra Z_n\ra 0$ is an a.e.s. of $A$-sets.

\begin{example} (Moore complex)
Let $X$ be a simplicial $A$-set.  Define a reduced, double-arrow complex in the following way.  For $n\geq 2$, set
\[(NX)_n=\bigcap_{i=0}^{n-2}\ker(X_n\xra{\partial_i}X_{n-1})\]
and set $(NX)_1=X_1, (NX)_0=X_0$.  Then the sequence
\[\cdots\da{\partial_{n+2}}{\partial_{n+1}}(N
X)_{n+1}\da{\partial_{n+1}}{\partial_n}(N X)_n\da{\partial_n}{\partial_{n-1}}(N
X)_{n-1}\da{\partial_{n-1}}{\partial_{n-2}}\cdots\da{\partial_1}{\partial_0}(N
X)_0\rra0\]
is a \emph{reduced}, double-arrow complex by the simplicial identity
$\partial_i\partial_j=\partial_{j-1}\partial_i$, for  $i<j$, and since the face map $\d_{n-1}:(NX)_n\ra (NX)_{n-1}$ is trivial.  Recall that the boundary maps of any reduced complex satisfy $sr=ss=0$ (see Definition \ref{d:dac}).  In this way we
obtain a functor $N:\SAset\ra\Darg(A)$ with $X\mapsto NX$.  We write $N_nX$ for $(NX)_n$.
\label{e:moore}
\end{example}

\begin{definition}\label{d:dac.homology}
Let $\cC$ be a pointed, concrete category closed under kernels and coequalizers.  The $n^{th}$ \emph{homology} of a 
double-arrow complex $X$ is the object $H_n(X)=\ker(\bar{r}_n)\cap\ker(\bar{s}_n)$ in $\cC$, where 
$\bar{r}_n,\bar{s}_n$ are induced by the universal property.  The following 
diagram places the homology
\[ \SelectTips{eu}{12}\xymatrix{
    & & H_n(X) \ar[d] \\
    X_{n+1} \ar@<.5ex>[r]^{r_{n+1}} \ar@<-.5ex>[r]_{s_{n+1}} & X_n \ar [r]
    \ar@<.1ex>[dr]^{r_{n}} \ar@<-1ex>[dr]_{s_{n}} &
    \coeq(r_{n+1},s_{n+1})\ar@<.5ex>[d]^{\bar{r}_{n}}
    \ar@<-.5ex>[d]_{\bar{s}_{n}}\\
    & & X_{n-1}
}\]
We say that $X$ is \emph{exact at $X_n$} when $H_n(X)=0$ and \emph{exact} when
it is exact at $X_n$ for all $n$.

We wish to stress the fact that the $n^{th}$ homology $H_n(X)$ is an object of $\cC$.  Thus, if $X$ is in $\Da(A)$, the homology is an $A$-set.
\end{definition}

\begin{example}\label{e:fibrant.homhomequal}
If $X$ is a fibrant simplicial $A$-set, then $\pi_n(X)=H_n(NX)$ where $N$ is the Moore functor of Example \ref{e:moore}.  As per Remark \ref{r:homotopycomp}, $\pi_n(X)$ may be computed from the diagram
\[N_{n+1}X\da{\d_{n+1}}{\d_n}N_nX\da{\d_n}{\d_{n-1}}N_{n-1}X\]
as
\[\pi_n(X)=\bigg(\ker(N_nX\xra{\d_n}N_{n-1}X)\bigcap\ker(N_nX\xra{\d_{n-1}}N_{n-1}X)\bigg)/\sim\]
where $\d_{n+1}x\sim\d_n x$ for all $x\in N_{n+1}X$.  Alternatively, $H_n(NX)$ is computed from the diagram by first applying the congruence $\coeq(\d_{n+1},\d_n)$, and then finding $\ker(\bar{\d}_n)\cap\ker(\bar{\d}_{n-1})$.  These computation differ only in the order in which we compute kernels and quotient by a congruence.  Since $X$ is fibrant, so that $\sim$ \emph{is} a congruence, these operations commute.
\end{example}

\begin{remark}\label{r:bpfree}
The use of two boundary morphisms in a double-arrow complex allow for a different definition of homology.  We could also define the $n^{th}$ 
\emph{homology} of a double-arrow complex to be the equalizer 
$H_n(X)=\eq(\bar{r}_n,\bar{s}_n)$.

It is interesting that this alternative definition of homology may be applied in categories (such as $\Sets$) which do not contain a zero object.  In this setting there is no such thing as a bounded 
double-arrow complex.  Rather there are double-arrow complexes of the form
\[\cdots\da{1}{1}\eq(r_m,s_m)\rra X_{m}\da{r_m}{s_m}\cdots
\da{r_{n}}{s_{n}}X_{n-1}\rra\coeq(r_n,s_n)\da{1}{1}\cdots\]
for some $m\geq n$; that is, complexes which are eventually constant.

By Example \ref{e:fibrant.homhomequal}, our definition of homology (\ref{d:dac.homology}) is compatible with the definition of homotopy groups.  This is the reason that homology in pointed categories makes use 
of kernels rather than equalizers.  Moreover, 
Example \ref{e:dac.abelian} below shows that the homology of a double-arrow complex 
specializes to the homology of a chain complex when considering abelian 
categories.  Although the homotopy groups of a simplicial set are independent of 
the choice of basepoint, the choice of a basepoint is still required so that one 
must work in $\SSets_*$ rather than $\SSets$.  Thus our 
double-arrow complexes must be ``based'' as well if we wish to preserve this 
link with homotopy.  We do not investigate a ``basepoint free'' homology any 
further in this thesis.
\end{remark}

\begin{example}\label{e:dac.abelian}
When $\cA$ is an abelian category, any chain complex
\[\cdots\xra{\d} A_{n+1}\xra{\d}A_n\xra{\d}A_{n-1}\xra{\d}\cdots\]
defines a double-arrow complex $A$ of the form
\[\cdots\da{\d}{0}A_{n+1}\da{\d}{0}A_n\da{\d}{0}A_{n-1}\da{\d}{0}\cdots\]
since the boundary maps satisfy the relation $\d^2=0$.  Moreover, $H_n(A)$ is 
the kernel of $A_n/\im(\d_{n+1})\xra{\overline{\d}_n} A_{n-1}$ which is 
equivalent to the usual definition of homology, namely $\ker(\d_n)/\im(\d_{n+1})$.
\end{example}

\begin{lemma}\label{l:dac.inducemap}
A morphism $f:X\ra Y$ in $\Da(\cC)$ induces morphisms $(f_n)_*:H_n(X)\ra H_n(Y)$ in $\cC$ for all $n$.  Namely, when $\cC=A\set$, the $(f_n)_*$ are $A$-set morphisms.  We use the notation $f_*:H_*(X)\ra H_*(Y)$ for the induced map.
\end{lemma}

\begin{myproof}
Let $g:X_n\ra\overline{X}_n$ and $h:Y_n\ra\overline{Y}_n$ be the coequalizers of
$r_{n+1}, s_{n+1}$ in their respective complexes.  Using the universal property
of colimits for $g$, we obtain a commutative diagram
\[ \SelectTips{eu}{12}\xymatrix{
    X_n \ar[r]^{f_n}\ar[d]_g & Y_n\ar[d]^h\\
    \overline{X}_n\ar@{-->}[r]_{\exists\bar{f}_n} & \overline{Y}_n
}\]
Now, the relations $r=\bar{r}g$ and $r=\bar{r}h$ together with the commutative
square above imply $f\bar{r}g=\bar{r}\bar{f}g$. Since $g$ is a surjection, we
have a commutative diagram
\[ \SelectTips{eu}{12}\xymatrix{
    \overline{X}_n \ar[r]^{\bar{f}_n}\ar[d]_{\bar{r}} &
    \overline{Y}_n\ar[d]^{\bar{r}}\\
    X_{n-1}\ar[r]_{f_{n-1}} & Y_{n-1}
}\]
so that $\bar{f}_n$ maps $\ker(\bar{r})\subseteq\overline{X}_n$ to
$\ker(\bar{r})\subseteq \overline{Y}_n$.  Doing likewise for $\bar{s}$ proves
the claim.
\qed
\end{myproof}

Let $f:X\ra Y$ be a morphism of double-arrow complexes.  When the induced map 
$f_*:H_*(X)\ra H_*(Y)$ is an isomorphism, we say that $f$ is a 
\emph{quasi-isomorphism}.  We will show that the induced map $f_*$ is 
admissible when $f$ is both admissible and surjective.  To prove this, we first 
recall a fact about equivalence relations that are associated with a coequalizer.

\begin{lemma}
Given a coequalizer $X\da{f}{g}Y\xra{h}\overline{Y}$ of $A$-sets, let $R$ denote the
congruence associated to the surjection $h$ so that $Y/R\cong \overline{Y}$.  If
$y,y'\in Y$ and $h(y)=h(y')$, then there is a finite (possibly empty) sequence
of elements $x_1,\ldots,x_n\in X$ and ordered pairs of maps $(p_i,q_i)$, $1\leq
i\leq n$, satisfying:
\begin{compactenum}
    \item for each $1\leq i\leq n$, either $p_i=f$ and $q_i=g$, or else $p_i=g$ 
        and $q_i=f$,
    \item $q_i=p_{i+1}$ for $1\leq i \leq n-1$,
    \item $p_1(x_1)=y$ and $q_n(x_n)=y'$.
\end{compactenum}
Here we say that the $x_i$ and $(p_i,q_i)$ are a \emph{zig-zag} connecting $y$
and $y'$.
    \label{l:zigzag}
\end{lemma}

\begin{myproof}
Define a relation on the elements of $Y$ by $y\sim y'$ whenever there is a zig-zag connecting $y$ and $y'$.  We first show that $\sim$ is a congruence on $Y$.  Writing $\widetilde{R}$ for the $A$-subset of $Y\times Y$ associated to the congruence $\sim$, it is apparent that $\widetilde{R}$ contains the image of the map $X\xra{f\times g}Y\times Y$ (any element $(y,y')=(f\times g)(x)$ has $y\sim y'$ using zig-zag $x$ with $(p_1,q_1)=(f,g)$). Since $R$ is the smallest congruence containing $\im(f\times g)$, if $\sim$ is a congruence, it must be the case that $\widetilde{R}=R$.

The relation $\sim$ is trivially reflexive using the empty zig-zag.  Now,
suppose $y\sim y'$ so there is zig-zag with elements $x_i\in X$ and maps
$(p_i,q_i)$, $1\leq i\leq n$, satisfying the conditions above.  Then the
sequence $x_{n-i}$ with maps $(q_{n-i},p_{n-i})$, $0\leq i\leq n-1$, form a
zig-zag connecting $y'$ to $y$; so $\sim$ is symmetric.  Next, if $y\sim y'$
and $y'\sim y''$, let $x_i$ with $(p_i,q_i)$, $1\leq i\leq n$, and $x'_j$ with
$(p'_j,q'_j)$, $1\leq j\leq m$, be the zig-zags connecting $y$ to $y'$ and $y'$
to $y''$ respectively.  Then the concatenation $x_1,\ldots,x_n,x'_1,\ldots,x'_m$
with maps $(p_1,q_1),\ldots,(p_n,q_n),(p'_1,q'_1),\ldots,(p'_m,q'_m)$ form a
zig-zag connection $y$ to $y''$.

Thus $\sim$ is an \emph{equivalence relation}.  What remains to show is that 
$\widetilde{R}\subseteq Y\times Y$ is also an $A$-subset, hence a \emph{congruence}.  
This follows immediately from the fact that $f$ and $g$ are $A$-set morphisms.  
Namely, if $\{x_i,(p_i,q_i)\}_{i=1}^n$ is a zig-zag from $y$ to $y'$, then 
$\{ax_i,(p_i,q_i)\}_{i=1}^n$ is a zig-zag from $ay$ to $ay'$.
\qed
\end{myproof}


\begin{proposition}
Let $f:X\ra Y$ be an admissible, surjective morphism of double-arrow complexes.
Then $f_*:H_n(X)\ra H_n(Y)$ is admissible for all $n$.
\label{p:adm.induce.adm}
\end{proposition}

\begin{myproof}
We have a diagram
\[ \SelectTips{eu}{12}\xymatrix{
    X_{n+1}\ar[r]^{f_{n+1}}\ar@<.5ex>[d]^{r} \ar@<-.5ex>[d]_{s} &
    Y_{n+1}\ar@<.5ex>[d]^{r} \ar@<-.5ex>[d]_{s} \\
    X_{n}\ar[r]^{f_{n}} \ar[d]_g & Y_{n}\ar[d]^h \\
    \overline{X}_{n}\ar[r]^{\bar{f}_{n}} & \overline{Y}_{n}
}\]
where $g$ and $h$ are the coequalizer maps.  Using Lemma \ref{l:dac.inducemap} 
and noting that $H_n(X)$ is the $A$-subset 
$\ker(\overline{r})\cap\ker(\overline{s})\subseteq \overline{X}_n$, the induced 
map $f_*:H_n(X)\ra H_n(Y)$ is simply the restriction 
$H_n(X)\hra\overline{X}_n\xra{\overline{f}_n}\overline{Y}_n$.  Therefore we need 
only show that $\bar{f}_n$ is admissible.

Suppose $\bar{x},\bar{x}'\in \overline{X}_n$ are such that 
$\bar{f}_n(\bar{x})=\bar{f}_n(\bar{x}')\neq 0$ and choose $x,x'\in X_n$ so that 
$g(x)=\bar{x}$ and $g(x')=\bar{x}'$.  Commutativity of the bottom square implies 
$hf_n(x)=hf_n(x')$, so there is a zig-zag connecting $f_n(x)$ to $f_n(x')$ with nonzero elements $y_i\in Y_{n+1}$ and maps $(p_i,q_i)$, $1\leq i\leq n$.  Since $f_{n+1}$ is surjective and admissible, we can lift each $y_i$ to a unique $x_i\in X_{n+1}$.  Moreover, since $f$ is a morphism of double-arrow complexes and each $p_i,q_i$ is one of the boundary maps $r$ or $s$, we have $fp_i(x_i)=p_if(x_i)=p_i(y_i)$.  In particular, $fp_1(x_1)=f(x)$ and $fq_n(x_n)=f(x')$ and, using the admissibility of $f$ once more, $p_1(x_1)=x$ and $q_n(x_n)=x'$.  Therefore, the $x_i\in X_{n+1}$ and maps $(p_i,q_i)$, $1\leq i\leq n$, are a zig-zag connecting $x$ to $x'$; that is $\bar{x}=\bar{x}'$.
\qed
\end{myproof}

\section{Projective resolutions}

Aside from being simpler to work with, double-arrow complexes provide a method 
for constructing projective resolutions.  In an abelian category, a projective 
resolution of an object $C$ is an exact sequence of projective objects 
$\cdots\ra P_2\ra P_1\ra P_0$ such that $C=\coker(P_1\ra P_0)$.  However, we are 
unable to replicate this in $A\set$ since not every $A$-set has a projective 
resolution using only admissible morphisms.  For example, consider $A=\F_1[x]$ and 
$X=(A\vee A)/(x^n\vee 0\sim 0\vee x^n)$.  The surjection $A\vee A\ra X$ is not 
admissible and has trivial kernel.  Projective resolutions given by 
double-arrow complexes provide a little more information.

\begin{definition}
Let $X$ be an $A$-set.  A \emph{projective resolution} for $X$ is an exact 
double-arrow complex of projective $A$-sets
\[\cdots\rra P_2\rra P_1\da{r_1}{s_1}P_0\]
together with a map $P_0\xra{\varepsilon} X$ that is the coequalizer of $r_1, 
s_1$.  The morphism $\varepsilon$ is called the \emph{augmentation} map.  A 
projective resolution is \emph{reduced} if it is reduced as a double-arrow 
complex (see Definition \ref{d:dac}).  Note that the previous definition of projective resolution is equivalent to the \emph{augmented complex}
\[\cdots\rra P_2\rra P_1\da{r_1}{s_1}P_0\da{\varepsilon}{0}X\da{}{}0\]
being exact.
\end{definition}

The \emph{length} of a projective resolution $P$ is, the least integer $\ell(P)
\geq 0$ such that $P_n=0$ for all $n>\ell(P)$.  When no such integer exists set
$\ell(P)=+\infty$.  Define the \emph{projective dimension} of an $A$-set $X$ to
be $pd(X)=\min\{\ell(P)\ |\ P\text{ is a projective resolution of } X\}$.

We now show that every $A$-set has a projective resolution.  Note that the use of kernels in the proof of Proposition \ref{p:projectiveresolution} is specific to our purposes and can be modified to use equalizers in order to obtain the same result for unpointed categories (see Remark \ref{r:bpfree}).  Afterward we show that every $A$-set has a reduced projective resolution.  A category $\cC$ has \emph{enough projectives} if for every $C$ in $\cC$, there is a projective $P$ and an epimorphism $P\ra C$.  Certainly $A\set$ has enough projectives for if $X$ is any $A$-set, we have a surjection $A[X]\ra X$ defined by $[x]\mapsto x$.

\begin{proposition}\label{p:projectiveresolution}
Every $A$-set $X$ has a projective resolution.
\end{proposition}

\begin{myproof}
Throughout this proof, the standard projection maps provided by a pullback will
be written $p_1,p_2$.  Also recall that we write $X/Y$ for the quotient of an
$A$-set $X$ by a congruence $Y$.

Find a surjection $P_0\xra{\epsilon}X$ with $P_0$ projective and let
$R_0$ be the pullback of $\varepsilon$ with itself.  We have the
exact sequence
\[R_0\da{p_1}{p_2}P_0\da{\varepsilon}{0}X\rra0.\]
Now, find a surjection $P_1\xra{h_1}R_0$ with $P_1$ projective and
define $r_1=p_1h_1$ and $s_1=p_1h_1$.  We now have a sequence
\[P_1\da{r_1}{s_1}P_0\da{\varepsilon}{0}X\rra0.\]
Once more, let $R_1$ be the pullback of $h_1$ with itself and
find a surjection $P_2\xra{h_2}R_1$.  Define maps $r_2=p_1h_2$
and $s_2=p_2h_2$ to obtain the sequence
\[P_2\da{r_2}{s_2}P_1\da{r_1}{s_1}P_0\da{\varepsilon}{0}X\rra0.\]
Repeating this process gives a sequence of projectives
\[\cdots\da{r_3}{s_3}P_2\da{r_2}{s_2}P_1\da{r_1}{s_1}P_0\da{\varepsilon}{0}X\rra0\]
which is a double-arrow complex since the coequalizer
\[R_n\da{p_1}{p_2}P_n\xra{h_n}R_{n-1}\]
implies
\[r_nr_{n+1}=p_1(h_np_1)h_{n+1}=p_1(h_np_2)h_{n+1}=r_ns_{n+1}.\]
A similar computation shows $s_nr_{n+1}=s_ns_{n+1}$.  To check exactness, let
$P$ denote the double-arrow complex of projectives and use
$P\da{\varepsilon}{0}X\rra0$ to denote the entire complex.  Note that
$R_n=\{(x,x')\in P_n\times P_n\ |\ h_n(x)=h_n(x')\}$ is a congruence on $P_n$, and since $h_{n}$ is
surjective, we have
\[\coeq(r_{n+1},s_{n+1})\cong
P_{n}/R_n\cong R_{n-1}.\]
For $n>0$, $H_n(P)$ is $\ker(\bar{r}_n)\cap\ker(\bar{s}_n)$ in
$\coeq(r_{n+1},s_{n+1})\da{\bar{r}_n}{\bar{s}_n}P_{n-1}$ and by construction, the induced maps
$\bar{r}_n, \bar{s}_n$ are simply the projection maps $p_1, p_2$.  So the
previous diagram becomes
\[R_{n-1}\da{p_1}{p_2}P_{n-1}\]
and $(x,y)\in\ker(p_1)\cap\ker(p_2)$ means $p_1(x,y)=x=0$ and $p_2(x,y)=y=0$ so
$(x,y)=(0,0)$.  Hence $H_n(P)=0$ for $n>0$.  That $H_0(P)=\coeq(r_1,s_1)=X$ is
clear.  \qed
\end{myproof}

\begin{remark}\label{r:res.simplify}
The projective resolution constructed in the previous proposition is quite
``large.''  The coequalizer in the diagram $X\da{f}{g}Y\xra{h} Z$ is equivalent
to $Y\times_Z Y\da{p_1}{p_2}Y\xra{h} Z$ even though the difference between $X$
and $Y\times_Z Y$ is (likely to be) very large.  Namely, $Y\times_Z Y$ is a
\emph{congruence} (see Section \ref{asubsetsandquotients}) whereas $X$
together with $f$ and $g$ \emph{generate} a congruence.  In fact
$Y\times_Z Y$ is the congruence generated by $X$, i.e.  there is no smaller
congruence containing the image of the map $X\xra{f\times g}Y\times Y$. 

Let notation be as in Proposition \ref{p:projectiveresolution}.  When
constructing a projective resolution of $X$, we do not need to work with
surjections $P_{n+1}\xra{h_{n+1}}R_n$, rather we may use coequalizer 
diagrams $P_{n+1}\da{r_{n+1}}{s_{n+1}}P_n$ where $R_n$ is the
congruence generated by $P_{n+1}$ (together with $r_{n+1}$ and $s_{n+1}$).

For a simple example of the difference between a ``large'' projective resolution (as in Proposition \ref{p:projectiveresolution}) and one that is much ``smaller,'' let $A=\F_1[x]$ and $X=A/x^nA$.  The surjection $A\xra{1}
X$ provides the pullback $A\times_X A$ which is an infinitely generated, free
$A$-set with generators
\[S=\{(1,1),(x^n,0),(0,x^n),(x^n,x^{n+i}),(x^{n+j},x^n)\ |\ i,j\geq 0\}.\]
Thus, Proposition \ref{p:projectiveresolution} yields the projective resolution
\[0\rra\vee_{s\in S}A\da{p_1}{p_2}A\rra0\]
for $X$ even though we know a smaller, much more intuitive resolution is
\[0\rra A\da{x^n}{0}A\rra0.\]
Many of the complications we encounter with $A$-sets do not arise in
abelian categories.  For example, every congruence on an $R$-module $M$ may be
identified with an $R$-submodule $M'\subseteq M$.  We can then find the
generators $x_1,\ldots,x_m$ of $M'$ \emph{as an $R$-module} and a surjection $\oplus_1^mR\ra M'$ to
avoid explicitly working with congruences.  The $A$-set $X$ in the example above is of this
form since $X$ is obtained from $A$ as the quotient by the ideal, hence $A$-subset, $x^nA$.  Of course, it is not always possible to realize a general $A$-set as the quotient of a projective $P$ by an $A$-subset.  In such a situation we are not
able to identify elements in a congruence on $Y$ with elements \emph{in} $Y$.
\end{remark}

\begin{proposition}
    Every $A$-set $X$ has a reduced projective resolution of the form:
\[\cdots\da{r_{n+1}}{0}P_n\da{r_n}{0}\cdots\da{r_3}{0}P_2\da{r_2}{0}P_1\da{r_1}{s_1}P_0.\]
\label{p:reduced.res}
\end{proposition}

\begin{myproof}
Begin by constructing a sequence $P_1\da{r_1}{s_1}P_0\da{\varepsilon}{0}X\rra0$
as in the proof of Proposition \ref{p:projectiveresolution}.  That is, $P_0$ and
$P_1$ are projective $A$-sets and $\coeq(r_1,s_1)=X$.  For $n\geq 1$, find a
projective $P_{n+1}$ with a surjection
$r_{n+1}:P_{n+1}\ra K_n=\ker(r_n)\cap\ker(s_n)$ and define $s_{n+1}:P_{n+1}\ra P_n$
to be the trivial map.  Then
\[\cdots\da{r_{n+1}}{0}P_n\da{r_n}{0}\cdots\da{r_3}{0}P_2\da{r_2}{0}P_1\da{r_1}{s_1}P_0\]
is a reduced projective resolution for $X$ since $\coeq(r_{n+1},0)=P_n/K_n$ and $\ker(\bar{r}_n)=\ker(r_n)/K_n$.  Note that when $n\geq 2$,
$K_n=\ker(r_n)$.
\qed
\end{myproof}

The reader should notice that reduced projective resolutions in $A\set$ are
not too different from projective resolutions in abelian categories.
After computing the initial coequalizer $P_1\da{r_1}{s_1}P_0\xra{\varepsilon}X$
the remainder of the (reduced projective) resolution is computed in the usual
way.

\subsection{Dold-Kan: A first attempt}

Recall the following defintions from \cite[IV.3.2]{GJ}.  Let $X$ be a simplicial set and obtain the \emph{$n$-truncated} simplicial
set $\tau_nX$ by removing \emph{all} cells above dimension $n$.  If
$\SSets_n$ denotes the category of $n$-truncated simplicial sets,
$\tau_n:\SSets\ra\SSets_n$ defines a functor whose left adjoint
$sk_n:\SSets_n\ra\SSets$ is called the \emph{$n$-skeleton}.  The functor $sk_n$
simply fills out an $n$-truncated simplicial set with the degenerate cells
necessary to make it a simplicial set.  More precisely, if $X$ is an $n$-truncated simplicial set, then for $k>n$ define
\[(sk_nX)_k=\coprod_{\eta:[k]\ra[m]}X_m\]
where the coproduct is taken over all surjective maps $\eta:[k]\ra[m]$ in the category $\D$.  All such summands consist of degenerate $k$-cells in $sk_nX$.  Given a summand $X_m\subseteq (sk_nX)_k$ corresponding to morphism $\eta$, the morphism $\eta$ has epi-monic factorization (see \ref{p:epimonic}) $\eta=\eta_{j_1}\cdots\eta_{j_t}$ consiting entirely of degeneracy maps; then $X_m$ is image of the composition $\sigma_{j_t}\cdots\sigma_{j_1}:X_m\subseteq(sk_nX)_{k-t}\ra (sk_nX)_k$.  In particular, we have $\pi_k(sk_n\tau_n X)=\pi_k(X)$ when $k<n$.  When $X$ is a trunctated, simplicial $A$-set, the previous formulas show that $sk_nX$ naturally has the structure of a simplicial $A$-set where the coproduct is the wedge sum $\vee$.  In this way we may consider the $n$-truncation and $n$-skeleton as functors acting on (truncated) simplicial $A$-sets.

The Moore functor $N:\SAset\ra\Darg(A)$ of Example \ref{e:moore} comprises half of the correspondence between simplicial $A$-sets and reduced, positively graded, double-arrow complexes.  We now show how to construct a sequence of truncated simplicial $A$-sets from an object of $\Darg(A)$ which is then used to define a left adjoint of $N$.  

\begin{example}
Let $X$ be a reduced, positive, double-arrow complex of $A$-sets, that is $X\in\Darg(A)$. We inductively define a sequence $X^{(0)},X^{(1)},\ldots$ of truncated simplicial $A$-sets where each $X^{(n)}$ is in $\SAset_n$.  Set $X^{(0)}=X_0$ and for $n>0$, define the cells of $X^{(n)}$ as follows:
\[X^{(n)}_{i}= \left\{
\begin{array}{ll}
    X^{(n-1)}_i & 0\leq i\leq n-1\\
    X_n\vee (sk_{n-1}X^{(n-1)})_n & i=n.
\end{array} \right. \]
Since $X^{(n-1)}$ is already an $(n-1)$-truncated simplicial $A$-set, all face
and degeneracy maps are defined and satisfy the simplicial identities everywhere
except at $X_n\subset X^{(n)}_n$.  Moreover, the image of the degeneracy maps
$\sigma_j:X^{(n)}_{n-1}\ra X^{(n)}_n$ misses $X_n$ so that we need only verify
the simplicial identities $\d_i\d_j=\d_{j-1}\d_i$ where $i<j$.  Define the
restrictions $X_n\hra X^{(n)}_n\xra{\d_i}X^{(n)}_{n-1}$ as follows:
\[\d_{i}= \left\{
\begin{array}{ll}
    0 & 0\leq i\leq n-2\\
    s_n & i=n-1\\
    r_n & i=n.
\end{array} \right. \]
We now verify the simplicial identity $\d_i\d_j=\d_{j-1}\d_i$ for $i<j$.  When
$j\leq n-2$, both sides evaluate to the trivial map.\vspace{2mm}\\
\underline{$j=n-1$}: we have $\d_i\d_j=\d_is_n$ and $\d_{j-1}\d_i=s_{n-1}\d_i$.
When $i<n-2$, both are trivial and when $i=n-2$, we have the equality
$\d_is_n=s_{n-1}s_n=0=s_{n-1}\d_{n-2}$ since $X$ is a reduced complex.\\
\underline{$j=n$}: here $\d_i\d_j=\d_ir_n$ and $\d_{j-1}\d_i=r_{n-1}\d_i$.
Again, when $i<n-2$, both evaluate to the trivial map.  When $i=n-2$, the
equality $\d_ir_n=s_{n-1}r_n=0=r_{n-1}\d_{n-2}$ holds since $X$ is reduced.  For
$i=n-1$, we have the equality $\d_ir_n=r_{n-1}r_n=r_{n-1}s_n=r_{n-1}\d_{n-1}$
since $X$ is a double-arrow complex.

Therefore, $X^{(n)}$ is an $n$-truncated, simplicial $A$-set.  From the
construction it is also clear that for every $n\geq 0$, we have an inclusion
$sk_nX^{(n)}\hra sk_{n+1}X^{(n+1)}$ and when $X$ is bounded above by $k$, we have
$sk_nX^{(n)}\cong sk_{n+1}X^{(n+1)}$ for all $n>k$.
\end{example}

\begin{definition}(Inverse Dold-Kan)
Let $X$ be an object in $\Darg(A)$ and
\[sk_0X^{(0)}\hra sk_1X^{(1)}\hra sk_2X^{(2)}\hra\cdots\]
be the sequence of simplicial $A$-sets obtained from the truncated, simplicial $A$-sets constructed above; that is
$X^{(n)}$ is $n$-truncated.  Then $KX=\mathrm{colim}\, sk_nX^{(n)}$ is also a simplicial
$A$-set from which we obtain the functor $K:\Darg(A)\ra\SAset$, $X\mapsto KX$.
Again, we write $K_nX$ for $(KX)_n$.  Also note that $KX$ is a \emph{split} simplicial $A$-set whose non-degenerate $n$-cells are contained in the summand
$X_n\subseteq K_nX$.
\label{d:inversedk}
\end{definition}


\begin{theorem}
The functors $N$ and $K$ form an adjunction
\[\Hom_{\SAset}(KX,Y)\cong\Hom_{\Darg(A)}(X,NY).\]
\end{theorem}

\begin{myproof}
For $n\geq 1$, we have $N_nKX=X_n\vee\im(\sigma_{n-1})$ since $\d_i\sigma_i=\id$ and $\d_i\sigma_{n-1}=\sigma_{n-2}\d_i=0$ for $0\leq i\leq n-2$.  Since $r,s:N_nKX\ra N_{n-1}KX$ are $\d_n=r_n$ and $\d_{n-1}=s_n$ respectively, there is an inclusion $X\hra NKX$ in $\Dar(A)$ given by the
identity map.  Moreover, if $f:KX\ra Y$ is a simplicial $A$-set morphism, the
composition
\[\Phi(f):X\hra NKX\xra{Nf} NY\]
is a morphism in $\Hom_{\Darg(A)}(X,NY)$.

Now, to define a map $\varepsilon:KNY\ra Y$ between
simplicial objects, it is enough to define it on the non-degenerate cells $N_nY\subseteq K_nNY$ in a way that is compatible with the face maps, since it is easily defined to be compatible with the degeneracy maps, i.e. $\varepsilon_n\sigma_j(y):=\sigma_j(y)$.  However, since the boundary maps of $NY$ are the restrictions $N_nY\xra{\varepsilon_n} Y_n\da{\d_n}{\d_{n-1}}Y_{n-1}$, it is clear that the $\varepsilon_n$ are compatible with the face maps $\d_i$ of $Y$.  Thus,  the $\varepsilon_n$ comprise a simplicial $A$-set map $\varepsilon:KNY\ra Y$ and if $g:X\ra NY$ is any map of (reduced)
double-arrow complexes, the composition
\[\Psi(g):KX\xra{Kg}KNY\ra Y\]
is a morphism in $\Hom_{\SAset}(KX,Y)$.

We now show the operations $\Phi$ and $\Psi$ are inverses.  Let $f:KX\ra Y$ be as above so that $(\Psi\circ\Phi)(f)$ is the composition
\[KX\hra KNKX\xra{KNf}KNY\xra{\varepsilon} Y.\]
Since this composition is a morphism of simplicial $A$-sets, we need only show
$(\Psi\circ\Phi)(f)$ and $f$ are equivalent when restricted to the
non-degenerate cells of $KX$.  Of course, $X_n\subseteq K_nX$ is the $A$-set of non-degenerate $n$-cells and the restriction of $(\Psi\circ\Phi)(f)_n$ is the composition
\begin{eqnarray*}
    X_n\hra X_n\vee sk_{n-1}X^{(n-1)}_n\ra
    (X_n\vee\im(\sigma_{n-1}))\vee(sk_{n-1}(KX)^{(n-1)}_n)&& \\
    \xra{K_nNf} N_nY \vee sk_{n-1}(NY)^{(n-1)}_n &\ra& Y_n
\end{eqnarray*}
which is precisely $x\mapsto f_n(x)$.  That is, the first two maps compose to be
the identity and $\im(K_nNf|_{X_n})\subseteq N_nY$.  Finally, let $g:X\ra NY$ be as
above so that $(\Phi\circ\Psi)(g)$ is the composition
\[X\hra NKX\xra{NKg} NKNY\ra NY\]
and at dimension $n$ we have
\[X_n\hra X_n\vee \im(\sigma_{n-1})\xra{N_nKg} N_nY\vee\im(\sigma_{n-1}) \ra
N_nY.\]
The composition of the first two map is precisely $g_n$, and the composition
$N_nY\hra N_nY\vee\im(\sigma_{n-1})\ra N_nY$ is, again, the identity map.
\qed
\end{myproof}

\section{$\Tor_*$}

Here we look at an example of derived functor of the tensor product functor $-\otimes_A Y:A\set\ra A\set$ defined by $X\mapsto X\otimes A$ (see Section \ref{tensorproduct}).  Of course, this extends to a functor $-\otimes_A Y:\SAset\ra\SAset$ by $X\mapsto X\otimes_A Y$ which is the simplicial $A$-set with $(X\otimes_A Y)_n=X_n\otimes Y$.  As usual, we write $\Tor_*^A(-,Y)(X)$, or simply $\Tor_*(-,Y)(X)$ for the left derived functor $L(-\otimes Y)(X)$ and $\Tor^A_n(-,Y)(X)$, or $\Tor_n(-,Y)(X)$ for the $n^{th}$ left derived functor $L_n(-\ot_A Y)(X)=\pi_n\circ L(-\otimes Y)(X)$.  The following is immediate from the right exactness of $-\ot_AY$.

\begin{lemma}\label{l:tor0}
$\Tor_0^A(X,Y)\cong X\ot_AY$.
\end{lemma}

Computing $\Tor_*$ is a difficult problem.  If $X$ is an $A$-set considered as a
constant simplicial $A$-set, we need to compute the homotopy groups of a
cofibrant replacement $QX$ for $X$.  In $\SAset$ it is simple to construct a level-wise
projective simplicial $A$-set $P$ with $\pi_0(P)=X$.  It is much more difficult
to construct $P$ so that it is aspherical.  The easiest way to compute 
homotopy groups is to consider the geometric realization $|P|$ and subsequently
$|P\ot_A Y|$ to determine the $n^{th}$ left derived functors.  We present a simple example exhibiting this process.

\begin{example}\label{e:tor1}
Let $X$ be an $A$-set and $a$ an element of $A$ such that the map $A\xra{a}A$ is
an injection.  We compute $\Tor_1^A(-,X)(A/aA)$ as follows.  Let $P$ be the
reduced free resolution $0\rra A\da{a}{0}A$ of $A/aA$ and $KP$ the associated
simplicial $A$-set (see Definition \ref{d:inversedk}).  Every non-degenerate 1-cell of $KP$ is an element of $x\in P_1$ and has $\d_0x=0$ and $\d_1x=ax$.  Since $a:A\ra A$ is an injection, $\d_1x\neq \d_1y$ when $x\neq y$ so that $|KP|$ has no cycles.  Note also that $KP$ is cofibrant since it is level-wise free.

Now, the double-arrow complex
\[P\ot_A X:\; 0\rra X\da{a}{0}X\rra0\]
determines the non-degenerate cell structure of $KP\ot_A X$.  By Lemma \ref{l:tor0}, we have $\Tor_0^A(-,X)(A/aA)\cong X/aX$. To compute $\Tor_1$, define $\d_0,\d_1:X\ra X$  by $\d_0x=0$ and $\d_1x=ax$ as above.  Again, every 1-cell lies in the connected component of 0 since $\d_0x=0$.  When $x\neq y$ are 1-cells, we have $\d_1x=\d_1y$ whenever $ax=ay$. Thus, all cycles in $|KP\ot X|$ are of the form
\[ \SelectTips{eu}{12}\xymatrix{
    {\scriptstyle 0} \ar@(ur,dr)^x@{-}}
    \hspace{10mm}
    \xymatrix{
    {\scriptstyle 0} \ar@/^/@{-}[rr]^x \ar@/_/@{-}[rr]_y & & {\scriptstyle
    ax=ay}
} \]
the left instance occurring when $ax=0$ and $x\neq 0$.  We then see that $\Tor_1(-,X)(A/aA)$ is the free, non-abelian group of order $\sum_{x\in X}\Big(|\d_1^{-1}x|-1\Big)$. This
computation is similar to the result for abelian groups which shows
$\Tor_1^{\Z}(\Z/m\Z,N)$ consists of all $m$-torsion elements in the abelian
group $N$.
\end{example}

\begin{remark}
One way to simplify computation is to work with homology groups rather than
homotopy groups.  Categorically this means computing with chain complexes rather
than simplicial objects and in abelian categories these notions are equivalent
(via Dold-Kan).

Let $k$ be a commutative ring with identity element.  The $k$-realization
functor $k[-]:A\set\ra k[A]\mod$ extends to a functor $\SAset\ra\Dop k[A]\mod$,
which we also call $k[-]$.  There are two functors $N,C:\Dop
k[A]\mod\ra\Chg(k[A])$ defined in the following way.  Let $M$ be a simplicial
$k[A]$-module.  The \emph{Moore} (or \emph{normalized}) chain complex $NM$ (not to be confused with the Moore functor of Example \ref{e:moore}) has
$N_0M=M_0$ and
\[N_nM=\bigcap_{i=0}^{n-1}\ker(\d_i:M_n\ra M_{n-1})\]
for $n>0$ with differential $N_nM\xra{d}N_{n-1}M$ defined by $d=(-1)^{n}\d_n$.
The \emph{unnormalized} chain complex $CM$ has $C_nM=M_n$ and differential
$C_nM\xra{d}C_{n-1}M$ the alternating sum $d=\sum_{i=0}^n(-1)^i\d_i$.  It is
shown in \cite[8.3.8]{WH} that $\pi_*(M)=H_*(NM)\cong H_*(CM)$.

Now, let $X$ be a simplicial $A$-set and $k[X]$ the associated simplicial
$k[A]$-module.  The \emph{simplicial homology of $X$ with coefficients in $k$}
is $H_*(X;k)=H_*(Ck[X])$.  Then the
inclusion $h:X\hra k[X]$ defined by $x\mapsto x$ is a morphism of simplicial
sets called the \emph{Hurewicz homomorphism}.  The induced map
\[h_*:\pi_*(X)\ra \pi_*(k[X])\cong H_*(Ck[X])=H_*(X;k)\]
provides the simplicial analogue to the well known relationship between the
homotopy and homology groups of $|X|$.  In particular, if $X$ is connected, then $H_1(X;\Z)$ is the abelianization of $\pi_1(X)$.
\end{remark}

\begin{example}
Continuing with Example \ref{e:tor1}, we can compute the homology groups of the chain complex $Ck[KP\ot_A X]$ of $k[A]$-modules.  Using Definition \ref{d:inversedk} together with the simplicial identities listed in Section \ref{simplicialsets}, we determine $Ck[KP\ot_A X]$ has the structure
\[\cdots\ra k[X]\oplus k[X]\oplus k[X]\xra{d_2}k[X]\oplus k[X]\xra{d_1}k[X]\ra0\]
where $d_1=\d_0-\d_1$ is
\[ (1,0)\mapsto 0, \quad (0,1)\mapsto -a\]
and $d_2=\d_0-\d_1+\d_2$ is
\[(1,0,0)\mapsto (1,0), \quad (0,1,0)\mapsto (1,0), \quad (0,0,1)\mapsto (0,0).\]
A simple computation of homology gives $H_0(KP\ot_AX;k)\cong k[X]/ak[X]\cong k[X/aX]$ and $H_1(KP\ot_AX;k)\cong\{f\in k[X]\ |\ af=0\}\cong \Tor_1^{k[A]}(k[X],k[A/aA])$.  As stated prior to this example, $H_1$ is the abelianization of $\Tor_1^A(-,X)(A/aA)$ and all elements $f\in H_1$ satisfy the singular condition that $af=0$ in $k[X]$.  For example, if $ax=ay$ in $X$ corresponds to a generator of $\Tor_1^A(-,X)(A/aA)$, then $a(x-y)=0$ in $k[X]$ and  $x-y\in H_1(KP\ot_AX;k)$.  Furthermore, if $ax_1=\cdots=ax_n=0$ in $X$, then $x_1,\ldots,x_n$ generate a free $k[A]$-submodule of $\Tor_1$ of rank $n$.

It is straightforward to determine $C_nk[KP\ot_AX]$, along with the boundary maps $d_n$, for $n>2$ in order to see that $H_n(KP\ot_AX;k)=0$ for $n\geq 2$.
\end{example}

\section{Extensions}\label{extensions}

We end the chapter with a small discussion on cohomology.  As $\Tor$ is the
natural start when considering left derived functors, $\Ext$ is the natural
start for right derived functors.  It is classical that when $M,N$ are two
$R$-modules, $\Ext_R^1(M,N)$ classifies the short exact sequences of the form
$0\ra N\ra X\ra M\ra 0$.  Such a s.e.s. is called an \emph{extension} of $M$ by
$N$.  If we hope to define $\Ext^1$ for $A$-sets, looking as such sequences is
a good place to start.

Let $X$ and $Y$ be $A$-sets.  Define an \emph{extension $\xi$ of $X$ by $Y$} to
be an admissible exact sequence of the form $0\ra Y\ra E\ra X\ra 0$.  Two
extensions $\xi,\xi'$ are equivalent when there is a commutative diagram
\[\SelectTips{eu}{12}\xymatrix{
    0 \ar[r] & Y \ar@{=}[d] \ar[r] & E \ar[d]^{\cong}\ar[r] & X\ar@{=}[d]\ar[r]
    & 0\\
    0 \ar[r] & Y \ar[r] & E' \ar[r] & X\ar[r] & 0
}\]
Let $Ext(X,Y)$ denote the set of all extensions of $X$ by $Y$.

Before attempting to describe $Ext$, we provide an idea of what an extension
looks like.  Suppose $0\ra Y\xra{i} E\xra{p} X\ra 0$ is an extension.  Then $X$
is the quotient of $E$ by an $A$-subset $Y$ which simply identifies $Y$ with 0;
so $E\cong Y\vee X$ \emph{as pointed sets} and we need only define an $A$-action on
$E$. Let $\cdot$ denote the $A$-action on $E$.  Since $Y$ is an $A$-subset, the
action on the $Y$ (set-theoretic) summand is determines the action on $i(Y)$.
If $x\in X$ and $ax\neq 0$ for some $a\in A$, then $a\cdot p^{-1}(x)=p^{-1}(ax)$
since both map to the same element in $X$, namely $ax$, and $p$ is admissible.
Therefore, the only freedom we have in defining an $A$-action on $E$ comes from
pairs $e\in\im(p^{-1}(X))$ and $a\in A$ such what $ap(e)=0$.  We simply need to extend the determined action to the pairs $(a,e)$.

\begin{proposition}
Let $X,Y$ be $A$-sets and define $Z=\{a[x]\subseteq A[X] \ | \ ax=0\}$.  There
is a one-to-one correspondence between extensions $0\ra Y\ra E\ra X\ra 0$ and
$A$-set morphisms $\varphi: Z\ra Y$ such that $\varphi(ab[x])=\varphi(a[bx])$
when $bx\neq 0$. \label{p:extensions}
\end{proposition}

\begin{myproof}
First, let $\xi$ denote the extension $0\ra Y\xra{i} E\xra{p} X\ra 0$ and let
$\cdot$ denote the $A$-action of $E$.  If $a[x]\in Z$ is nonzero, then $ax=0$ in $X$ and
$p(a\cdot p^{-1}(x))=ax=0$ so that $a\cdot p^{-1}(x)\in i(Y)$ may be identified
with an element of $Y$.  Define $\varphi:Z\ra Y$ by $\varphi(a[x])=i^{-1}(a\cdot
p^{-1}(x))$ if $a[x]\neq 0$ and $\varphi(0)=0$.  Notice $b\cdot p^{-1}(x)=p^{-1}(bx)$ whenever $bx\neq 0$ so that
\[\varphi(ab[x])=i^{-1}((ab)\cdot p^{-1}(x))=i^{-1}(a\cdot
p^{-1}(bx))=\varphi(a[bx])\]
as desired.  Since $i^{-1}$ is an $A$-set map on $i(Y)$, when $ab\neq 0$ in $A$ we have
\[\varphi(ab[x])=i^{-1}((ab)\cdot p^{-1}(x))=bi^{-1}(a\cdot
p^{-1}(x))=b\varphi(a[x]).\]
Now, if $ab=0$ and say $a[x]\in Z$, then $\varphi(ab[x])=\varphi(0[x])=0$ and
\[b\varphi(a[x])=bi^{-1}(a\cdot p^{-1}(x))=i^{-1}((ab)\cdot p^{-1}(x))i^{-1}(0\cdot p^{-1}(x))=0\]
as well.  This shows $\varphi$ is an $A$-set morphism and defines a function $\Phi$ which assigns to each extension
$\xi$ a map $\varphi=\Phi(\xi)$.

Conversely, let $\varphi:Z\ra Y$ be an $A$-set morphism satisfying
$\varphi(ab[x])=\varphi(a[bx])$ when $bx\neq 0$.  From the pushout diagram
\[\SelectTips{eu}{12}\xymatrix{
    Z \ar[r] \ar[d]_{\varphi} & A[X]\ar[d]_{\lrcorner} \\
Y \ar[r] & Y\vee_Z A[X]
}\]
define $E=(Y\amalg_Z A[X])/\sim$ where $\sim$ is the equivalence relation
generated by the relations  $ab[x]\sim a[bx]$ whenever $bx\neq 0$.  Let $i$ be
the composition $Y\ra Y\vee_Z A[X]\xra{\pi} E$.  To show $i$ is monic we need
only show $\pi$ is monic since the first map of the composition is the pushout
of the inclusion $Z\ra A[X]$, hence monic.  Let $y_1,y_2\in Y$ be such that
$\pi(y_1)=\pi(y_2)$ in $E$.  This can only happen when $y_1=\varphi(a[x])=a[x]$,
$y_2=\varphi(b[x'])=b[x']$ where $a[x],b[x']\in Z$ and $a[x]\sim b[x']$.  Now
$a[x]\sim b[x']$ if there is a $c\in A$ such that $a=bc$ and $x'=cx$, so
\[a[x]=bc[x]\sim b[cx]=b[x']\]
or $b=ac$ and $x=cx'$, so
\[b[x']=ac[x']\sim a[cx']=a[x].\]
In the first case, $\varphi(bc[x])=\varphi(b[cx])$ by assumption on $\varphi$ so
$y_1=\varphi(a[x])=\varphi(b[x'])=y_2$.  The latter case is identical, hence $i$
is monic.  In fact, $i$ is just the inclusion $y\mapsto y$.

Let $f$ be the composition $A[X]\ra Y\amalg_Z A[X]\ra E$ and notice that
$E= i(Y)\cup\im(f)$ with $i$ defined as above.  Let $p:E\ra X$ be a map defined
via the restrictions $ i(Y)\xra{0} X$ and $\im(f)\ra X$ given by $a[x]\mapsto
ax$.  Since $ i(Y)\cap\im(f)\subseteq Z$, the map is well defined and is an
admissible $A$-set morphism with $\ker(p)= i(Y)\cong Y$.  Hence, $E/Y\cong X$
and fits into an a.e.s. $0\ra Y\xra{i}E\xra{p}X\ra0$.  Note this constructions
provides a function $\Xi$ which assigns to each map $\varphi:Z\ra Y$ an
extension $\xi=\Xi(\varphi)$.

We now proceed to show that there is a one-to-one correspondence between
isomorphism classes of extensions and maps $\varphi:Z\ra Y$.  To this end we
show $\Phi$ and $\Xi$ are inverses. We first show $\xi\cong
(\Xi\circ\Phi)(\xi)$.

Let $\xi$ be the extension $0\ra Y\xra{i} E\xra{p} X\ra 0$, $\varphi=\Phi(\xi)$
and $\xi'=\Xi(\varphi)$ the extension $0\ra Y\ra E'\ra X\ra 0$.  Again, let
$\cdot$ denote the $A$-action of $E$.  Define $\psi:E'\ra E$ by the restrictions
$\psi|_Y=i$ and $\psi([x])=p^{-1}(x)$ otherwise.  It is clear that $\psi_Y$ is
an $A$-set morphism and $\psi(a[x])=a\cdot \psi([x])$ whenever $ax\neq 0$ in
$X$.  Therefore, we simply need to show that $\psi(a[x])=a\cdot\psi([x])$
whenever $[x]\in E'\backslash Y$ and  $ax=0$.  In this case
$a[x]=\varphi(a[x])=i^{-1}(a\cdot p^{-1}(x))$ in $E'=\Xi(\varphi)$ and
\[\psi(a[x])=\psi(i^{-1}(a\cdot p^{-1}(x)))=a\cdot p^{-1}(x)=a\cdot \psi([x]).\]
That $\varphi$ is surjective is clear and since $i$, $p^{-1}$ are monic,
$\varphi$ is an isomorphism.

Conversely, given a map $\varphi:Z\ra Y$ with $\varphi(ab[x])=\varphi(a[bx])$
when $bx\neq 0$ in $X$, let $\xi=\Xi(\varphi)$ be the extension $0\ra Y\ra E\ra
X\ra 0$, and $\varphi'=\Phi(\xi)$.  Noting that $Y\ra E$ is just the inclusion
$y\mapsto y$ and $E\ra X$ is the projection $e\mapsto e$ when $e\in E\backslash
Y$ and $e\mapsto 0$ otherwise, we have $\varphi'(a[x])=a\cdot p^{-1}(x)=a\cdot
x=\varphi(a[x])$ when $ax\neq 0$ in $X$.
\qed
\end{myproof}

One reason simplicial $A$-sets are needed when computing left derived functors
is that not every $A$-set had a projective resolution by admissible morphisms.
However, every $A$-set \emph{does} have an injective resolution by admissible morphisms.  That is, if $X$ is an $A$-set, it is shown in \cite[Thm. 6]{Bert} that $E=\pHom(A,X)$ is an injective $A$-set with action $a\cdot f(a')=f(aa')$ for $a\in A$ and $f\in E$; then $X\hra E$ defined by $x\mapsto (1\mapsto x)$ is an inclusion.  Since every injection is admissible and cokernels exist, we can construct
an injective resolution
\[0\ra E_0\xra{f_0} E_1\xra{f_1} E_2\xra{f_2}\cdots\]
where $X=\ker(E_0\ra E_1)$ and $E_{n+1}$ is an injective containing the cokernel
of the inclusion $\coker(f_{n-1})\hra E_{n}$.  However, applying $\Hom_A(Y,-)$ and
computing $H_1$ in the usual way, i.e. kernel modulo image, does not provide the
classification for extensions found above.

Whenever $\cC$ is a model category, $\cC^{op}$ is as well by dualizing the
definitions for fibration, etc..  Following this reasoning the correct setting for computing
right derived functors may be found by dualizing the definitions for the model
structure on $\SAset$ to $(\SAset)^{op}=\D A\set$, the category of cosimplicial
$A$-sets.  However, we point out that such a model
structure provides homotopy \emph{groups}, but the description of extensions given in Proposition \ref{p:extensions}
shows that $Ext(Y,X)$ is an $A$-set and generally \emph{not} a group.  Therefore
we should not expect to obtain $Ext(Y,X)$ explicitly as a ($n^{th}$) derived
functor.  Instead, in light of the $\Tor_1$ computation, we should hope to find
a group for which the elements of $Ext(Y,X)$ are the generators.

There is a way to produce monoid extensions via a ``cocyle condition.''  The
following construction is an adaptation of the cosimplicial Hochschild object
constructed in \cite[9.1]{WH}.  Note that a \emph{cosimplicial}  $A$-set $X$ is
a functor $X:\D\ra A\set$ and a morphism of two cosimplicial $A$-sets is a
natural transformation.  The cosimplicial identities involving the coface
($\d^i$) and codegeneracy ($\sigma_j$) maps are equivalent to identities for the
face ($\varepsilon_i$) and degeneracy ($\eta_j$) maps of $\D$ listed in Section
\ref{simplicialsets}.

Let $A$ be a monoid and $X$ an $A$-set.  Define $\Hom(A^{\land*}, X)$ to be
the cosimplicial $\F_1$-set with $\Hom(A^{\land*}, X)^n=\pHom(A^{\land n},X)$
for $n\geq 1$ and $\Hom(A^{\land 0}, X)^0 = X$, together with coface and
codegeneracy maps given by
\begin{eqnarray*}
    (\d^if)(a_0a_1,\ldots a_n) &=& \left\{
\begin{array}{ll}
a_0f(a_1,\ldots,a_n) & \text{if } i = 0,\\
f(a_0,\ldots,a_{i-1}a_i,\ldots,a_n) & \text{if } 0<i<n,\\
a_nf(a_0,\ldots,a_{n-1}) & \text{if } i=n.
\end{array} \right.\vspace{2mm}\\
(\sigma^jf)(a_1,\ldots,a_{n-1})&=& f(a_1,\ldots,a_i,1,a_{i+1},\ldots,a_{n-1})
\end{eqnarray*}
In the case where $R$ is a $k$-algebra and $M$ is an $R$-module, the
cosimplicial $k$-module $\Hom(R^{\ot *},M)$ may be considered a cochain complex
with coboundary maps the alternating sums $d^n=\sum_i(-1)^i\d^i$, $n\geq 0$. For
an $A$-set $X$, the coboundary map $d^n$ is not defined for
$\Hom(A^{\land*},X)^n$, but we can make sense of the kernel on a
$\F_1$-subset.

\begin{definition}
Write  $C^n(A,X)$ for the collection of all $f$ in $\Hom(A^{\land*},X)^n$ which
satisfy the following.  For each $(a_0,a_1,\ldots,a_n)$ in $A^{\land n}$,
there is at most one odd and one even value of $i$ for which
$(\d^if)(a_0,a_1,\ldots,a_n)$ may be nonzero.  We call an element $f\in
C^n(A,X)$ a \emph{$n$-cocyle} if for every $(a_0,a_1,\ldots,a_n)\in A^{\land
n}$, we have
\[(\d^{2j}f)(a_0,a_1,\ldots,a_n)=(\d^{2k+1}f)(a_0,a_1,\ldots,a_n)\] where $0\leq
2j,2k+1\leq n $ are the values of $i$ for which $(\d^if)(a_0,a_1,\ldots,a_n)$
may be nonzero.
\end{definition}

\begin{lemma}
If $f\in \pHom(A^{\land*},X)^2$ satisfies $f(a,b)=0$ whenever $ab\neq
0$ in $A$, then $f\in C^2(A,X)$.
\label{l:cocycle}
\end{lemma}

\begin{myproof}
Let $f\pHom(A^{\land*},X)^2$ satisfy $f(a,b)=0$ when $ab\neq 0$.  For any $a,b,c\in A$, we need to show that $(\d^if)(a,b,c)=0$ for at least one
odd and one even value of $i$.  This simply requires us to check all possible
cases which are summarized in the table below.  Note also that $(a,b)=0$ in $A\land A$ when either $a=0$ or $b=0$ (or both), and $f(0)=0$ for any pointed set map.

First suppose that $ab=0$.  Then $(\d^1f)(a,b,c)=f(ab,c)=f(0,c)=0$.  When $bc=0$, we have $(\d^2f)(a,b,c)=f(a,bc)=f(a,0)=0$ and when $bc\neq0$, we have $(\d^0f)(a,b,c)=af(b,c)=0$ by assumption on $f$.

Now suppose $ab\neq 0$. It is immediate that $(\d^3f)(a,b,c)=cf(a,b)=0$, again by assumption on $f$.  Finally, as above, when $bc=0$, $\d^2f=0$ and when $bc\neq 0$, then $\d^0f=0$.  The following table provides a brief summary.
    \begin{center}
        \begin{tabular}{ccc|cccc}
            $ab$ & $bc$ & $ac$ & $\d^0f$ & $\d^1f$ & $\d^2f$ & $\d^3f$\\
            \hline
            0 & 0 & 0 & $*$ & 0 & 0 & $*$   \\
            0 & 0 & 1 & $*$ & 0 & 0 & $*$   \\
            0 & 1 & 0 & 0 & 0 & $*$   & $*$   \\
            0 & 1 & 1 & 0 & 0 & $*$   & $*$   \\
            1 & 0 & 0 & $*$ & $*$   & 0 & 0 \\
            1 & 0 & 1 & $*$ & $*$   & 0 & 0 \\
            1 & 1 & 0 & 0 & $*$   &  $*$  & 0 \\
            1 & 1 & 1 & 0 & 0 & 0 & 0 \\
        \end{tabular}
    \end{center}
An entry is 0 when the value is zero, 1 when the value is non-zero and $*$
otherwise.
\qed
\end{myproof}

Let $A$ be a monoid and $X$ an $A$-set.  A \emph{square zero extension} of $A$ by $X$ is a noncommutative monoid $E$ and an admissible, surjective monoid morphism $\varepsilon:E\ra A$ such that $\ker(\varepsilon)$ is a square zero ideal (i.e. $ab=0$ for all $a,b\in\ker(\varepsilon)$) together with an $A$-set isomorphism $X\cong\ker(\varepsilon)$.  A \emph{square zero extension} is \emph{commutative} when $E$ is a (commutative) monoid.

\begin{proposition}
Let $A$ be a monoid, $X$ an $A$-set and
\[K=\{f\in\Hom(A^{\land*},X)^2\ |\ f(a,b)=0\text{ whenever } ab\neq 0\text{
in } X\}\subseteq C^2(A,X).\]
Then the 2-cocyles in $K$ are in one-to-one correspondence with the square zero extensions of $A$ by $X$.  Furthermore, if $K_c\subseteq K$ consists of all elements $f$ of $K$ satisfying $f(a,b)=f(b,a)$, then the \emph{commutative} square zero extensions are in one-to-one correpondence with the elements of $K_c$.
\end{proposition}

\begin{myproof}
Let $A$ be a monoid, $X$ an $A$-set and $\xi: X\hra E\ra A$ an extension of $A$ by $X$.  Write $\cdot$ for the multiplication in $E$.  Similar to Proposition \ref{p:extensions}, we have $E=A\vee X$ \emph{as pointed sets} and the multiplicative structure of $E$ is determined everywhere except for the set $Z=\{(a,b)\in A\land A\ |\ ab=0\}$.  For fixed $f$ in $K$, we show that the multiplication
\begin{eqnarray*}
    a\cdot b &:=& \left\{
\begin{array}{ll}
ab & \text{if } ab\neq 0,\\
f(a,b) & \text{if } ab=0.\\
\end{array} \right.
\end{eqnarray*}
defines an extension $\xi$ by defining the products $a\cdot b$ in $E$ for every pair $(a,b)$ in $Z$.  That is, once $f$ is verified to complete a multiplicative structure on the set $A\vee X$, the surjective monoid map $E\ra A$ sending every element of $X\subseteq E$ to 0 defines an extension.

It suffices to show that the multiplication given by $f$ is associative since 1 is not a zero divisor of $A$ and commutativity is given by assumption.  We must consider cases each of which follows from a row of the table given in the proof Lemma \ref{l:cocycle}.  The properties $f$ satisfies are:
\begin{compactenum}
\item $af(b,c)=f(ab,c)$ when $bc=0$, $ab\neq 0$,
\item $f(a,bc)=cf(a,b)$ when $ab=0$, $bc\neq 0$,
\item $f(ab,c)=f(a,bc)$ when $ac=0$, $ab, bc\neq 0$, and
\item $af(b,c)=cf(a,b)$ when $ab,ac,bc=0$.
\end{compactenum}
Now, let $a,b,c$ be elements of $A$ and consider the product $a\cdot(b\cdot c)$ in $E$.  If $bc=0$ and $ab\neq 0$ in $A$, then 
\[a\cdot(b\cdot c)=af(b,c)=f(ab,c)=(ab)\cdot c=(a\cdot b)\cdot c\]
by (i).  If $ab=0$ as well, then
\[a\cdot(b\cdot c)=af(b,c)=f(a,b)c=(a\cdot b)\cdot c\]
by (iv) and the commutativity of $X$ as an $A$-set.  On the other hand, when $bc\neq 0$ in $A$ and $ab=0$, we have
\[a\cdot(b\cdot c)=a\cdot (bc)=f(a,bc)=f(a,b)c=(a\cdot b)\cdot c\]
by (ii) and the commtutativity of $X$ as an $A$-set.  Finally, when $ab\neq 0$, we have
\[a\cdot(b\cdot c)=a\cdot(bc)=f(a,bc)=f(ab,c)=(ab)\cdot c=(a\cdot b)\cdot c\]
by (iii).  Hence, $f$ defines an extension $\xi=\Phi(f)$.  Conversely, if we are given an extension $\xi: X\ra E\ra A$, we obtain a 2-cocycle $f=\Psi(\xi):A\land A\ra X$ of $K$ by
\begin{eqnarray*}
    f(a,b) &:=& \left\{
\begin{array}{ll}
0 & \text{if } ab\neq 0,\\
a\cdot b & \text{if } ab=0.\\
\end{array} \right.
\end{eqnarray*}
Of course, when $ab=0$ in $A$, it must be that $a\cdot b\in\ker(E\ra A)\cong X$.  Hence, it is clear that this $f$ is a 2-cocycle and an element of $K$.

It also straightforward to show that the operations $\Phi$ and $\Psi$ are inverses.  If $f$ is in $K$, $\Phi(f):X\hra E\ra A$ is an extension where $a\cdot b:=f(a,b)$ in $E$ whenever $ab=0$ in $A$, and
\begin{eqnarray*}
    (\Psi\circ\Phi)(f)(a,b) &:=& \left\{
\begin{array}{ll}
0 & \text{if } ab\neq 0,\\
a\cdot b=f(a,b) & \text{if } ab=0\\
\end{array} \right.
\end{eqnarray*}
is precisely $f$.  Conversely, given an extension $\xi:X\hra E\ra A$, $\Psi(\xi):A\land A\ra X$ is the 2-cocycle of $K$ defined by
\begin{eqnarray*}
    \Psi(\xi) &:=& \left\{
\begin{array}{ll}
0 & \text{if } ab\neq 0,\\
a\cdot b & \text{if } ab=0.\\
\end{array} \right.
\end{eqnarray*}
Write $X\hra E'\ra A$ for the extension $(\Phi\circ\Psi)(\xi)$ and $\cdot'$ for the multiplication of $E'$.  Then $a\cdot' b=\Psi(\xi)(a,b)=a\cdot b$ is precisely the multiplicative structure of $E$.
\qed
\end{myproof}

\chapter{Geometry of monoids}\label{geometry.monoids}

The topics discussed in the chapter all require the language of monoid schemes.
Perhaps the reader familiar with commutative ring theory noticed the absence of
a discussion on normalization in Chapter 2.  Of course the reason is that the
normalization of a monoid need not be a monoid, but rather a monoid scheme.  The
chapter as a whole presents some basic results on the Picard group of a monoid
scheme and the Weil divisor class group of a normal monoid scheme.  These
results are taken from \cite{FW}.  For more information on monoid schemes see
\cite{cls}, \cite{chww}, \cite{chww-p} and \cite{Vezz}.

Recall that the set of prime ideals of a monoid $A$ is written $\MSpec(A)$ and
forms a topological space using the Zariski topology.  On $\MSpec(A)$ the
structure sheaf $\cA$ has stalk $A_{\frakp}$ at the point $\frakp$.  An
\emph{affine monoid scheme} is the topological space $\MSpec(A)$ together with
its structure sheaf $\cA$.  A \emph{monoid scheme} is a topological space $X$
together with a sheaf of monoids $\cF$ that is locally affine.

\section{Normalization}\label{sec:normalize}

If $A$ is a cancellative monoid, its normalization is the integral closure of
$A$ in its group completion $A_0$ and is universal for maps $A\ra B$ with $B$
normal (see Lemma \ref{pcuniversalnormal}). In contrast, consider the problem of
defining the normalization of a non-cancellative monoid $A$, which should be
something which has a kind of universal property for morphisms $A\to B$ with $B$
normal.

We will restrict ourselves to the case when the monoid $A$ is {\it partially
cancellative}  (or {\it pc}), i.e., a quotient $A=C/I$ of a cancellative monoid
$C$ (\cite[1.3, 1.20]{chww-p}). One advantage is that $A/\frakp$ is cancellative
for every prime ideal $\frakp$ of a pc monoid, and the normalization
$(A/\frakp)_\nor$ of $A/\frakp$ exists.

\begin{lemma}\label{pcuniversalnormal}
If $A$ is a pc monoid and $f:A\to B$ is a morphism with $B$ normal, then $f$
factors through the normalization of $A/\frakp$, where $\frakp=\ker(f)$.
\end{lemma}

\begin{myproof}
The morphism $A/\frakp\to B$ of cancellative monoids induces a homomorphism
$f_0: (A/\frakp)_0\to B_0$ of their group completions. If $a\in(A/\frakp)_0$
belongs to $(A/\frakp)_\nor$ then there is an $n$ so that $a^n\in A/\frakp$.
Then $b=f_0(a)\in B_0$ satisfies $b^n\in B$, so $b\in B$. Thus $f_0$ restricts
to a map $(A/\frakp)_\nor\to B$.
\end{myproof}

\begin{remark}\label{xz=yz}
The non-pc monoid $A=\langle x,y,z | xz=yz\rangle$ is non-cancellative, reduced
(\ref{r:nil(A)}) and even seminormal (see after Example \ref{norm.axes}), yet
has no one obvious notion of normalization in either the above sense or in the
sense of Definition \ref{normaldefn} below, since $0$ is a prime ideal. We have
restricted to pc monoids in order to avoid these issues.
%
\end{remark}

Thus the collection of maps $A\to(A/\frakp)_\nor$ has a versal property: every
morphism $A\ra B$ with $B$ normal factors through one of these maps. However, a
strict universal property is not possible within the category of monoids because
monoids are local. This is illustrated by the monoid $A=\langle x_1,x_2 |
x_1x_2=0\rangle$; see Example \ref{norm.axes} below. Following the example of
algebraic geometry, we will pass to the category of (pc) monoid schemes, where
the normalization exists.

\begin{definition}\label{normaldefn}
Let $A$ be a pc monoid. The {\it normalization} $X_\nor$ of $X=\MSpec(A)$ is the
disjoint union of the monoid schemes $\MSpec((A/\frakp)_\nor)$ as $\frakp$ runs
over the minimal primes of $A$. By abuse of notation, we will refer to $X_\nor$
as the normalization of $A$.

This notion is stable under localization: the normalization of
$U=\MSpec(A[1/s])$ is an open subscheme of the normalization of $\MSpec(A)$; by
Lemma \ref{l:local.normal}, its components are $\MSpec$ of the normalizations of
the $(A/\frakp)[1/s]$ for those minimal primes $\frakp$ of A not containing $s$.

If $X$ is a pc monoid scheme, covered by affine opens $U_i$, one can glue the
normalizations $\widetilde{U}_i$ to obtain a normal monoid scheme $X_\nor$,
called the {\it normalization} of $X$.
\end{definition}

\begin{remark}
The normalization $X_\nor$ is a normal monoid scheme: the stalks of $\cA_\nor$
are normal monoids. It has the universal property that for every connected
normal monoid scheme $Z$, every $Z\to X$ dominant on a component factors
uniquely through $X_\nor\to X$. As this is exactly like
\cite[Ex.\,II.3.8]{Hart}, we omit the details.
\end{remark}

Recall that the (categorical) product $A\times B$ of two pointed monoids is the
set-theoretic product with slotwise product and basepoint $(0,0)$.

\begin{lemma}
Let $A$ be a pc monoid. The monoid of global sections $H^0(X_\nor,\cA_\nor)$ of
the normalization of $A$ is the product of the pointed monoids $(A/\frakp)_\nor$
as $\frakp$ runs over the minimal primes of $A$.
\end{lemma}

\begin{myproof}
For any sheaf $\cF$ on a disjoint union $X=\coprod X_i$, $H^0(X,\cF)=\prod
H^0(X_i,\cF)$ by the sheaf axiom.
\end{myproof}

\begin{example}\label{norm.axes}
The normalization of $A=\langle x_1,x_2 | x_1x_2=0\rangle$ is the disjoint union
of the affine lines $\langle x_i\rangle$. The monoid of its global sections is
$\langle x_1\rangle\times\langle x_2\rangle$, and is generated by
$(1,0),(0,1),(x_1,1),(1,x_2)$.
\end{example}

\goodbreak
{\it Seminormalization}

Recall from \cite[1.7]{chww-p} that a reduced monoid $A$ is {\it seminormal} if
whenever $b,c\in A$ satisfy $b^3=c^2$ there is an $a\in A$ such that $a^2=b$ and
$a^3=c$. Any normal monoid is seminormal, and $\langle x,y|xy=0\rangle$ is
seminormal but not normal. The passage from monoids to seminormal monoids (and
monoid schemes) was critical in \cite{chww-p} for understanding the behaviour of
cyclic bar constructions under the resolution of singularities of a pc monoid
scheme.

The {\it seminormalization} of a monoid $A$ is a seminormal monoid $A_\sn$,
together with an injective map $A_\red\to A_\sn$ such that every $b\in A_\sn$
has $b^n\in A_\red$ for all $n\gg0$. If it exists, it is unique up to
isomorphism, and any monoid map $A\to C$ with $C$ seminormal factors uniquely
through $A_\sn$; see \cite[1.11]{chww-p}. In particular, the seminormalization
of $A$ lies between $A$ and its normalization,  i.e.,
$\MSpec(A)_\nor\to\MSpec(A)$ factors through $\MSpec(A_\sn)$.

We shall restrict ourselves to the seminormalization of pc monoids (and monoid
schemes). By \cite[1.15]{chww-p}, if $A$ is a pc monoid, the seminormalization
of $A$  exists and is a pc monoid. When $A$ is cancellative, $A_\sn$ is easy to
construct.

\begin{example}\label{sn.cancellative}
When $A$ is cancellative, $A_\sn =\{ b\in A_0: b^n\in A\text{ for }n\gg0\}$;
this is a submonoid of $A_\nor$, and $A_\nor=(A_\sn)_\nor$. Since the
normalization of a cancellative monoid induces a homeomorphism on the
topological spaces $\MSpec$ \cite[1.6.1]{chww}, so does the seminormalization.
\end{example}

If $A$ has more than one minimal prime, then $\MSpec(A)_\nor\to\MSpec(A)$ cannot
be a bijection. However, we do have the following result.

\begin{lemma}\label{sn.homeo}
For every pc monoid $A$, $\MSpec(A_\sn)\to\MSpec(A)$ is a homeomorphism of the
underlying topological spaces.
\end{lemma}

\begin{myproof}
Write $A=C/I$ for a cancellative monoid $C$, so $\MSpec(A)$ is the closed
subspace of $\MSpec(C)$ defined by $I$. By \cite[1.14]{chww-p},
$A_\sn=C_\sn/(IC_\sn)$. Thus $\MSpec(A_\sn)$ is the closed subspace of
$\MSpec(C_\sn)$ defined by $I$. Since $\MSpec(C_\sn)\to\MSpec(C)$ is a
homeomorphism (by \ref{sn.cancellative}), the result follows.
\end{myproof}

The seminormalization of any pc monoid scheme may be constructed by glueing,
since the seminormalization of $A$ commutes with localization
\cite[1.13]{chww-p}. Thus if $X$ is a pc monoid scheme then there are canonical
maps
\[X_\nor \to X_\sn\to X_\red\to X, \]
and $X_\sn\to X$ is a homeomorphism by Lemma \ref{sn.homeo}.  We will use
$X_\sn$ to discuss the Picard group $\Pic(X)$ in Proposition \ref{Pic.Xsn}
below.

\bigskip
\section{Weil divisors}\label{sec:Weil}

Although the theory of Weil divisors is already interesting for normal monoids,
it is useful to state it for normal monoid schemes.

Let $X$ be a normal monoid scheme with generic monoid $A_0$. Corollary
\ref{c:DVmonoid} states that the stalk $\cA_x$ is a DV monoid (see after Lemma
\ref{l:wedge.normal}) for every height one point $x$ of $X$. When $X$ is
separated, a discrete valuation on $A_0$ uniquely determines a point $x$ of $X$
\cite[8.9]{chww}.

By a {\it Weil divisor} on $X$ we mean an element of the free abelian group
$\Div(X)$ generated by the height one points of $X$. We define the divisor of
$a\in A_0^\times$ to be the sum, taken over all height one points of $X$:
\[\div(a) = \sum_x v_x(a) x. \]
When $X=\MSpec(A)$ is of finite type, there are only finitely many prime ideals
in $A$, so this is a finite sum. Divisors of the form $\div(a)$ are called {\it
principal divisors}. Since $v_x(ab)=v_x(a)+v_x(b)$, the function
$\div:A_0^\times\to\Div(X)$ is a group homomorphism, and the principal divisors
form a subgroup of $\Div(X)$.

\begin{definition}
The {\it Weil divisor class group} of $X$, written as $\Cl(X)$, is the quotient
of $\Div(X)$ by the subgroup of principal divisors.
\end{definition}

\begin{lemma}\label{Cl.sequence}
If $X$ is a normal monoid scheme of finite type, there is an exact sequence
\[1 \to \cA(X)^\times \to A_0^\times \map{\div} \Div(X) \to \Cl(X) \to 0. \]
\end{lemma}

\begin{myproof}
We may suppose that $X$ is connected. It suffices to show that if $a\in
A_0^\times$ has $\div(a)=0$ then $a\in\cA(X)^\times$. This follows from Theorem
\ref{t:intersect.Ap}: when $X=\MSpec(A)$, $A$ is the intersection of the $A_x$.
\end{myproof}

\begin{example}
(Cf.\,\cite[II.6.5.2]{Hart}) Let $A$ be the submonoid of $\Z^2_*$ generated by
$x=(1,0)$, $y=(1,2)$ and $z=(1,1)$, and set $X=\MSpec(A)$. (This is the toric
monoid scheme $xy=z^2$.) Then $A$ has exactly two prime ideals of height one:
$p_1=(x,z)$ and $p_2=(y,z)$. Since $\div(x)=2p_1$ and $\div(z)=p_1+p_2$, we see
that $\Cl(X)=\Z/2$.
\end{example}

\begin{example}
If $X$ is the non-separated monoid scheme obtained by gluing together $n+1$
copies of $\A^1$ along the common (open) generic point, then $\Cl(X)=\Z^n$, as
we see from Lemma \ref{Cl.sequence}.  \end{example}

If $U$ is an open subscheme of $X$, with complement $Z$, the standard argument
\cite[II.6.5]{Hart} shows that there is a surjection $\Cl(X)\to\Cl(U)$, that it
is an isomorphism if $Z$ has codimension $\ge2$, and that if $Z$ is the closure
of a height one point $z$ then there is an exact sequence
\[\Z\; \map{z} \Cl(X) \to \Cl(U) \to 0. \]

\begin{proposition}
$\Cl(X_1\times X_2) = \Cl(X_1)\oplus\Cl(X_2).$
\end{proposition}

\begin{myproof}
By \cite[3.1]{chww}, the product monoid scheme exists, and its underlying
topological space is the product. Thus a codimension one point of $X_1\times
X_2$ is either of the form $x_1\times X_2$ or $X_1\times x_2$. Hence
$\Div(X_1\times X_2) \cong \Div(X_1)\oplus\Div(X_2)$. It follows from Lemma
\ref{l:wedge.normal} that $X_1\times X_2$ is normal, and the pointed monoid at
its generic point is the smash product of the pointed monoids $A_1$ and $A_2$ of
$X_1$ and $X_2$ at their generic points. If $a_i\in A_i$ then the principal
divisor of $a_1\wedge a_2$ is $\div(a_1)+\div(a_2)$. Thus
\[\Cl(X_1\times X_2)=
\frac{\Div(X_1)\oplus\Div(X_2)}{\div(A^{\times}_1)\oplus\div(A^{\times}_2)}
\cong \Cl(X_1)\oplus\Cl(X_2).   \qedhere \]
\end{myproof}

\begin{example}\label{normal.v.toric}
By \cite[4.5]{chww}, any connected separated normal monoid scheme $X$ decomposes
as the product of a toric monoid scheme $X_\Delta$ and $\MSpec(U_*)$ for some
finite abelian group $U$. ($U$ is the group of global units of $X$.) Since $U_*$
has no height one primes, $\Div(X)=\Div(X_\Delta)$ and the Weil class group of
$X$ is $\Cl(X_\Delta)$, the Weil class group of the associated toric monoid
scheme.

By construction \cite[4.2]{chww}, the points of $X_\Delta$ correspond to the
cones of the fan $\Delta$ and the height one points of $X_\Delta$ correspond to
the edges in the fan. Thus our Weil divisors correspond naturally to what Fulton
calls a ``$T$-Weil divisor'' on the associated toric variety $X_k$ (over a field
$k$) in \cite[3.3]{F}. Since the group completion $A_0$ is the free abelian
group $M$ associated to $\Delta$, it follows from \cite[3.4]{F} that our Weil
divisor class group $\Cl(X_\Delta)$ is isomorphic to the Weil divisor class
group $\Cl(X_k)$ of associated toric variety.
\end{example}

\newpage
\section{Invertible sheaves}\label{sec:Pic}

Let $X$ be a monoid scheme with structure sheaf $\cA$. An {\it invertible sheaf}
on $X$ is a sheaf $\cL$ of $\cA$-sets which is locally isomorphic to $\cA$ in
the Zariski topology. If $\cL_1,\cL_2$ are invertible sheaves, their smash
product is the sheafification of the presheaf $U\mapsto
\cL_1(U)\wedge_{\cA(U)}\cL_2(U)$; it is again an invertible sheaf. Similarly,
$\cL^{-1}$ is the sheafification of $U\mapsto \Hom_{\cA}(\cL(U),\cA(U))$, and
evaluation $\cL\wedge_{\cA}\cL^{-1}\map{\sim}\cA$ is an isomorphism. Thus the
set of isomorphism classes of invertible sheaves on $X$ is a group under the
smash product.

\begin{definition}\label{def:Pic}
The {\it Picard group} $\Pic(X)$ is the group of isomorphism classes of
invertible sheaves on $X$.
\end{definition}

Since a monoid $A$ has a unique maximal ideal (the non-units), an invertible
sheaf on $\MSpec(A)$ is just an $A$-set isomorphic to $A$. This proves:

\begin{lemma}\label{Pic.affine}
For every affine monoid scheme $X=\MSpec(A)$, $\Pic(X)=0$.
\end{lemma}

For any monoid $A$, the group of $A$-set automorphisms of $A$ is canonically
isomorphic to $A^\times$. Since the subsheaf $\Gamma$ of generators of an
invertible sheaf $\cL$ is a torsor for $\cA^\times$, and
$\cL=\cA\wedge_{\cA^\times}\Gamma$, this proves:

\begin{lemma}\label{Pic=H1}
$\Pic(X) \cong H^1(X,\cA^\times)$.
\end{lemma}

Recall that a morphism $f:Y\to X$ of monoid schemes is {\it affine} if
$f^{-1}(U)$ is affine for every affine open $U$ in $X$; see \cite[6.2]{chww}.

\begin{proposition}\label{affine.f_*}
If $f:Y\to X$ is an affine morphism of monoid schemes, then the direct image
$f_*$ is an exact functor from sheaves (of abelian groups) on $Y$ to sheaves on
$X$. In particular, $H^*(Y,\cL)\cong H^*(X,p_*\cL)$ for every sheaf $\cL$ on
$Y$.
\end{proposition}

\begin{myproof}
Suppose that $0\to\cL'\to\cL\to\cL''\to0$ is an exact sequence of sheaves on
$Y$. Fix an affine open  $U=\MSpec(A)$ of $X$ with closed point $x\in X$. Then
$f^{-1}(U)=\MSpec(B)$ for some monoid $B$. If $y\in Y$ is the unique closed
point of $\MSpec(B)$ the stalk sequence $0\to\cL'_y\to\cL_y\to\cL''_y\to0$ is
exact. Since this is the stalk sequence at $x$ of $0\to f_*\cL'\to f_*\cL\to
f_*\cL''\to0$, the direct image sequence is exact.
\end{myproof}

Here is an application, showing one way in which monoid schemes differ from
schemes. Let $T$ denote the free (pointed) monoid on one generator $t$, and let
$\A^1$ denote $\MSpec(T)$. Then $A\wedge T$ is the analogue of a polynomial ring
over $A$, and $X\times\A^1$ is the monoid scheme which is locally
$\MSpec(A)\times\A^1=\MSpec(A\wedge T)$; see \cite[3.1]{chww}. Thus
$p:X\times\A^1\to X$ is affine, and $f_*\cA_Y^\times=\cA_X^\times$. From
Proposition \ref{affine.f_*} we deduce

\begin{corollary}\label{Pic.HI}
For every monoid scheme $X$, $\Pic(X)\cong\Pic(X\times\A^1)$.
\end{corollary}

\section{Cartier divisors}\label{sec:Cartier}

Let $(X,\cA)$ be a cancellative monoid scheme. We write $A_0$ for the stalk of
$\cA$ at the generic point of $X$, and $\cA_0$ for the associated constant
sheaf. A {\it Cartier divisor} on $X$ is a global section of the sheaf of groups
$\cA_0^\times/\cA^\times$. On each affine open $U$, it is given by an $a_U\in
A_0^\times$ up to a unit in $\cA(U)^\times$, and we have the usual
representation as $\{(U,a_U)\}$ with $a_U/a_V$ in $\cA(U\cap V)^\times$. We
write $\Cart(X)$ for the group of Cartier divisors on $X$. The {\it principal}
Cartier divisors, i.e., those represented by some $a\in A_0^\times$, form a
subgroup of $\Cart(X)$.

\begin{proposition}\label{Cart.Pic}
Let $X$ be a cancellative monoid scheme. Then the map $D\mapsto\cL(D)$ defines
an isomorphism between the group of Cartier divisors modulo principal divisors
and $\Pic(X)$.
\end{proposition}

\begin{myproof}
Consider the short exact sequence of sheaves of abelian groups
\[1 \to\cA^\times \to \cA_0^\times \to \cA_0^\times/\cA^\times \to 1. \]
Since $\cA_0^\times$ is constant and $X$ is irreducible we have
$H^1(X,\cA_0^\times)=0$ \cite[III.2.5]{Hart}. By Lemma \ref{Pic=H1}, the
cohomology sequence becomes:
\[0\to\cA(X)^\times \to \cA_0^\times \map{\div} \Cart(X) \map{\delta}\Pic(X)\to
0. \qedhere \]
\end{myproof}

\begin{example}\label{ex:L(D)}
If $D$ is a Cartier divisor on a cancellative monoid scheme $X$, represented by
$\{(U,a_U)\}$, we define a subsheaf $\cL(D)$ of the constant sheaf $\cA_0$ by
letting its restriction to $U$ be generated by $a_U^{-1}$. This is well defined
because $a_U^{-1}$ and $a_V^{-1}$ generate the same subsheaf on $U\cap V$. The
usual argument \cite[II.6.13]{Hart} shows that $D\mapsto\cL(D)$ defines an
isomorphism from $\Cart(X)$ to the group of invertible subsheaves of
$A_0^\times$. By inspection, the map $\delta$ in \ref{Cart.Pic} sends $D$ to
$\cL(D)$.
\end{example}

\begin{lemma}
If $X$ is a normal monoid scheme of finite type, $\Pic(X)$ is a subgroup of
$\Cl(X)$.
\end{lemma}

\begin{myproof}
Every Cartier divisor $D=\{(U,a_U)\}$ determines a Weil divisor; the restriction
of $D$ to $U$ is the divisor of $a_U$. It is easy to see that this makes the
Cartier divisors into a subgroup of the Weil divisor class group $D(X)$, under
which principal Cartier divisors are identified with principal Weil divisors.
This proves the result.
\end{myproof}

\begin{theorem}
Let $X$ be a separated connected monoid scheme. If $X$ is locally factorial then
every Weil divisor is a Cartier divisor, and $\Pic(X)=\Cl(X)$.
\end{theorem}

\begin{myproof}
By Example \ref{e:normal}, $X$ is normal since factorial monoids are normal.
Thus $\Pic(X)$ is a subgroup of $\Cl(X)$, and it suffices to show that every
Weil divisor $D=\sum n_i x_i$ is a Cartier divisor. For each affine open $U$,
and each point $x_i$ in $U$, let $p_i$ be the generator of the prime ideals
associated to $x_i$; then the divisor of $a_U=\prod p_i^{n_i}$ is the
restriction of $D$ to $U$, and $D=\{(U,a_U)\}$.
\end{myproof}

\begin{lemma}\label{Pic.Pn}
For the projective space monoid scheme $\P^n$ we have
$$\Pic(\P^n)=\Cl(\P^n)=\Z.$$
\end{lemma}

\begin{remark}
This calculation of $\Pic(\P^n)$ formed the starting point of our investigation.
We learned it from Vezzani (personal communication), but it is also found in
\cite{cls} and \cite{GHS}. Related calculations are in \cite{Hutt} and
\cite{Sz}.
\end{remark}

\begin{myproof}
Since $\P^n$ is locally factorial, $\Pic(\P^n)=\Cl(\P^n)$. By definition, $\P^n$
is $\MProj$ of the free abelian monoid on $\{ x_0,...,x_n\}$, and $A_0$ is the
free abelian group with the $x_i/x_0$ as basis ($i=1,...,n$). On the other hand,
$\Div(\P^n)$ is the free abelian group on the generic points $[x_i]$ of the
$V(x_i)$. Since $\div(x_i/x_0)=[x_i]-[x_0]$, the result follows.
\end{myproof}

Before proceeding, we recall a definition from \cite{chww}.  If $X=\MSpec(A)$ is
an affine monoid scheme, define its \emph{$k$-realization} $X_k$ to be
$\Spec(k[A])$ where $k[A]$ is the $k$-realization of the monoid $A$ (see Section
\ref{krealization}).  The (affine monoid scheme) $k$-realization functor has
left adjoint $\Spec(R)\mapsto\MSpec(R,\times)$ where $(R,\times)$ is the
underlying multiplicative monoid of the $k$-algebra $R$.  If $\mathbf{MSch}$
denotes the category of monoid schemes, the adjunction
\[\Hom(\Spec(R),X_k)\cong\Hom_{\mathbf{MSch}}(\MSpec(R,\times),X)\]
defines a functor $\Spec(R)\mapsto \MSpec((R,\times),X)$ represented by $X_k$.
If $X$ is any monoid scheme and $k$ a ring, define a contravariant functor $F_X$
on the category of affine $k$-schemes to be the (Zariski) sheafification of the
presheaf
\[\Spec(R)\mapsto\Hom_{\mathbf{MSch}}(\MSpec(R,\times),X).\]
It is shown in \cite[5.2]{chww} that $F_X$ is represented by a $k$-scheme $X_k$
which is the \emph{$k$-realization} of $X$.  Note that this agrees with the 
previous definition when $X=\MSpec(A)$.

Let $\Delta$ be a fan and $X$ the toric monoid scheme associated to $\Delta$ by
\cite[4.2]{chww}.  and $X_k$ the usual toric variety associated to $\Delta$ over
some field $k$.  ($X_k$ is the $k$-realization $X_k$ of $X$.) As pointed out in
Example \ref{normal.v.toric}, our Weil divisors correspond to the $T$-Weil
divisors of the toric variety $X_k$ and $\Cl(X)\cong\Cl(X_k)$.  Moreover, our
Cartier divisors on $X$ correspond to the $T$-Cartier divisors of
\cite[3.3]{F}). Given this dictionary, the following result is established by
Fulton in \cite[3.4]{F}.

\begin{theorem}\label{thm:toric}
Let $X$ and $X_k$ denote the toric monoid scheme and toric variety (over $k$)
associated to a given fan. Then
$\Pic(X) \cong \Pic(X_k).$

Moreover, $\Pic(X)$ is free abelian if $\Delta$ contains a cone of maximal
dimension.
\end{theorem}

\section{$\Pic$ of pc monoid schemes}

In this section, we derive some results about the Picard group of pc monoid
schemes.
When $X$ is a pc monoid scheme, we can form the reduced monoid scheme
$X_\red=(X,\cA_\red)$ using Remark \ref{r:nil(A)}: the stalk of $\cA_\red$ at
$x$ is $\cA_x/\nil(\cA_x)$.  Since $\cA^\times=\cA_\red^\times$, the map
$X_\red\to X$ induces an isomorphism $\Pic(X)\cong\Pic(X_\red)$.

We will use the constructions of normalization and seminormalization given in
Section \ref{sec:normalize}.

\begin{proposition}\label{Pic.Xsn}
If $X$ is a pc monoid scheme, the canonical map $X_\sn\to X$ induces an
isomorphism $\Pic(X)\cong\Pic(X_\sn)$.
\end{proposition}

\begin{myproof}
Since $X_\red$ and $X$ have the same underlying space, it suffices by Lemma
\ref{Pic=H1} to assume that $X$ is reduced and show that the inclusion
$\cA^\times\to\cA_\sn^\times$ is an isomorphism. It suffices to work stalkwise,
so we are reduced to showing that if $A$ is reduced then $A^\times\to
A_\sn^\times$ is an isomorphism. If $b\in A_\sn^\times$ then both $b^n$ and
$(1/b)^n$ are in $A$ for large $n$, and hence both $b=b^{n+1}b^{-n}$ and
$b^{-1}=b^n(1/b)^{1+n}$ are in $A$, so $b\in A^\times$.
\end{myproof}

\begin{lemma}\label{Pic.Xnor}
Let $X$ be a cancellative seminormal monoid scheme and $p:X_\nor\to X$ its
normalization. If $\cH$ denotes the sheaf $p_*(\cA_\nor^\times)/\cA^\times$ on
$X$, there is an exact sequence
\[1\to\cA(X)^\times\to\cA_\nor(X_\nor)^\times \to H^0(X,\cH) \to \Pic(X)
\map{p^*} \Pic(X_\nor) \to H^1(X,\cH). \]
\end{lemma}

\begin{myproof}
At each point $x\in X$, the stalk $A=\cA_x$ is a submonoid of its normalization
$A_\nor=p_*(\cA_\nor)_x$ (by Lemma \ref{l:local.normal}) and we have an exact
sequence of sheaves on $X$:
\[1 \to \cA^\times \to p_*(A_\nor^\times) \to \cH\to 1. \]
Since $p$ is affine, Proposition \ref{affine.f_*} implies that
$\cA_\nor(X_\nor)^\times=H^0(X,p_*\cA_\nor^\times)$ and
$\Pic(X_\nor)=H^1(X,p_*A_\nor^\times)$, and the associated cohomology sequence
is the displayed sequence.
\end{myproof}

Here are two examples showing that $\Pic(X)\to\Pic(X_\nor)$ need not be an
isomorphism when $X$ is seminormal and cancellative.

\begin{example}
Let $A_+$ (resp., $A_-$) be the submonoid of the free monoid $B=\langle
x,y\rangle$ generated by $\{ x,y^2,xy\}$ (resp., $\{ x,y^{-2},xy^{-1}\}$). These
are seminormal but not normal. If $X$ is the monoid scheme obtained by gluing
the $U_\pm=\MSpec(A_\pm)$ together along $\MSpec(\langle x,y^2,y^{-2}\rangle)$
then it is easy to see that $\Pic(X)=\Z$, with a generator represented by
$(U_+,y^2)$ and $(U_-,1)$.  The normalization $X_\nor$ is the toric monoid
scheme $\A^1\times\P^1$, and $\Pic(X_\nor)\cong\Z$, with a generator represented
by $(U_+,y)$ and $(U_-,1)$. Thus $\Pic(X)\to\Pic(X_\nor)$ is an injection with
cokernel $\Z/2$.
\end{example}

\begin{example}
Let $U$ be an abelian group and $A_x$ the submonoid of $B=U_*\wedge\langle
x\rangle$ consisting of $0,1$ and all terms $ux^n$ with $u\in U$ and $n>0$.
Note that $U_*\land \langle x\rangle$ my be written $U\ot\F_1[x]$ or simply
$U[x]$ (see Section \ref{tensoralgebras}).  Then $A_x$ is seminormal and $B$ is
its normalization.  Let $X$ be obtained by gluing $\MSpec(A_x)$ and
$\MSpec(A_{1/x})$ together along their common generic point,
$\MSpec(U_*\wedge\langle x,1/x\rangle)$.  The normalization of $X$ is
$X_\nor=\MSpec(U_*)\times\P^1$, and $\Pic(X_\nor)=\Z$ by Example
\ref{normal.v.toric} and Lemma \ref{Pic.Pn}.  Because
$p_*(\cA_\nor^\times)/\cA^\times$ is a skyscraper sheaf with stalk $U$ at the
two closed points, we see from Lemma \ref{Pic.Xnor} that $\Pic(X)=\Z\times U$.
Thus $\Pic(X)\to\Pic(X_\nor)$ is a surjection with kernel $U$.
\end{example}

Finally, we consider the case when $X$ is reduced pc monoid scheme which is not
cancellative. We may suppose that $X$ is of finite type, so that the stalk at a
closed point is an affine open $\MSpec(A)$ with minimal points
$\frakp_1,...,\frakp_r$, $r>1$. Then the closure $X'$ of the point $\frakp_1$ is
a cancellative seminormal monoid scheme. Let $X''$ denote the closure of the
remaining minimal points of $X$, and set $X'''=X'\cap X''$. Then we have the
exact Mayer-Vietoris sequence of sheaves on $X$:
\[0 \to \cA_X \to p'_*\cA_{X'}\times p''_*\cA_{X''} \to p'''_*\cA_{X'''} \to 1.
\]
Because the immersions are affine, Proposition \ref{affine.f_*} yields the exact
sequence
\begin{equation}\label{eq:MV}
1\to\cA(X)^\times \to \cA_{X'}^\times \times \cA_{X''}^\times
\to \cA_{X'''}^\times \to
\Pic(X)\to \Pic(X')\times\Pic(X'')\to \Pic(X''').
\end{equation}
The Picard group may then be determined by induction on $r$ and $\dim(X)$.

\begin{example}
If $X$ is obtained by gluing together $X_1,..., X_n$ at a common generic point,
then \eqref{eq:MV} yields $\Pic(X)=\oplus\Pic(X_i)$.
\end{example}

\end{document}